\documentstyle[12pt]{article}
\begin{document}
\pagestyle{empty}
\renewcommand{\thefootnote}{\fnsymbol{footnote}}
\def\lsim{\raise0.3ex\hbox{$<$\kern-0.75em\raise-1.1ex\hbox{$\sim$}}}
\def\gsim{\raise0.3ex\hbox{$>$\kern-0.75em\raise-1.1ex\hbox{$\sim$}}}
\def\noi{\noindent}
\def\sq{\hbox {\rlap{$\sqcap$}$\sqcup$}}
\def\R{ {\rm R \kern -.31cm I \kern .15cm}}
\def\C{ {\rm C \kern -.15cm \vrule width.5pt \kern .12cm}}
\def\Z{ {\rm Z \kern -.27cm \angle \kern .02cm}}
\def\N{ {\rm N \kern -.26cm \vrule width.4pt \kern .10cm}}
\def\1{{\rm 1\mskip-4.5mu l} }
\def\lsim{\raise0.3ex\hbox{$<$\kern-0.75em\raise-1.1ex\hbox{$\sim$}}}
\def\gsim{\raise0.3ex\hbox{$>$\kern-0.75em\raise-1.1ex\hbox{$\sim$}}}
\vbox to 2 truecm {}
\centerline{\Large \bf Long Range Scattering and Modified}
\vskip 3 truemm  
\centerline{\Large \bf Wave Operators for some Hartree Type Equations III}
\vskip 3 truemm  
\centerline{\Large \bf Gevrey spaces and low dimensions\footnote{Work supported in
part by NATO Collaborative Linkage Grant 976047.}}
\vskip 1 truecm
\centerline{\bf J. Ginibre}
\centerline{Laboratoire de Physique Th\'eorique\footnote{Unit\'e Mixte de Recherche (CNRS) UMR
8627.}}  \centerline{Universit\'e de Paris XI, B\^atiment 210,
F-91405 Orsay Cedex, France}

\vskip 5 truemm
\centerline{\bf G. Velo}
\centerline{Dipartimento di Fisica, Universit\`a di
Bologna}
\centerline{ and INFN, Sezione di Bologna, Italy }

\vskip 1 truecm 
\begin{abstract}
We study the theory of scattering for a class of Hartree type equations with long range
interactions in arbitrary space dimension $n \geq 1$, including the case of Hartree equations with
time dependent potential $V(t, x) = \kappa \ t^{\mu - \gamma} |x|^{-\mu}$ with $0 < \gamma \leq 1$
and $0 < \mu < n$. This includes the case of potential $V(x) = \kappa \ |x|^{- \gamma}$, and can be
extended to the limiting case of nonlinear Schr\"odinger equations with cubic nonlinearity $\kappa
\  t^{n - \gamma} |u|^2u$. Using Gevrey spaces of asymptotic states and solutions, we prove the
existence of modified local wave operators at infinity with no size restriction on the data and
we determine the asymptotic behaviour in time of solutions in the range of the wave operators,
thereby extending the results of previous papers which covered the range $0 < \gamma \leq 1$ but
only $0 < \mu \leq n-2$, and were therefore restricted to space dimension $n \geq 3$.
   \end{abstract}
\vskip 3 truecm
\noi AMS Classification : Primary 35P25. Secondary 35B40, 35Q40, 81U99.  \par
\noi Key words : Long range scattering, modified wave operators, Hartree equation, Gevrey spaces. 
\vskip 1 truecm

\noindent LPT Orsay 00-08 \par
\noindent January 2000 \par

  \newpage
\pagestyle{plain}
\baselineskip=24 pt

\section{Introduction}
\hspace*{\parindent} This is the third paper where we study the theory of scattering and more
precisely the existence of modified wave operators for a class of long range Hartree type equations
$$i \partial_t u  + {1 \over 2} \Delta u = \widetilde{g} (|u|^2)u \eqno(1.1)$$
\noi where $u$ is a complex function defined in space time ${I\hskip-1truemm R}^{n+1}$, $\Delta$
is the Laplacian in ${I\hskip-1truemm R}^n$, and
$$\widetilde{g} (|u|^2) = \kappa \ t^{\mu - \gamma} \  |\nabla |^{\mu - n} \ |u|^2 \eqno(1.2)$$
\noi with $|\nabla | = (- \Delta)^{1/2}$, $\kappa \in {I\hskip-1truemm R}$, $0 < \gamma \leq 1$ and
$0 < \mu \leq n$. For $\mu < n$, the operator $|\nabla|^{\mu - n}$ can be represented by the
convolution in ${I\hskip-1truemm R}^n$ 
$$|\nabla|^{\mu - n} \ f = C_{n, \mu} \ |x|^{- \mu} * f \eqno(1.3)$$
\noi so that (1.2) is a Hartree type interaction with potential $V(x) = C |x|^{-
\mu}$. The more standard Hartree equation corresponds to the case $\gamma = \mu$. In that case,
the nonlinearity $\widetilde{g}(|u|^2)$ becomes
$$\widetilde{g}(|u|^2) = V * |u|^2 = \kappa |x|^{-\gamma} * |u|^2 \eqno(1.4)$$
\noi with a suitable redefinition of $\kappa$. \par

A large amount of work has been devoted to the theory of scattering for the Hartree equation
(1.1) with nonlinearity (1.4) as well as with similar nonlinearities with more general potentials.
As in the case of the linear Schr\"odinger equation, one must
distinguish the short range case, corresponding to $\gamma > 1$, from the long range case
corresponding to $\gamma \leq 1$. In the short range case, it is known that the (ordinary) wave
operators exist in suitable function spaces for $\gamma > 1$ \cite{14r}. Furthermore for repulsive
interactions, namely for $\kappa \geq 0$, it is known that all solutions in suitable spaces admit
asymptotic states in $L^2$ for $\gamma > 1$, and that asymptotic completeness holds for $\gamma >
4/3$ \cite{12r}. In the long range case $\gamma \leq 1$, the ordinary wave operators are known not
to exist in any reasonable sense \cite{12r}, and should be replaced by
modified wave operators including a suitable phase in their definition, as is the case for the
linear Schr\"odinger equation. A well developed theory of long range scattering exists for the
latter. See for instance \cite{1r} for a recent treatment and for an extensive
bibliography. In contrast with that situation, only partial results are available for the
Hartree equation. For small solutions (or equivalently small asymptotic states) the existence of
modified wave operators has been proved in the critical case $\gamma = 1$
\cite{2r}. On the other hand, it has been shown, first in the critical case $\gamma = 1$
and then in the whole range $0 < \gamma \leq 1$ \cite{6r,7r,9r,10r} that the global
solutions of the Hartree equation (1.1) with nonlinearity (1.4) and with small initial data
exhibit an asymptotic behaviour as $t \to \pm \infty$ of the expected scattering type
characterized by scattering states $u_{\pm}$ and including suitable phase factors that are
typical of long range scattering. \par

In the previous two papers of this series \cite{4r} \cite {5r} (hereafter referred to as I and II)
we proved the existence of modified wave operators in suitable spaces for the equation (1.1)
with nonlinearity (1.2), and we gave a description of the asymptotic behaviour in time of
solutions in the ranges of those operators, with no size restriction on the data, first for
$1/2 < \gamma < 1$ in I and then in the whole range $0 < \gamma \leq 1$ in II. The method is
an extension of the energy method used in \cite{6r,9r,10r}, and uses in particular an
auxiliary system of equations introduced in \cite{9r} to study the asymptotic behaviour of
small solutions. The spaces of initial data, namely in the present case of asymptotic
states, are Sobolev spaces of finite order. However there occurs a loss of derivatives in
the auxiliary system, which has to be compensated for by the smoothing effect of the
operator $|\nabla |^{\mu - n}$ in (1.2). This is done in the framework of Sobolev spaces at
the expense of assuming $\mu \leq n - 2$, which in particular restrict the space dimension
to $n \geq 3$. In the present paper, we overcome that difficulty by treating the problem in
Gevrey spaces \cite{13r}, following and extending the method used in \cite{7r} to treat the
case of small solutions. This makes it possible to cover the whole range $0 < \mu \leq n$,
and in particular the case of dimensions 1 and 2 and the case of cubic nonlinear
Schr\"odinger (NLS) equations (with time dependent nonlinearity). More precisely we use
Gevrey classes $G^{1/\nu}$ of order $1/\nu$ with $0 < \nu \leq 1$, and the method applies
under the condition $\mu \leq n - 2 + 2 \nu$. In particular for cubic NLS equations we need
$\nu = 1$, namely spaces of analytic functions. The previous restriction on $\mu$ and $\nu$
can still be weakened and has been weakened in \cite{7r} at the expense of introducing
parabolic terms in the auxiliary system of equations. However those terms introduce a
priviledged orientation of time, which is inconvenient for the study of scattering theory,
where we like to go back and forth from finite to infinite time, and we shall not make use
of that extension here. The origin of the derivative loss and the mechanism by which that
loss is overcome in Gevrey spaces, which is the same as in \cite{7r}, will be described in
Section 3 after sufficient technical material has been introduced, namely after Lemma 3.4.
\par

The construction of the modified wave operators is in its principle the same as in II and will
be recalled in Section 2 below, which is mostly a summary of Section 2 of II. It involves the
study of the same auxiliary system of equations as in II for an amplitude $w$ and a phase
$\varphi$ which replace the original function $u$, and the definition of the same modified
asymptotic dynamics for that system as in II. \par

We now give a brief outline of the contents of this paper. A more technical description will be
given at the end of Section 2. In Section 3 we define the relevant Gevrey spaces and derive the
basic estimates in those spaces that are needed to study the auxiliary system. In Section 4 we
prove the existence of the large time dynamics associated with that system and some preliminary
asymptotic properties of that dynamics. In Section 5 we study the asymptotic dynamics and we
prove the existence of asymptotic states for the previously constructed solutions of the
auxiliary system. In Section 6 we construct the local wave operators at infinity for the
auxiliary system by solving the Cauchy problem for that system with infinite initial time. We
then come back from the auxiliary system to the original equation (1.1) for $u$ and construct
the (local) wave operators (at infinity) for $u$ in Section 7, where the main result is stated
as Proposition 7.5. \par

We have tried to make this paper self-contained and at the same time to keep duplication with I
and II to a minimum. Duplication occurs in Section 2, as already mentioned, and in part of
Section 7. The more technical sections 3 to 6 follow the same pattern as in II, but there is
almost no duplication because the functional framework is significantly different. \par

We conclude this section by giving some general notation which will be used freely throughout this
paper. We shall work mostly in Fourier space. We denote by $*$ the convolution in ${I\hskip-1truemm
R}^n$, by $F$ the Fourier transform, and by $\widehat{u} = Fu$ the Fourier transform of $u$. We
denote by $\parallel \cdot \parallel_r$ the norm in $L^r \equiv L^r ({I\hskip-1truemm R}^n)$ and by
$<\cdot , \cdot >$ the scalar product in $L^2$. For any interval $I$ and
any Banach space $X$, we denote by ${\cal C}(I, X)$ the space of strongly continuous functions
from $I$ to $X$, by $L^{\infty}(I, X)$ (resp. $L_{loc}^{\infty}(I, X))$ the space of
measurable essentially bounded (resp. locally essentially bounded) functions from $I$ to $X$, and
by $L^2(I,X)$ (resp. $L^2_{loc}(I,X)$, resp. $L^2_{\rho}(I, X))$ the space of measurable functions
$u$ from $I$ to $X$ such that $\parallel u(\cdot );X \parallel$ belongs to $L^2(I)$ (resp.
$L^2_{loc}(I)$, resp. $L_{\rho}^2(I))$, where $L_{\rho}^2(I)$ is the weighted space $L^2(I, \rho
(t) dt)$ for some positive function $\rho$. For
real numbers $a$ and $b$, we use the notation $a \vee b = {\rm Max} (a, b)$ and $a \wedge b = {\rm
Min}(a,b)$. In the estimates of solutions of the relevant
equations, we shall use the letter $C$ to denote constants, possibly different from an
estimate to the next, depending on various pa\-ra\-me\-ters such as $\gamma$, but not on the
solutions themselves or on their initial data. Those constants will be bounded in $\gamma$
for $\gamma$ away from zero. We shall use the notation $A(a_1,a_2,\cdots)$ for es\-ti\-ma\-ting
functions, also possibly different from an estimate to the next, depending in addition on
suitable norms $a_1, a_2, \cdots$ of the solutions or of their initial data. If (p.q) is a double
inequality, we denote by (p.qa) and (p.qb) the first and second inequality in (p.q). Finally,
Item (p.q) of I or II will be referred to as Item (I.p.q.) or (II.p.q). Additional notation will be
given when needed. \par

In all this paper, we assume that $0 < \mu \leq n$ and $0 < \gamma \leq 1$.

\section{Heuristics}
\hspace*{\parindent} In this section, we describe in heuristic terms the construction of the modified
wave ope\-ra\-tors for the equation (1.1). That construction is the same as that performed in II,
and this section is mostly a summary of Section II.2, which we include in order to make this paper
self-contained. \par

The problem that we address is that of classifying the possible asymptotic behaviours of
the solutions of (1.1) by relating them to a set of model functions ${\cal V} = \{v = v(u_+)\}$
parametrized by some data $u_+$ and with suitably chosen and preferably simple asymptotic
behaviour in time. For each $v \in {\cal V}$, one tries to construct a solution $u$ of (1.1) such
that $u(t)$ behaves as $v(t)$ when $t \to \infty$ in a suitable sense. The map $\Omega : u_+ \to
u$ thereby obtained classifies the asymptotic behaviours of solutions of (1.1) and is a
preliminary version of the wave operator for positive time. A similar question can be asked for
$t \to - \infty$. From now on we restrict our attention to positive time. \par

In the short range case corresponding to $\gamma > 1$ in (1.2), the previous scheme can be
implemented by taking for ${\cal V}$ the set ${\cal V} = \{ v = U(t) u_+\}$ of solutions of the
equation
$$i \partial_t v + {1\over 2} \Delta v = 0 \quad , \eqno(2.1)$$

\noi with $U(t)$ being the unitary group
$$U(t) = \exp \left ( i(t/2)\Delta \right ) \quad . \eqno(2.2)$$

\noi The initial data $u_+$ for $v$ is called the asymptotic state for $u$. \par

In the long range case corresponding to $\gamma \leq 1$ in (1.2), the previous set is known
to be inadequate and has to be replaced by a better set of model functions obtained by
modifying the previous ones by a suitable phase. The modification that we use requires
additional structure of $U(t)$. In fact $U(t)$ can be written as
$$U(t) = M(t) \ D(t) \ F \ M(t) \eqno(2.3)$$

\noi where $M(t)$ is the operator of multiplication by the function 
$$M(t) = \exp \left ( i  x^2/2t \right ) \quad , \eqno(2.4)$$

\noi and $D(t)$ is the dilation operator defined by 
$$\left ( D(t) \ f \right ) (x) = (it)^{-n/2} \ f(x/t) \quad . \eqno(2.5)$$

\noi Let now $\varphi^{(0)} = \varphi^{(0)}(x,t)$ be a real function of space time and let
$z^{(0)}(x, t) = \exp (-i \varphi^{(0)}(x,t))$. We replace $v(t) = U(t) u_+$ by the modified free
evolution \cite{16r} \cite{17r}
$$v(t) = M(t) \ D(t) \ z^{(0)}(t) \ w_+   \eqno(2.6)$$

\noi where $w_+ = Fu_+$. In order to allow for easy comparison of $u$ with $v$, it is then
convenient to represent $u$ in terms of a phase factor $z(t) = \exp (-i \varphi (t))$ and of an
amplitude $w(t)$ in such a way that asymptotically $\varphi (t)$ behaves as $\varphi^{(0)}(t)$ and
$w(t)$ tends to $w_+$. This is done by writing $u$ in the form \cite{8r} \cite{9r}
$$u(t) = M(t) \ D(t) \ z(t) \ w(t) \equiv \left ( \Lambda (w, \varphi ) \right ) (t) \quad .
\eqno(2.7)$$

The construction of the wave operators for $u$ proceeds by first constructing the wave
o\-pe\-ra\-tors for the pair $(w , \varphi )$ and then recovering the wave operators for
$u$ therefrom by the use of (2.7). The evolution equation for $(w, \varphi )$ is obtained by substituting
(2.7) into the equation (1.1). One obtains the equation
$$\left ( i \partial_t + (2t^2)^{-1} \Delta - D^*\widetilde{g}D \right ) z w = 0 \eqno(2.8)$$

\noi for $zw$, with
$$\widetilde{g} \equiv \widetilde{g} \left ( |u|^2 \right ) = \widetilde{g} \left ( |Dw|^2 \right )
\quad , \eqno(2.9)$$

\noi or equivalently, by expanding the derivatives in (2.8),
$$\left \{ i \partial_t + (2t^2)^{-1} \Delta - i (2t^2)^{-1} \left ( 2 \nabla \varphi \cdot
\nabla + (\Delta \varphi ) \right ) \right \} w$$
$$+ \left \{ \partial_t \varphi - (2t^2)^{-1} \ |\nabla \varphi|^2 - D^*\widetilde{g}D \right \} w
= 0 \quad . \eqno(2.10)$$

We are now in the situation of a gauge theory. The equation (2.8) or (2.10) is invariant under the
gauge transformation $(w, \varphi ) \to (w \exp (i \omega ), \varphi + \omega )$, where $\omega$ is
an arbitrary function of space time, and the original gauge invariant equation is not sufficient
to provide evolution equations for the two gauge dependent quantities $w$ and $\varphi$. At this
point we arbitrarily add the Hamilton-Jacobi equation as a gauge condition. This yields a system
of evolution equations for $(w, \varphi )$, namely

$$\hskip 2.5 truecm \left \{ \begin{array}{ll} \partial_t w = i(2t^2)^{-1} \Delta w + (2t^2)^{-1}
\left ( 2 \nabla \varphi \cdot \nabla + (\Delta \varphi ) \right ) w &\hskip 4 truecm (2.11) \\
& \\
\partial_t \varphi = (2t^2)^{-1} \ |\nabla \varphi |^2 + t^{- \gamma} \ g_0 (w, w)
&\hskip 4 truecm (2.12) \end{array} \right .$$

\noi where we have defined
$$ g_0(w_1, w_2) = \kappa \ {\rm Re} \ |\nabla|^{\mu - n} \ w_1 \ \bar{w}_2 \eqno(2.13)$$

\noi and rewritten the nonlinear interaction term in (2.10) as
$$D^*\widetilde{g} \left ( |Dw|^2 \right ) D = t^{-\gamma} \ g_0(w, w) \quad . $$

\noi The gauge freedom in (2.11)-(2.12) is now reduced to that given by an arbitrary function
of space only. It will be shown in Section 4 that the Cauchy problem for
the system (2.11)-(2.12) is locally well-posed in a neighborhood of infinity in time. The
solutions thereby obtained behave asymptotically as $w(t) = O(1)$ and $\varphi (t) \cong O(t^{1-
\gamma})$ as $t \to \infty$, a behaviour that is immediately seen to be compatible with
(2.11)-(2.12). \par

We next study the asymptotic behaviour of the solutions of the auxiliary system
(2.11)-(2.12) in more detail and try to construct wave operators for that system. For that
purpose, we need to choose a set of model functions playing the role of $v$, in the spirit
of (2.6). We proceed as  follows. Let $p \geq 0$ be an integer. We write
$$\hskip 4 truecm \left \{ \begin{array}{ll} w = \displaystyle{\sum\limits_{0 \leq m \leq p}} w_m +
q_{p+1} \equiv W_p + q_{p+1} &\hskip 5 truecm (2.14) \\ & \\
\varphi = \displaystyle{\sum\limits_{0 \leq m \leq p}} \varphi_m + \psi_{p+1} \equiv \phi_p +
\psi_{p+1} &\hskip 5 truecm (2.15) \end{array} \right .$$

\noi with the understanding that asymptotically in $t$
$$w_m(t) = O \left ( t^{-m\gamma} \right ) \quad , \quad q_{p+1}(t) = o \left ( t^{-p\gamma} \right
) \quad , \eqno(2.16)$$
$$\varphi_m (t) = O \left ( t^{1 - (m+1)\gamma} \right ) \quad , \quad \psi_{p+1} (t) = o \left (
t^{1 - (p+1)\gamma} \right ) \quad . \eqno(2.17)$$

\noi Substituting (2.14)-(2.15) into (2.11)-(2.12) and identifying the various powers of
$t^{-\gamma}$ yields the following system of equations for $(w_m, \varphi_m)$~: 
$$\hskip 1 truecm\left \{ \begin{array}{ll}\partial_t \ w_{m+1} = \left ( 2t^2 \right )^{-1}
\displaystyle{\sum\limits_{0 \leq j \leq m}} \left ( 2 \nabla \varphi_j \cdot \nabla + (\Delta
\varphi_j ) \right ) w_{m-j} &\hskip 2 truecm(2.18) \\  & \\ \partial_t \ \varphi_{m+1} = \left (
2t^2 \right )^{-1} \displaystyle{\sum\limits_{0 \leq j \leq m}} \nabla \varphi_j \cdot \nabla
\varphi_{m-j} + t^{- \gamma} \displaystyle{\sum\limits_{0 \leq j \leq m+1}} g_0 \left ( w_j ,
w_{m+1-j} \right ) &\hskip 2 truecm (2.19) \end{array}\right . $$

\noi for $m + 1 \geq 0$. Here it is understood that $w_j = 0$ and $\varphi_j = 0$ for $j < 0$. We
supplement that system with the initial conditions $$\hskip 2.5 truecm \left \{ \begin{array}{ll}
w_0(\infty ) = w_+ \qquad , \quad w_m(\infty ) = 0 \quad \hbox{for} \ m \geq 1 &\hskip 4.5 truecm
(2.20) \\ & \\ \varphi_m (1) = 0 \qquad \hbox{for} \  0 \leq m \leq p \quad . &\hskip 4.5 truecm
(2.21) \end{array} \right . $$

\noi The system (2.18)-(2.19) with the initial conditions (2.20)-(2.21) can be solved by
successive integrations~: knowing $(w_j, \varphi_j)$ for $0 \leq j \leq m$, one constructs
successively $w_{m+1}$ by integrating (2.18) between $t$ and $\infty$, and then $\varphi_{m+1}$
by integrating (2.19) between 1 and $t$. \par

If $(p+1)\gamma < 1$, that method of resolution reproduces the asymptotic behaviour in time (2.16)
(2.17) which was used in the first place to provide a heuristic derivation of the system
(2.18)-(2.19). For sufficiently large $p$, $\phi_p$ is a sufficiently good approximation for
$\varphi$ to ensure that $\psi_{p+1}$ has a limit as $t \to \infty$. In fact by comparing the
system (2.18)-(2.19) with (2.11)-(2.12), one finds that $\partial_t \ \psi_{p+1}$ is
essentially of the same order in $t$ as $\partial_t \ \varphi_{p+1}$, namely $\partial_t \
\psi_{p+1} \cong O(t^{-(p+2)\gamma})$, which is integrable at infinity for $(p + 2)\gamma >
1$. In this way every solution $(w, \varphi )$ of the system (2.11)-(2.12) as obtained
previously has asymptotic states consisting of $w_+ = \lim\limits_{t \to \infty} w(t)$ and
$\psi_+ = \lim\limits_{t \to \infty} \psi_{p+1}(t)$. \par

Conversely, under the condition $(p+2)\gamma > 1$, we shall be able to solve the system
(2.11)-(2.12) by looking for solutions in the form (2.14)-(2.15) with the additional initial
condition $\psi_{p+1} (\infty ) = \psi_+$, thereby getting a solution which is asymptotic to
$(W_p, \phi_p + \psi_+)$ with  $$w - W_p \cong O\left (t^{-(p+1)\gamma} \right ) \quad ,
\quad \varphi - \phi_p - \psi_+ \cong O     \left ( t^{1-(p+2)\gamma} \right ) \quad .
\eqno(2.22)$$

This allows to define a map $\Omega_0 : (w_+, \psi_+) \to (w, \varphi )$ which is essentially
the wave operator for $(w, \varphi )$. \par

We next discuss the gauge covariance properties of $\Omega_0$. Two solutions $(w, \varphi )$ and
$(w' , \varphi ')$ of the system (2.11)-(2.12) will be said to be gauge equivalent if they
give rise to the same $u$ through (2.7), namely if $w \exp (- i \varphi ) = w' \exp (- i
\varphi ')$. If $(w, \varphi )$ and $(w' , \varphi ')$ are two gauge equivalent solutions,
one can show easily that the difference $\varphi_- = \varphi ' - \varphi$ has a limit
$\omega$ when $t \to \infty$ and that $w'_+ = w_+ \exp (i \omega )$. Under that condition,
it turns out that the phases $\{\varphi_j\}$ and $\phi_p$ (but not the amplitudes) obtained
by solving (2.18)-(2.19) are gauge invariant, namely $\varphi_m = \varphi '_m$ for $0 \leq m
\leq p$ and therefore $\phi_p = \phi '_p$, so that $\psi '_+ = \psi_+ + \omega$. It is then
natural to define gauge equivalence of asymptotic states $(w_+, \psi_+)$ and $(w'_+ , \psi
'_+)$ by the condition $w_+ \exp (- i \psi_+) = w'_+ \exp (-i \psi '_+)$ and the previous
result can be rephrased as the statement that gauge equivalent solutions of (2.11)-(2.12) in
${\cal R} (\Omega_0)$ have gauge equivalent asymptotic states. Conversely, we shall show
that gauge equivalent asymptotic states have gauge equivalent images under $\Omega_0$. Here
however we meet with a technical problem coming from the construction of $\Omega_0$ itself.
For given $(w_+ , \psi_+)$ we construct $(w, \varphi )$ in practice as follows. We take a
(large) finite time $t_0$ and we define a solution $(w_{t_0}, \varphi_{t_0})$ of the system
(2.11)-(2.12) by imposing a suitable initial condition at $t_0$, depending on $(w_+,
\psi_+)$, and solving the Cauchy problem with finite initial time. We then let $t_0$ tend to
infinity and obtain $(w, \varphi )$ as the limit of $(w_{t_0}, \varphi_{t_0})$. The simplest
way to prove the gauge equivalence of two solutions $(w, \varphi )$ and $(w', \varphi ')$
obtained in this way from gauge equivalent $(w_+, \psi_+ )$ and $(w'_+, \psi '_+)$ consists
in using an initial condition at $t_0$ which already ensures that $(w_{t_0}, \varphi_{t_0})$
and $(w'_{t_0}, \varphi '_{t_0})$ are gauge equivalent. However the natural choice
$(w_{t_0}(t_0), \varphi_{t_0} (t_0)) = (W_p(t_0), \phi_p (t_0) + \psi_+)$ does not satisfy
that requirement as soon as $p \geq 1$ because $\phi_p(t_0)$ is gauge invariant while $W_p
(t_0) \exp (- \psi_+ )$ is not. In order to overcome that difficulty, we introduce a new
amplitude $V$ and a new phase $\chi$ defined by solving the transport equations
$$\hskip 4 truecm \left \{ \begin{array}{ll} \partial_t V = (2t^2)^{-1} \left ( 2 \nabla \phi_{p-1}
\cdot \nabla + \left ( \Delta \phi_{p-1} \right ) \right ) V &\hskip 3.8 truecm(2.23) \\ & \\
\partial_t \chi = t^{-2} \ \nabla \phi_{p-1} \cdot \nabla \chi &\hskip 3.8 truecm(2.24) \end{array}
\right . $$

\noi with initial condition
$$V(\infty ) = w_+ \qquad , \qquad \chi (\infty ) = \psi_+ \quad . \eqno(2.25)$$

\noi It follows from (2.23) (2.24) that $V \exp (- i \chi )$ satisfies the same transport
equation as $V$, now with gauge invariant initial condition $(V \exp (- i \chi ))(\infty ) =
w_+ \exp (-i \psi_+)$, and is therefore gauge invariant. Furthermore, $(V, \chi )$ is a
sufficiently good approximation of $(W_p , \psi_+)$ in the sense that
$$V(t) - W_p(t) \cong O \left (t^{-(p+1)\gamma} \right ) \quad , \quad \chi (t) - \psi_+ \cong
O(t^{-\gamma}) \quad . \eqno(2.26)$$

\noi One then takes $(w_{t_0} (t_0), \varphi_{t_0} (t_0)) = (V(t_0), \phi_p (t_0) + \chi
(t_0))$ as an initial condition at time $t_0$, thereby ensuring that $(w_{t_0},
\varphi_{t_0})$ and $(w'_{t_0}, \varphi '_{t_0})$ are gauge equivalent. That equivalence is
easily seen to be preserved in the limit $t_0 \to \infty$. Furthermore, the estimates (2.26)
ensure that the asymptotic properties (2.22) are preserved by the modified construction. As a
consequence of the previous discussion, the map $\Omega_0$ is gauge covariant, namely induces
an injective map of gauge equivalence classes of asymptotic states $(w_+, \psi_+)$ to gauge
equivalence classes of solutions $(w, \varphi )$ of the system (2.11)-(2.12). \par

The wave operator for $u$ is obtained from $\Omega_0$ just defined and from $\Lambda$ defined by
(2.7). From the previous discussion it follows that the map $\Lambda \circ \Omega_0 : (w_+, \psi_+ )
\to u$ is injective from gauge equivalence classes of asymptotic states $(w_+ , \psi_+)$ to
solutions of (1.1). In order to define a wave operator for $u$ involving only the asymptotic state
$u_+$ but not an arbitrary phase $\psi_+$, we choose a representative in each equivalence class
$(w_+, \psi_+)$, namely we define the wave operator for $u$ as the map $\Omega : u_+ \to u =
(\Lambda \circ \Omega_0) (Fu_+, 0)$. Since each equivalence class of asymptotic states contains at
most one element with $\psi_+ = 0$, the map $\Omega$ is again injective. \par

The previous heuristic discussion was based in part on a number of time decay estimates in terms
of negative powers of $t$. In practice however two complications occur, namely (i) for integer
$\gamma^{-1}$, some of the estimates involve logarithmic factors in time, and (ii) the use of
Gevrey spaces requires that of norms defined by integrals over time involving a convergence factor
which eventually produces a small loss in the time decays. Both difficulties are handled by
introducing suitable estimating functions of time, some of which are defined by integral
representations and generalize in a natural way a similar family of functions defined in II.\par

In the same way as in I, the system (2.11)-(2.12) can be rewritten as a system of equations
for $w$ and for $s = \nabla \varphi$, from which $\varphi$ can then  be recovered by (2.12),
thereby leading to a slightly more general theory since the system for $(w, s)$ can be
studied without even assuming that $s$ is a gradient. For simplicity, and in the same way as
in II, we shall not follow that track. However, we shall use systematically the notation $s
= \nabla \varphi$, and for the purposes of estimation, we shall supplement the system
(2.11)-(2.12) with the equation satisfied by $s$, which is simply the gradient of (2.12),
namely

$$\partial_t s = t^{-2} \ s \cdot \nabla s + t^{-\gamma} \ \nabla g_0 (w, w) \quad . \eqno(2.27)$$ 

We are now in a position to describe in more detail the contents of the technical sections 3-7 of
this paper. In Section 3, we introduce the relevant Gevrey spaces and derive the basic estimates
in those spaces that are needed to study the system (2.11)-(2.12) (Lemmas 3.4-3.7), we
explain in passing the mechanism by which those spaces make it possible to overcome the
derivative loss in (2.11)-(2.12) for $\mu > n - 2$ (after Lemma 3.4) and finally we
introduce the estimating functions of time mentioned above and obtain some estimates for
them. In Section 4, we prove that the Cauchy problem for the system (2.11)-(2.12) is
well-posed for large time, with large but finite initial time (Proposition 4.1), we prove
the existence of a limit for $w(t)$ as $t \to \infty$ for the solutions thereby obtained
(Proposition 4.2) and we derive a uniqueness result of solutions with prescribed asymptotic
behaviour (Proposition 4.3). In Section 5, we study the asymptotic behaviour in time of the
solutions obtained in Section 4. We derive a number of properties and estimates for the
solutions of the asymptotic system (2.18)-(2.19), defined inductively (Proposition 5.2). We
then obtain asymptotic estimates on the approximation of the solutions of the system
(2.11)-(2.12) by the asymptotic functions $(W_m, \phi_m)$ defined by (2.14)-(2.15), and in
particular we complete the proof of the existence of asymptotic states for those solutions
(Proposition 5.3). In Section 6, we study the Cauchy problem with infinite initial time,
first for the transport equations (2.23) (Proposition 6.1) and (2.24) (Proposition 6.2), and
then for the system (2.11)-(2.12). For a given solution $(V, \chi )$ of the system
(2.23)-(2.24) and a given (large) $t_0$, we construct a solution $(w_{t_0}, \varphi_{t_0})$
of the system (2.11)-(2.12) which coincides with $(V, \phi_p + \chi )$ at $t_0$ and we
estimate it uniformly in $t_0$ (Proposition 6.4). We then prove that when $t_0 \to \infty$,
$(w_{t_0}, \varphi_{t_0})$ has a limit $(w, \varphi )$ which is asymptotic both to $(V,
\phi_p + \chi )$ and to $(W_p , \phi_p + \psi_+)$ (Proposition 6.5). In Section 7, we
exploit the results of Section 6 to construct the wave operators for the equation (1.1) and
to describe the asymptotic behaviour of solutions in their range. We first prove that the
local wave operator at infinity for the system (2.11)-(2.12) defined through Proposition 6.5
in Definition 7.1 is gauge covariant in the sense of Definitions 7.2 and 7.3 in the best
form that can be expected with the available regularity (Propositions 7.2 and 7.3). With the
help of some information on the Cauchy problem for (1.1) at finite time (Proposition 7.1),
we then define the wave operator $\Omega : u_+ \to u$ (Definition 7.4), and we prove that it
is injective.  We then collect all the available information on $\Omega$ and on solutions of
(1.1) in its range in Proposition 7.5, which contains the main results of this paper.
Finally some side results relevant for the definition and properties of the Gevrey spaces
used here are collected in two Appendices.

\section{Gevrey spaces and preliminary estimates}
\hspace*{\parindent} In this section, we define the Gevrey spaces where we shall study the
auxiliary system (2.11)-(2.12) and we derive a number of energy type estimates which hold in
those spaces and play an essential role in that study. We then introduce some estimating
functions of time generalizing those of II and we derive a number of estimates for them. \par

The relevant spaces will be defined with the help of the functions
$$f_0(\xi ) = \exp \left ( \rho |\xi |^{\nu} \right ) \ , \ f(\xi ) = \exp \left ( \rho (|\xi |^{\nu}
\vee 1) \right ) \eqno(3.1)$$  
\noindent where $0 < \nu \leq 1$, $\rho$ is a positive parameter to be specified later, and $\xi
\in {I\hskip-1truemm R}^n$. The dependence of $f$ on $\rho$ will always be omitted in the notation.
\par

In all this paper, one could use instead of $f_0$ the function
$$\widetilde{f}(\xi ) = \sum_{j\geq 0} \ (j!)^{-1/\nu} \ \rho^{j/\nu} |\xi|^j \eqno(3.2)$$
\noi which satisfies the same basic estimates and would yield essentially the same results. The
function $\widetilde{f}$ is also convenient in order to relate the definition of the Gevrey
spaces $K_{\rho}^k$ and $Y_{\rho}^{\ell}$ (see (3.8) (3.9) below) to more standard definitions.
Those points are discussed in Appendix A. \par

The functions $f_0$ and $f$ satisfy the following estimates. \\

\noi {\bf Lemma 3.1.} {\it Let $\xi$, $\eta \in {I\hskip-1truemm R}^n$. Then $f$ satisfies the
estimates~:
$$f( \xi ) \leq f(\xi - \eta ) \ f(\eta ) \qquad \hbox{\it for all $\xi$ , $\eta$} \ ,\eqno(3.3)$$
$$f( \xi ) \leq f(\xi - \eta ) \ f_0 (\eta )^{\nu} \qquad \hbox{\it for} \ |\xi | \wedge |\eta| \leq
|\xi - \eta | \quad , \eqno(3.4)$$
$$|f(\xi ) - f(\eta ) | \ |\eta |^{1 - \nu} \leq |\xi - \eta |^{1- \nu} \ f(\xi - \eta ) \ f(\eta
) \qquad \hbox{\it for all $\xi$, $\eta$} \quad , \eqno(3.5)$$
$$|f(\xi ) - f(\eta ) | \ |\eta |^{1 - \nu} \leq C |\xi - \eta |^{1- \nu} \ f_0(\xi - \eta )^{\nu} \
f(\eta ) \qquad \hbox{\it for} \ |\xi| \wedge |\xi - \eta | \leq |\eta | \quad , 
\eqno(3.6)$$ 
$$|f(\xi ) - f(\eta ) | \ |\eta |^{1 - \nu} \leq C |\xi - \eta |^{1- \nu} \ f(\xi - \eta ) \
f_0(\eta )^{\nu} \qquad \hbox{\it for} \ |\xi| \wedge |\eta | \leq |\xi - \eta | \quad . 
\eqno(3.7)$$ 

In (3.6) and (3.7), one can take $C = 1$, except in the region $|\xi | \leq |\xi - \eta | \leq |\eta
|$ where $C = 2^{1 - \nu}$. \par

The function $f_0$ satisfies the same estimates as $f$.} \\

\noi {\bf Proof.} We first prove the estimates for $f_0$ \par

(3.3) follows from the fact that $|\xi|^{\nu} \leq |\xi - \eta |^{\nu} + |\eta |^{\nu}$. \\

(3.4) is obvious for $|\xi | \leq |\xi - \eta |$. For $|\eta | \leq |\xi - \eta | \leq |\xi |$, we
estimate
$$|\xi |^{\nu} \leq |\xi - \eta |^{\nu} + \nu |\xi - \eta |^{\nu - 1} \ |\eta | \leq |\xi - \eta
|^{\nu} + \nu |\eta |^{\nu} \quad .$$

(3.5) follows from (3.3) for $|\eta| \leq |\xi - \eta |$ and from (3.6) with $C = 1$ for $|\xi -
\eta | \leq |\xi | \wedge |\eta |$. For $|\xi | \leq |\xi - \eta | \leq |\eta |$, we estimate
$$\left ( f_0(\eta ) - f_0 (\xi )\right ) |\eta |^{1 - \nu} \leq f_0(\eta ) \ |\eta |^{1 - \nu} \left
( 1 - \exp (- \rho |\eta |^{\nu} ) \right ) \leq f_0(\eta ) \ \rho \ |\eta |$$
$$\leq f_0 (\eta ) \left ( |\eta |/|\xi - \eta | \right ) |\xi - \eta |^{1 - \nu} \ \rho|\xi -
\eta |^{\nu}$$
\noi and the result follows from the fact that $|\eta | \leq 2 |\xi - \eta |$ and $\rho |\xi -
\eta |^{\nu} \leq e^{-1} f_0(\xi - \eta )$. \\

(3.6) is obvious for $|\xi | \leq |\eta | \leq |\xi - \eta |$, with $C = 1$. \par

For $|\xi | \leq |\xi - \eta | \leq |\eta |$, we estimate
$$\left ( f_0(\eta ) - f_0(\xi ) \right ) |\eta |^{1 - \nu} \leq 2^{1 - \nu} \ |\xi - \eta |^{1 -
\nu} \ f_0(\eta )$$
\noi which yields (3.6) with $C = 2^{1 - \nu}$ since $|\eta | \leq 2 |\xi - \eta |$. \par

The really crucial case is the case $|\xi - \eta | \leq |\xi | \wedge |\eta |$. \par

For $|\xi - \eta | \leq |\xi | \leq | \eta |$, we estimate
$$\left ( f_0(\eta ) - f_0(\xi ) \right ) |\eta |^{1-\nu} \leq f_0(\eta ) \ |\eta |^{1 - \nu} \left
( 1 - \exp (- \rho \nu |\xi - \eta | \ |\xi |^{\nu - 1}) \right )$$ 
\noi where we have used
$$|\eta |^{\nu} - |\xi |^{\nu} \leq \nu \ |\xi |^{\nu - 1} \ |\xi - \eta | \quad ,$$
$$\cdots \leq f_0 (\eta ) \left ( |\eta |/|\xi |\right )^{1 - \nu} \ \rho \ \nu \ |\xi - \eta |$$
$$\leq f_0(\eta ) |\xi - \eta |^{1 - \nu} \ 2^{1 - \nu} \ e^{-1} \ f_0(\xi - \eta )^{\nu}$$
\noi since $|\eta | \leq 2 |\xi |$ and $\rho \nu |\xi - \eta |^{\nu} \leq e^{-1} f_0 (\xi - \eta
)^{\nu}$. This proves (3.6) with $C = 1$ in that case. \par

For $|\xi - \eta | \leq |\eta | \leq |\xi |$, we estimate similarly
$$\left ( f_0(\xi ) - f_0(\eta ) \right ) |\eta |^{1 - \nu} \leq f_0 (\eta ) \left \{ |\eta |^{1-
\nu} \left ( \exp (\rho \nu |\xi - \eta |\ |\eta |^{\nu-1} ) - 1 \right ) \right \} \quad .$$
\noi Now for fixed $|\xi - \eta |$, the last bracket is a decreasing function of $|\eta |$, and is
therefore bounded by its value for $|\eta | = |\xi - \eta |$, so that 
$$\cdots \leq f_0(\eta ) \ |\xi - \eta |^{1 - \nu} \left ( f_0(\xi - \eta )^{\nu} - 1 \right )$$
\noi which proves (3.6) with $C = 1$ in that case. \\

(3.7) is obvious for $|\xi | \leq |\eta | \leq |\xi - \eta |$ and follows from (3.4) with $C = 1$
for $|\eta | \leq |\xi | \wedge |\xi - \eta |$. For $|\xi | \leq |\xi - \eta | \leq |\eta |$, we
estimate
$$\left ( f_0(\eta ) - f_0(\xi ) \right ) |\eta |^{1 - \nu} \leq 2^{1 - \nu} |\xi - \eta |^{1 -
\nu} \ f_0(\eta )$$
\noi and (3.7) with $C = 2^{1 - \nu}$ follows from 
$$|\eta |^{\nu} \leq \nu |\eta |^{\nu} + (1 - \nu) 2^{\nu}\  |\xi - \eta |^{\nu} \leq \nu |\eta
|^{\nu} + |\xi - \eta |^{\nu} \quad .$$

The estimates for $f$ follow from those for $f_0$. This is obvious for (3.3) (3.4). For (3.5) (3.6)
(3.7), it follows from the fact that for all $\xi$, $\eta$ and all $a > 0$
$$|f_0(\xi ) \vee a - f_0 (\eta ) \vee a| \leq |f_0(\xi ) -f_0(\eta )| \quad .$$\par
\hfill $\sq$\par

We now turn to the definition of the spaces where we shall solve the system (2.11)-(2.12).
For any tempered distribution $u$ in ${I\hskip-1truemm R}^n$ with $\widehat{u} \in L_{loc}^1
({I\hskip-1truemm R}^n)$, we define $u_{>\atop <}$ by $\widehat{u}_> (\xi ) =
\widehat{u}( \xi )$ for $|\xi | > 1$, $\widehat{u}_>(\xi ) = 0$ for $|\xi | \leq 1$, $\widehat{u}_<
(\xi ) = 0$ for $|\xi | > 1$, $\widehat{u}_< (\xi ) = \widehat{u}(\xi )$ for $|\xi | \leq 1$.
Similarly, for $\xi \in {I\hskip-1truemm R}^n$ and $m \in {I\hskip-1truemm R}$, we define $|\xi
|_>^m$ and $|\xi |_<^m$ to be equal to $|\xi |^m$ for $|\xi | > 1$ and $|\xi | \leq 1$ respectively,
and zero otherwise. Occasionally we shall make the separation between low and high $|\xi |$ at some
value $a \not= 1$. In that case we shall denote by $u_{<a}$ and $u_{>a}$ the corresponding
components of $u$. \par

Let now $\rho > 0$, $k \in {I\hskip-1truemm R}$, $\ell \in {I\hskip-1truemm R}$ and $0 \leq
\ell_< < n/2$. Starting from Lemma 3.4 below (but not until then), we shall assume in addition
that $\ell_< > n/2 - \mu$. We define
$$K_{\rho}^k = \left \{ w:|w|_k^2 \ \equiv \ \parallel w;K^k_{\rho}\parallel^2 \ \equiv \ \parallel
|\xi |^k \ f(\xi ) \ \widehat{w}_> (\xi ) \parallel_2^2 + \parallel f(\xi ) \ \widehat{w}_<(\xi ) 
\parallel_2^2 \ < \infty \right \} \ , \eqno(3.8)$$ 
$$Y_{\rho}^{\ell} = \Big \{ \varphi : \widehat{\varphi} \in L_{loc}^1({I\hskip-1truemm R}^n) \
\hbox{and} \ |\varphi |_{\ell}^2 \ \equiv \ \parallel \varphi; Y^{\ell}_{\rho}\parallel^2 \ \equiv
\ \parallel |\xi |^{\ell + 2} \ f(\xi ) \ \widehat{\varphi}_> (\xi ) \parallel_2^2 $$ 
$$ + \parallel |\xi |^{\ell_<} \ f(\xi ) \ \widehat{\varphi}_< (\xi ) \parallel_2^2 \ < \infty \Big
\} \quad . \eqno(3.9)$$    
\noindent The apparent ambiguity in the notation $|\cdot |_b$ will be lifted by the fact that the
symbol $b$ will always contain the letter $k$ when referring to $K_{\rho}^k$ spaces and the letter
$\ell$ when referring to $Y_{\rho}^{\ell}$ spaces. The spaces $K_{\rho}^k$ and $Y_{\rho}^{\ell}$
are Hilbert spaces and satisfy the embeddings $K_{\rho}^k \subset K_{\rho}^{k'}$ for $k \geq k'$
and $Y_{\rho}^{\ell} \subset Y_{\rho}^{\ell '}$ for $\ell \geq \ell '$, $K_{\rho}^k \subset
K_{\rho '}^k$ and $Y_{\rho}^{\ell} \subset Y_{\rho '}^{\ell}$, for $\rho \geq \rho '$. \\

\noi {\bf Remark 3.1.} The norms in the spaces $K_{\rho}^k$ and $Y_{\rho}^{\ell}$ are both of the
form $\parallel f_1\ f \ \widehat{u}\parallel_2$ with
$$f_1 = |\xi |_>^{k_>} + |\xi|_<^{k_<} \eqno(3.10)$$
\noi where $(k_>, k_<) = (k,0)$ for $K_{\rho}^k$ and $(k_>,k_<) = (\ell + 2, \ell_<)$ for
$Y_{\rho}^{\ell}$. In particular this implies that 
$$\left \{ \begin{array}{l} \parallel w; K_{\rho}^k\parallel \  = \ \parallel F^{-1} (f\widehat{w});
K_0^k\parallel \\ \\ \parallel \varphi ; Y_{\rho}^{\ell} \parallel \ = \ \parallel F^{-1} (f
\widehat{\varphi});Y_0^{\ell}\parallel \quad . \end{array}\right . \eqno(3.11)$$

If $k_> \geq k_<$, namely if $k \geq 0$ for $K_{\rho}^k$ and if $\ell + 2 \geq \ell_<$ for
$Y_{\rho}^{\ell}$ (which will always be the case in the applications), one can omit either or both
of the upper and lower restrictions in the definition of $K_{\rho}^k$ or $Y_{\rho}^{\ell}$, thereby
obtaining equivalent norms uniformly in $\rho$. In fact, in that case 
$$|\xi |_>^{2k_>} + |\xi |_<^{2k_<} \leq |\xi |^{2k_>} + |\xi |^{2k_<} \leq 2 \left ( |\xi
|_>^{2k_>} + |\xi |_<^{2k_<} \right ) \ .$$
Furthermore, the relations (3.11) are preserved under thoses changes. \\

We shall use the spaces $K_{\rho}^k$ and $Y_{\rho}^{\ell}$ with a time dependent parameter $\rho \in
{\cal C}^1({I\hskip-1truemm R}^+, {I\hskip-1truemm R}^+)$. The form (3.8)-(3.9) of the norms
has been chosen so as to ensure that for fixed (time independent) $w$ and $\varphi$, the
following relations hold
$$\left \{ \begin{array}{l} \displaystyle{{d \over dt}} \ |w|_k^2 = 2 \rho ' |w|_{k+ \nu /2}^2
\quad , \\ \\ \displaystyle{{d \over dt}} \ |\varphi |_{\ell}^2 = 2 \rho ' |\varphi |_{\ell + \nu
/2}^2 \quad ,\end{array}\right . \eqno(3.12)$$
\noindent where $\rho ' = d\rho /dt$. \par

We shall look for $w$ as complex $K_{\rho}^k$ valued functions of time and for $\varphi$ as real
$Y_{\rho}^{\ell}$ valued functions of time. More precisely, we shall look for $(w, \varphi )$ such
that for some interval $I \subset [1, \infty )$
$$(w, \varphi ) \in {\cal X}_{\rho, loc}^{k,\ell} (I) \equiv {\cal C} \left (I, K_{\rho}^k \oplus
Y_{\rho}^{\ell} \right )  \cap L_{loc}^2 \left ( I, K_{\rho}^{k+ \nu /2} \oplus Y_{\rho}^{\ell + \nu
/ 2} \right ) \eqno(3.13)$$
\noi by which is meant, especially as regards continuity, that 
$$\left ( F^{-1} f \widehat{w}, F^{-1} f \widehat{\varphi} \right ) \in
{\cal X}_{0,loc}^{k,\ell}(I) \equiv {\cal C} \left ( I, K_0^k \oplus Y_0^{\ell} \right ) \cap
L_{loc}^2 \left ( I, K_0^{k + \nu /2} \oplus Y_0^{\ell + \nu /2} \right )$$
\noi in the usual sense. In particular, when taking norms such as $|w(t)|_k$ or $|\varphi (t)
|_{\ell}$ with time dependent $\rho$, it will always be understood that $\rho$ in the definition
of the relevant space is taken at the same value of the time as $w$ or $\varphi$. \par

We shall also need global versions of the previous spaces, especially when the interval $I$ is
unbounded. The definition of those global versions will require assumptions on $\rho$ that are
irrelevant for the considerations of this section and will be postponed until the beginning of
Section 4. \par

We shall need the following elementary estimates. \\

\noi {\bf Lemma 3.2.} {\it Let $m \in {I\hskip-1truemm R}$. The following estimates hold~:
$$\parallel |\xi |^m \widehat{(u_1u_2)}_>\parallel_2 \ \leq \ C \parallel <\xi >^{k_1}
\widehat{u}_1\parallel_2 \ \parallel <\xi >^{k_2} \widehat{u}_2 \parallel_2 \eqno(3.14)$$
\noi for $k_1$, $k_2 \geq m \vee 0$ and $k_1 + k_2 > m + n/2$, where $<\cdot> = (1 + |\cdot
|^2)^{1/2}$, $$\parallel |\xi |^m \widehat{(u_1u_2)}_<\parallel_2 \ \leq \ C \parallel
u_1\parallel_2 \ \parallel u_2\parallel_2 \eqno(3.15)$$
\noi for $m > - n/2$.} \\

\noi {\bf Proof.} (3.14) for $m \geq 0$ follows from (3.14) (3.15) for $m = 0$ and from 
$$\parallel | \xi |^m \widehat{(u_1u_2)} \parallel_2 \ \leq \ C \left ( \parallel (|\xi |^m
|\widehat{u}_1|) * |\widehat{u}_2 | \parallel_2 \ + \ \parallel |\widehat{u}_1| * (|\xi |^m
|\widehat{u}_2|) \parallel_2 \right ) \ .$$  

For $m < 0$, we estimate by the H\"older and Young inequalities 
$$\parallel | \xi |^m \widehat{(u_1u_2)}_> \parallel_2 \ \leq \ C \parallel |\xi |^m_>
\parallel_s \ \parallel \widehat{u}_1 \parallel_{\bar{r}_1} \ \parallel 
\widehat{u}_2 \parallel_{\bar{r}_2}$$
\noi with $1/s + 1/\bar{r}_1 + 1 /\bar{r}_2 = 3/2$ , $n/s < |m|$ and $1 \leq \bar{r}_1$,
$\bar{r}_2 \leq 2$. The last two norms are estimated by $\parallel < \xi >^{k_i}
\widehat{u}_i\parallel_2$ provided $k_i > n/\bar{r}_i - n/2$. One can find $\bar{r}_1$,
$\bar{r}_2$ satisfying all previous conditions under the assumptions made on $k_i$.  \par

For $m = 0$, we apply the same argument with $s = \infty$. \\

For the proof of (3.15) we estimate by the Schwarz and Young inequalities 
$$\parallel |\xi |^m \widehat{(u_1u_2)}_<\parallel_2 \ \leq \ \parallel |\xi |_<^m\parallel_2 \
\parallel u_1\parallel_2 \ \parallel u_2\parallel_2$$ \noindent for $m > - n/2$. \par
\hfill $\sq$ \par

For future reference, we also state the following elementary inequalities
$$\parallel |\xi |^m \widehat{\varphi}_<\parallel_1 \ \leq \ C \parallel |\xi |^{\ell_<} \ 
\widehat{\varphi}_<\parallel_2 \ \hbox{for all} \ m \geq 0 \ , \eqno(3.16)$$
$$\parallel |\xi |^m \widehat{\varphi}_>\parallel_1 \ \leq \ C \parallel |\xi |^{\ell +2} \ 
\widehat{\varphi}_>\parallel_2 \ \hbox{for } \ \ell + 2 > m + n/2 \ .  \eqno(3.17)$$
 \noi In what follows, we
shall repeatedly estimate norms such as $\parallel |\xi |^m f\widehat{u_1u_2}\parallel_2$ with $m
\geq 0$. For that purpose, using (3.3), we shall write 
$$|\xi |^m f |\widehat{u_1u_2}| \leq |\xi |^m \int d\eta \ f(\xi ) |\widehat{u}_1 (\xi - \eta )| \
|\widehat{u}_2(\eta ) |$$
$$\leq 2^m \int d\eta \left ( |\xi - \eta |^m + |\eta |^m \right ) f(\xi - \eta ) \ f(\eta )
|\widehat{u}_1(\xi - \eta )| \ |\widehat{u}_2 (\eta ) |$$
$$= 2^m \left \{ \left ( |\cdot |^m f |\widehat{u}_1| \right ) * \left ( f|\widehat{u}_2| \right )
+ (1 \leftrightarrow 2 ) \right \} \ . \eqno(3.18)$$ 
\noi That inequality will be often combined with restrictions to low or high values of $|\xi
|$, either in the product $(u_1u_2)$ or in $u_1$ or $u_2$ separately.  \par

The next lemma states that under suitable assumptions on $k$ and $\ell$, $Y_{\rho}^{\ell}$ is an
algebra under ordinary multiplication and acts boundedly on $K_{\rho}^k$ by multiplication.\\

\noi {\bf Lemma 3.3.} {\it Let $\ell + 2 > n/2$ and $0 \leq k \leq \ell + 2$. Then there exist
constants $C_1$ and $C_2$, independent of $\rho$, such that
$$| \varphi \psi |_{\ell} \leq C_1 |\varphi |_{\ell} \ |\psi |_{\ell} \quad  \hbox{for all
$\varphi$,} \ \psi \in Y_{\rho}^{\ell} \ , \eqno(3.19)$$ 
$$| \varphi w |_{k} \leq C_2 |\varphi |_{\ell} \ |w |_{k} \quad  \hbox{for all $\varphi \in
Y_{\rho}^{\ell}$,} \ w \in K_{\rho}^{k} \ . \eqno(3.20)$$
\noi In particular
$$|(\exp (- i \varphi ) - 1) w|_k \leq C_2 \ C_1^{-1} \left ( \exp (C_1 |\varphi |_{\ell} ) - 1
\right ) |w|_k \eqno(3.21)$$ \noi for all $\varphi \in Y^{\ell}_{\rho}$, $w \in K_{\rho}^k$.} \\

\noi {\bf Proof.} From the definitions (3.8) (3.9) of the norms and in particular from (3.11), and
from (3.3), it follows that (3.19) (3.20) for $\rho = 0$ imply the same estimates for arbitrary
$\rho > 0$ with the same constants. We therefore restrict our attention to the case $\rho = 0$. In
that case (3.19) (3.20) are almost standard properties of Sobolev spaces, except for the presence
of $|\xi |_<^{\ell_<}$ with possibly $\ell_< > 0$ in (3.9). We give a proof for completeness. \\

(3.19). From the definition (3.9) with $\rho = 0$, it follows that it is sufficient to estimate
$\parallel | \xi |^{\ell + 2} (\widehat{\varphi \psi})_{>2}\parallel_2$ and $\parallel | \xi
|^{\ell_<} \ \widehat{\varphi \psi}\parallel_2$. We estimate
$$\parallel | \xi |^{\ell + 2} (\widehat{\varphi \psi})_{>2}\parallel_2 \ \leq C \Big \{
\parallel \left ( |\xi |^{\ell + 2} \ |\widehat{\varphi}_>| \right ) * |\widehat{\psi}|
\parallel_2$$  
$$+ \ \parallel |\widehat{\varphi}_<| * |\widehat{\psi}_>|\parallel_2 \ + (\varphi
\leftrightarrow \psi ) \Big \}$$
$$\leq C \Big \{ \parallel |\xi |^{\ell + 2} \ \widehat{\varphi}_> \parallel_2 \ \parallel
\widehat{\psi}\parallel_1 \ + \ \parallel \widehat{\varphi}_< \parallel_1 \ \parallel
\widehat{\psi}_>\parallel_2 + (\varphi \leftrightarrow \psi ) \Big \}$$
$$\leq C |\varphi |_{\ell} \ |\psi |_{\ell}$$
\noi by (3.3), by the Young inequality and by (3.16) (3.17). The lower restriction $|\xi | > 2$ in
$\varphi \psi$ implies that there is no $\varphi_< \ \psi_<$ contribution. \par

On the other hand
$$\parallel |\xi |^{\ell <} \ \widehat{\varphi \psi}\parallel_2 \ \leq \ C \Big \{ \parallel
|\xi |^{\ell_<}\ \widehat{\varphi}\parallel_2 \ \parallel \widehat{\psi} \parallel_1 \ + (\varphi
\leftrightarrow \psi ) \Big \} \leq C |\varphi |_{\ell}\  |\psi |_{\ell}$$
  \noi by the Young
inequality and (3.16) (3.17). \\

(3.20). We estimate similarly by (3.3) and the Young inequality
$$\parallel <\xi >^k \widehat{\varphi w}\parallel_2 \ \leq \ C \Big \{ \parallel
\widehat{\varphi}\parallel_1 \ \parallel <\xi >^k \widehat{w}\parallel_2 \ + \ \parallel
\widehat{\varphi}_<\parallel_1 \ \parallel \widehat{w}\parallel_2$$ $$+ \ \parallel (|\xi |^k
|\widehat{\varphi}_>|) * |\widehat{w}|\parallel_2 \Big \}$$ \noi and the result follows from
(3.16) (3.17) and from Lemma 3.2 with $m = 0$, $k_1 = \ell + 2 - k$ and $k_2 = k$. \\

(3.21) follows immediately from a repeated application of (3.19) (3.20). \par
\hfill $\sq$ \par

\noi {\bf Remark 3.2.} In Lemma 3.3 we have used only (3.3) from Lemma 3.1. For $\nu < 1$, by
u\-sing (3.4), one can obtain more general results. In particular $K_{\rho}^k$ and $Y_{\rho}^{\ell}$
are algebras under multiplication for any $k$ and $\ell$ in ${I\hskip-1truemm R}^+$ (see Appendix B).
Similarly in what follows, we shall use only (3.3) and (3.5). For $\nu < 1$, by using in addition
(3.4) and (3.6), one could generalize some of the results by weakening the assumptions made on $k$
and $\ell$. That extension however would not hold uniformly in $\nu$ for $\nu \to 1$, and we shall
therefore refrain from following that track. \\

We now turn to the derivation of the basic estimates needed to study the auxiliary system
(2.11)-(2.12). For that purpose, we shall use a regularization. Let $j \in {\cal
C}^1({I\hskip-1truemm R}^n, {I\hskip-1truemm R})$ with $0 \leq j \leq 1$ and $j(0) = 1$. We
denote by $j_{\varepsilon}$ both the function $j_{\varepsilon}(\xi ) = j(\varepsilon \xi )$
and the operator of multiplication by that function in Fourier space variables, and by
$J_{\varepsilon}$ the operator $F^*j_{\varepsilon}F$. \par

We recall that $g_0$ is defined by (2.13). We shall use systematically the notation $s= \nabla
\varphi$ and whenever convenient, start from the equation (2.27) satisfied by $s$ instead of
the equation (2.12) satisfied by $\varphi$. \par

We first state the basic estimates. \\   

\noi {\bf Lemma 3.4.} {\it Let $m \geq 0$. The following estimates hold~: 
$$\left |{\rm Re} \ < j_{\varepsilon} |\xi |^k f \widehat{w}, j_{\varepsilon} |\xi |^k
f(\widehat{s\cdot \nabla w})> \right | + \left | {\rm Re} \ < j_{\varepsilon}
f\widehat{w}, j_{\varepsilon} f(\widehat{s\cdot \nabla w})> \right | \leq C |w|_{k+ \nu /2}^2
\ |\varphi |_{\ell} \eqno(3.22)$$ \noi uniformly in $\varepsilon$, for $\ell > n/2 - \nu$ , $k \geq
\nu /2$ , $\ell + 1 \geq k - \nu /2$.
$$\parallel |\xi |^m f(\widehat{s\cdot \nabla w})\parallel_2 \ \leq C |\varphi|_{\ell} \ |w|_k
\eqno(3.23)$$
\noi for $k + \ell > m + n/2$ , $k \geq m + 1$ , $\ell + 1 \geq m$.
$$\parallel |\xi |^m f((\widehat{\nabla\cdot  s)w})\parallel_2 \ \leq C |\varphi|_{\ell} \ |w|_k
\eqno(3.24)$$
\noi for $k + \ell > m + n/2$ , $k \geq m$ , $\ell \geq m$.
$$\left |{\rm Re} \ < j_{\varepsilon} |\xi |^{\ell ' +1} f \widehat{s'_>}, j_{\varepsilon} |\xi
|^{\ell ' + 1} f(\widehat{s\cdot \nabla s'})_> > \right | \leq C |\varphi '|_{\ell '+ \nu /2}^2
\ |\varphi |_{\ell} \eqno(3.25)$$
\noi uniformly in $\varepsilon$, for $\ell > n/2 - \nu$ , $\ell ' + 1 \geq \nu /2$ , $\ell \geq \ell
' - \nu /2$.
$$\parallel |\xi |^m f(\widehat{s\cdot \nabla s'})_>\parallel_2 \ \leq C |\varphi|_{\ell} \
|\varphi '|_{\ell '} \eqno(3.26)$$
\noi for $\ell + \ell ' > m - 1 + n/2$ , $\ell + 1 \geq m$ , $\ell ' \geq m$.
$$\parallel |\xi |^{\ell_<} f(\widehat{s\cdot s'})\parallel_2 \ \leq C |\varphi|_{\ell} \
|\varphi '|_{\ell '} \eqno(3.27)$$
\noi for $\ell$, $\ell ' > n/2 - 2$.
$$\parallel |\xi |^m f\widehat{g_0(w_1}w_2)_>\parallel_2 \ \leq C \left \{ |w_1|_{k_1} \
|w_2|_{k_2} + |w_1|_{k'_1} \ |w_2|_{k'_2} \right \} \eqno(3.28)$$
\noi for $(k_1 + k_2) \wedge (k'_1 + k'_2) > \beta + n/2$ , $k_1 \wedge k'_2 \geq \beta \vee 0$ ,
$k_2 \wedge k'_1 \geq 0$, where $\beta = m + \mu - n$.
$$\parallel |\xi |^{\ell_<} \widehat{g_0(w_1}w_2)_<\parallel_2 \ \leq C \parallel w_1\parallel_2 \
\parallel w_2 \parallel_2 \eqno(3.29)$$
\noi for $\ell_< > n/2 - \mu$. \par

The constants $C$ in (3.22)-(3.29) can be taken independent of $\rho$. }\\

\noi {\bf Proof.} (3.22). We have to estimate
$${\rm Im} \int d\xi \ d\eta \ j_{\varepsilon}(\xi ) \ |\xi |^k \ f(\xi ) \ \bar{\widehat{w}}(\xi )
\ j_{\varepsilon}(\xi ) \ |\xi |^k \ f(\xi ) \ \widehat{s} (\xi - \eta ) \cdot \eta \ \widehat{w}
(\eta ) \eqno(3.30)$$
\noi and a similar expression with $k = 0$. We consider only the former one. The proof for the
latter is similar and simpler. We split the domain of integration into three regions, namely 
$$|\xi - \eta | \leq |\xi | \wedge |\eta |\ , \ |\xi | \leq |\xi - \eta | \wedge |\eta | \
\hbox{and} \ |\eta | \leq |\xi | \wedge |\xi - \eta |$$
\noi and correspondingly the integral (3.30) is written as the sum $I_1 + I_2 + I_3$ of three
terms which we estimate successively. \par

\noi \underbar{Region} $|\xi - \eta | \leq |\xi | \wedge |\eta |$, \underbar{estimate of} $I_1$.
\par
In this region we decompose the integrand according to the identity 
$$j_{\varepsilon} (\xi ) |\xi |^k f(\xi ) \eta = j_{\varepsilon} (\xi ) |\xi |^k \left ( f(\xi ) -
f(\eta ) \right ) \eta + \left ( j_{\varepsilon}(\xi )\  |\xi |^k - j_{\varepsilon}(\eta )
\ |\eta |^k \right ) f(\eta ) \eta + j_{\varepsilon} (\eta )\  |\eta |^k \ f(\eta )
\eta\eqno(3.31)$$  \noi and correspondingly $I_1$ is written as the sum $I_{1,1} + I_{1,2} +
I_{1,3}$ of three terms which we estimate successively. \par

\noi \underbar{Estimate of} $I_{1,1}$. From (3.5) of Lemma 3.1, we obtain 
$$\left | j_{\varepsilon} (\xi ) |\xi |^{2k} (f(\xi ) - f(\eta ))\eta \right | \leq C |\xi |^{k+ \nu
/2} \ |\xi - \eta |^{1 - \nu} \ |\eta |^{k+ \nu /2} f(\xi - \eta ) f(\eta )$$
\noi and therefore by the Schwarz and Young inequalities
$$|I_{1,1}| \leq C |w|_{k + \nu /2}^2 \ \parallel |\xi |^{1 - \nu} f\widehat{s} \parallel_1$$
\noi so that by (3.16) (3.17)
$$|I_{1,1}| \leq C |w|_{k + \nu /2}^2 \ |\varphi |_{\ell} \quad . \eqno(3.32)$$

\noi  \underbar{Estimate of} $I_{1,2}$. We rewrite
$$\left ( j_{\varepsilon} (\xi ) \ |\xi |^{k} - j_{\varepsilon} (\eta ) \ |\eta |^k \right ) \eta =
j_{\varepsilon} (\xi ) \left ( |\xi |^k - |\eta |^k \right ) \eta + |\eta |^k \left (
j_{\varepsilon}(\xi) - j_{\varepsilon}(\eta ) \right ) \eta \ . \eqno(3.33)$$
\noi We estimate
$$\left | | \xi |^k - |\eta |^k \right | \  |\eta | \leq k \ 2^{|k-1|} \ |\xi - \eta| \ |\eta |^k
\eqno(3.34)$$
\noi and we rewrite
$$\left ( j_{\varepsilon} (\xi ) - j_{\varepsilon} (\eta ) \right ) \eta = j_{\varepsilon}(\xi )
\xi -  j_{\varepsilon}(\eta ) \eta - j_{\varepsilon}(\xi ) (\xi - \eta )$$ 
$$= \int_0^1 d\theta \left \{ (\xi - \eta ) \cdot \nabla j_{\varepsilon} (\xi_{\theta})
\xi_{\theta} + (j_{\varepsilon} (\xi_{\theta}) - j_{\varepsilon}(\xi )) (\xi - \eta ) \right \}
\eqno(3.35)$$    
\noi with $\xi_{\theta} = \theta \xi + (1 - \theta )\eta$, so that 
$$| j_{\varepsilon}(\xi ) - j_{\varepsilon}(\eta )| \ |\eta | \leq \left ( \parallel |\cdot|
\nabla j \parallel_{\infty} \ + 2 \right ) |\xi - \eta | \quad . \eqno(3.36)$$
\noi Comparing (3.33) (3.34) (3.36), we obtain
$$\left | j_{\varepsilon}(\xi ) |\xi |^k - j_{\varepsilon}(\eta ) |\eta |^k \right |\  | \eta |
\leq C |\xi - \eta| \ |\eta|^k \leq C |\xi |^{\nu /2} \ |\xi - \eta |^{1 - \nu} \ |\eta |^{k+ \nu/2}
\eqno(3.37)$$ \noi from which we obtain as previously $$| I_{1,2}| \leq C |w|_{k+ \nu /2}^2 \
|\varphi |_{\ell} \quad .\eqno(3.38)$$

\noi \underbar{Estimate of} $I_{1,3}$. By reality and symmetry, which is the Fourier space version
of integration by parts, $I_{1,3}$ can be rewritten as 
$$I_{1,3} = (i/2) \int_{|\xi - \eta | \leq |\xi | \wedge |\eta |} d \xi \ d\eta \ j_{\varepsilon}(
\xi ) \ |\xi|^k \ f(\xi ) \ \bar{\widehat{w}} (\xi ) (\xi - \eta ) \cdot \widehat{s} (\xi - \eta )
\ j_{\varepsilon}(\eta )\  |\eta|^k\  f(\eta ) \ \widehat{w}(\eta) \ .$$
\noi Using the second inequality in (3.37) we obtain as previously
$$|I_{1,3}| \leq C |w|_{k+ \nu /2}^2 \ |s|_{\ell } \quad . \eqno(3.39)$$
\noi We now turn to the contribution of the region $|\xi | \leq |\xi - \eta | \wedge |\eta |$.
\par

\noi \underbar{Estimate of} $I_2$. Using (3.3) and
$$|\xi |^{2k} |\eta | \leq C |\xi |^{k+\nu /2} \ |\xi - \eta |^{1 - \nu } \ |\eta |^{k+ \nu /2}
\eqno(3.40)$$
\noi we obtain as previously
$$|I_2| \leq C |w|_{k+ \nu /2}^2 \ |\varphi |_{\ell} \quad . \eqno(3.41)$$
\noi We finally consider the region $|\eta | \leq |\xi | \wedge |\xi - \eta |$. \par

\noi \underbar{Estimate of} $|I_3|$. Using (3.3) and
$$|\xi |^{2k} \ |\eta | \leq C |\xi |^{k+ \nu /2} \ |\xi - \eta |^{k - \nu/2 + 1 - \theta}\ |\eta
|^{\theta}$$
\noi with $0 \leq \theta \leq 1$ and decomposing $s = s_> + s_<$, we
estimate 
$$|I_3| \leq C |w|_{k+ \nu /2} \Big \{ \parallel (|\xi |^{k- \nu/2 + 1} f|\widehat{s}_<|) *
(f|\widehat{w}|)\parallel_2$$ 
$$+ \parallel (|\xi |^{k-\nu /2 + 1 - \theta} f|\widehat{s}_>|) * (|\xi |^{\theta}
f|\widehat{w}|)\parallel_2 \Big \} \ .\eqno(3.42)$$
\noi Using the Young inequality for the term in $s_<$ and Lemma 3.2 for the term in $s_>$, we
obtain  
$$|I_3| \leq C |w|_{k+ \nu /2} \Big \{ \parallel |\xi |^{k- \nu/2 + 1} f \ \widehat{s}_<
\parallel_1\  \parallel  f \widehat{w}\parallel_2 + \parallel |\xi |^{\ell + 1} f\widehat{s}_>
\parallel_2 \ \parallel |\xi |^{k+\nu /2} f\widehat{w}\parallel_2 \Big \}$$
\noi with $\ell > n/2 - \nu$ , $k + \nu /2 \geq \theta$ and $\ell \geq k - \nu /2 - \theta$. We
choose $\theta = 1 \wedge (k + \nu /2)$ and the last two conditions reduce to $\ell + 1 \geq k -
\nu /2$. Using (3.16), we obtain
$$|I_3| \leq C |w|_{k+ \nu /2}^2 \ |\varphi |_{\ell} \quad . \eqno(3.43)$$
\noi Collecting (3.32) (3.38) (3.39) (3.41) (3.43) yields (3.22). \par

The constants $C$ appearing in the proof of (3.22) are independent of $\rho$. This can be checked
explicitly on each case. More generally, it is a consequence of the following facts. We treat $f$
only through the inequalities (3.3) and (3.5) of Lemma 3.1, until we end up with integrals of the
type (3.30) with however $f$ appearing only through the product $f(\xi ) f(\xi - \eta )f(\eta )$
(possibly after adding one missing $f$ and using the fact that $f \geq 1$). It then follows from
(3.11) that an estimate of the type (3.22) of such a quantity for $\rho = 0$ implies the same
estimate for any $\rho > 0$ with the same constant. The same argument applies to all subsequent
estimates of Lemma 3.4 that contain $f$. Therefore from now on and in the same way as in the
proof of Lemma 3.3, we omit $f$ in the proofs.\\

(3.23). We estimate 
$$\parallel |\xi |^m (\widehat{s\cdot \nabla w})\parallel_2 \ \leq C \left \{ \parallel |\widehat{s}|
* (|\xi |^{m+1} |\widehat{w}|)\parallel_2 + \parallel (|\xi|^m |\widehat{s}|) * (|\xi |\ 
|\widehat{w}|)\parallel_2 \right \} \ .$$
\noi We decompose $s = s_< + s_>$, we estimate the contribution of $s_<$ by the Young inequality
and the contribution of $s_>$ by Lemma 3.2, thereby obtaining
$$\cdots \leq C \Big \{ \parallel \widehat{s}_< \parallel_1 \left (\parallel |\xi |^{m+1} \
\widehat{w}\parallel_2 \ + \ \parallel |\xi | \widehat{w} \parallel_2 \right ) + \ \parallel |\xi
|^{\ell + 1} \ \widehat{s}_> \parallel_2 \ \parallel |\xi |^k \widehat{w}\parallel_2 \Big \}$$
\noi under the condition stated on $k$, $\ell$, $m$. (3.23) then follows by (3.16). \\

(3.24). We estimate
$$\parallel |\xi |^m (\widehat{(\nabla \cdot s)w }) \parallel_2 \ \leq C \Big \{ \parallel
(|\xi|^{m+1}|\widehat{s}|) * |\widehat{w}|\parallel_2 + \parallel (|\xi| \ 
|\widehat{s}|) * (|\xi |^m |\widehat{w}|)\parallel_2 \Big \} \ .$$
\noi Proceeding as above, we obtain
$$\cdots \leq C \Big \{ \parallel \widehat{s}_< \parallel_1 \left (\parallel |\xi |^{m} \
\widehat{w}\parallel_2 \ + \ \parallel \widehat{w} \parallel_2 \right )
+ \ \parallel |\xi |^{\ell + 1} \ \widehat{s}_> \parallel_2 \ \parallel |\xi |^k
\widehat{w}\parallel_2 \Big \}$$
\noi from which (3.24) follows by the same argument as above. \\

(3.25) and (3.26). We decompose $s' = s'_> + s'_<$. The contribution of $s'_>$ to (3.25) and (3.26)
is estimated by (3.22) and (3.24) respectively, by replacing $w$ by $s'_>$ and $k$ by $\ell ' + 1$.
In order to complete the proof it remains to estimate $\parallel |\xi |_>^m (\widehat{s\cdot
\nabla s'_<})\parallel_2$ with $m = \ell ' + 1 - \nu /2$ for (3.25) and with general $m$ for
(3.26), namely to estimate the $L^2$ norm in $\xi$ of the integral
$$J = |\xi |_>^m \ f(\xi ) \int d \eta \ \widehat{s} (\xi - \eta ) \cdot \eta \ \widehat{s}'_< (\eta
)\quad .$$
\noi For that purpose, we decompose $s = s_{>1/2} + s_{<1/2}$ and correspondingly $J = J_> + J_<$.
In $J_>$ we have $|\eta | \leq 1 \leq 2 |\xi - \eta |$ and therefore $|\xi | \leq 3 |\xi - \eta
|$, so that
$$\begin{array}{ll} \parallel J_>\parallel_2 &\leq C \parallel (|\xi |^m \
|\widehat{s}_{>1/2}|) * (|\eta |\  |\widehat{s}'_<|)\parallel_2 \\ \\ &\leq C |s|_{\ell} \
\parallel \widehat{s}'_<\parallel_1 \end{array} \eqno(3.44)$$ \noi by the Young inequality,
provided $m \leq \ell + 1$, a condition which appears explicitly in (3.25) and which reduces
to $\ell \geq \ell ' - \nu/2$ in (3.26) for $m = \ell ' + 1 - \nu /2$. \par

In $J_<$, we have $|\xi - \eta | \leq 1/2$, $|\xi | \geq 1$ and $|\eta | \leq 1$, and therefore
$1/2 \leq |\eta | \leq 1$ and $|\xi | \leq 3/2$, so that
$$\parallel J_< \parallel_2 \ \leq C \parallel \widehat{s}_{<1/2} \parallel_1 \ \parallel
s'_{>1/2} \parallel_2 \quad . \eqno(3.45)$$

\noi (3.25) and (3.26) now follow from (3.22) (3.24), (3.44) (3.16) and (3.45). \\

(3.27) follows immediately from
$$\parallel |\xi |^{\ell_<} (\widehat{s\cdot s'}) \parallel_2 \ \leq C \left \{ \parallel
(|\xi|^{\ell<}\  |\widehat{s}|) * |\widehat{s}'| \parallel_2 \ + (s \leftrightarrow s')
\right \}$$ $$\leq C \left \{ \parallel |\xi |^{\ell_<} \ \widehat{s} \parallel_2 \
\parallel \widehat{s}'\parallel_1 \ + (s \leftrightarrow s') \right \}$$ \noi and from
(3.16) (3.17). \\

(3.28) follows from (3.14) either directly for $\beta \leq 0$, or through the inequality
$$\parallel |\xi |^{m} \ \widehat{g_0(w_1}w_2)_> \parallel_2 \ \leq C \left \{ \parallel
(|\xi |^{\beta} \ |\widehat{w}_1|) * (|\widehat{w}_2|)\parallel_2 \ + (1 \leftrightarrow 2)
\right \}$$ \noi for $\beta \geq 0$. \\

(3.29) follows from (3.15) with $m = \beta$.\par
\hfill $\sq$\par

We now explain the origin of the derivative loss in the system (2.11)-(2.12) for\break
\noindent $\lambda \equiv \mu - n+2 >0$ and the mechanism by which that loss is overcome
through the use of the spaces ${\cal X}_{\rho , loc}^{k, \ell}$ defined by (3.13). If we try
to solve the system (2.11)-(2.12) by the energy method in a space like ${\cal C}(I, H^k
\oplus \dot{H}^{\ell + 2}$), we have to estimate in particular
$$\left \{ \begin{array}{l}
\partial_t \parallel \partial^k w \parallel_2^2 \ = 2 {\rm Re} <\partial^kw,
\partial^k\ \partial_t \ w > \\ \\
\partial_t \parallel \partial^{\ell + 2} \varphi \parallel_2^2 \ = 2 <\partial^{\ell + 2}
\ \varphi , \partial^{\ell + 2} \ \partial_t \ \varphi > \quad . \end{array} \right .\eqno(3.46)$$
 
\noi The term with $\Delta \phi$ from $\partial_tw$ forces us to apply $k+2$ derivatives to
$\varphi$ and requires therefore $\ell \geq k$, while the term with $g_0$ from $\partial_t\varphi$
forces us to apply $\ell + 2$ derivatives to $g_0$ or equivalently $\ell + \lambda$ derivatives to
$|w|^2$ and requires therefore $k \geq \ell + \lambda$. The terms with $\nabla \varphi \cdot \nabla
w$ from $\partial_tw$ and with $|\nabla \varphi |^2$ from $\partial_t\varphi$ can be handled
essentially under the same assumptions, possibly after an integration by parts. The method
therefore applies only if $\lambda \leq 0$, namely $\mu \leq n - 2$, which is the case treated in I
and II. \par

If we try instead to solve the same problem in the space ${\cal X}_{\rho, loc}^{k,\ell}$ with time
dependent $\rho$, we have by (3.12)
$$\left \{ \begin{array}{l}
\partial_t |w|_k^2 = 2 \rho ' |w|_{k+ \nu /2}^2 + 2 {\rm Re} \ <w, \partial_tw>_k \\ \\
\partial_t |\varphi |_{\ell}^2 = 2 \rho ' |\varphi |_{\ell + \nu /2}^2 + 2 <\varphi, \partial_t
\varphi >_{\ell} \end{array} \right .\eqno(3.47)$$

\noi where $<\cdot , \cdot >_k$ and $<\cdot , \cdot >_{\ell}$ denote the scalar products in
$K_{\rho}^k$ and $Y_{\rho}^{\ell}$. If $\rho '$ has a favourable sign, namely if $\rho$ decreases
away from the initial time, the terms containing $\rho '$ provide a control of the norm in
$L^2(K_{\rho}^{k+\nu /2} \oplus Y_{\rho}^{\ell + \nu /2})$, and it suffices to control the scalar
products at most quadratically in terms of the norms in $K_{\rho}^{k+\nu /2}$ and
$Y_{\rho}^{\ell + \nu /2}$. In the term with $\Delta \phi$ from $\partial_tw$, after
distributing the function $f$ with the help of (3.3) and shifting $\nu /2$ derivatives on the
first vector of the scalar product, it suffices to apply $k - \nu /2$ derivatives on $\Delta
\phi$ while one is allowed to use the norm $|\varphi |_{\ell + \nu /2}$. This requires only
$\ell \geq k - \nu$. Similarly in the term with $g_0$ from $\partial_t \varphi$, it suffices to
apply $\ell + 2 - \nu /2$ derivatives to $g_0$ or equivalently $\ell + \lambda - \nu /2$
derivatives to $|w|^2$, while one is allowed to use $|w|_{k + \nu /2}$. This requires only $k
\geq \ell + \lambda - \nu$. The two conditions on $(k, \ell )$ are compatible provided $\lambda
\leq 2 \nu$, which allows for $\mu \leq n$ under that condition. It remains to estimate the
terms $\nabla \varphi \cdot \nabla w$ from $\partial_tw$ and $|\nabla \varphi |^2$ from
$\partial_t \varphi$, which in the Sobolev case requires the integration by part of one
derivative. Here however, by the same argument as above, it suffices to integrate by parts $1 -
\nu$ derivative. Now it turns out that integration by parts of $1 - \nu$ derivative is exactly what
is allowed by the inequality (3.5), which is exploited through (3.31) to derive the estimates (3.22)
(3.25) where that integration by parts occurs. Actually the inequality (3.5) is optimal in the
dangerous part $|\xi - \eta | \ll |\xi | \sim |\eta |$ of the region where it is used, and more
precisely when $|\eta | \to \infty$ for fixed $|\xi - \eta |$. The conditions $\ell \geq k - \nu$
and $k \geq \ell + \lambda - \nu$ and therefore their consequence $\lambda \leq 2 \nu$ will appear
from the next lemma onward as the most important part of the condition (3.48) and will propagate
throughout this paper (except in Section 5) up to the main and final results of Propositions 6.5 and
7.5. \par

We now exploit Lemma 3.4 to derive energy like estimates for the solutions of the auxiliary
system (2.11)-(2.12). In the following three lemmas, $I$ is an interval contained in $[1,
\infty )$, $\rho$ is a nonnegative continuous and piecewise ${\cal C}^1$ function defined in
$I$. We shall be interested in solutions $(w, \varphi )$ in spaces of the type ${\cal
X}_{\rho , loc}^{k , \ell}(I)$ for suitable values of $k$ and $\ell$. The estimates will
hold in integrated form in any compact subinterval of $I$ under the available regularity,
but will be stated in differential form for brevity. \\

\noi {\bf Lemma 3.5.} {\it Let $k$, $\ell$ satisfy
$$\left \{ \begin{array}{l} \ell > n/2 - \nu \quad , \quad k \geq \nu /2 \quad ,
\quad \ell \geq k - \nu \quad , \\ \\ k \geq \ell + \lambda - \nu \quad , \quad 2
k > \ell + \lambda - \nu + n/2 \quad ,\end{array}\right . \eqno(3.48)$$
\noi where $\lambda = \mu - n + 2$. \par

Let $(w, \varphi ) \in {\cal X}_{\rho , loc}^{k,\ell}(I)$ be a solution of the system
(2.11)-(2.12). Then the following estimates hold~:} $$\left | \partial_t|w|_k^2 - 2 \rho '
|w|_{k + \nu /2}^2 \right | \leq C \ t^{-2} \left \{ |w|_{k+ \nu /2}^2 \ |\varphi |_{\ell} +
|w|_{k + \nu /2} \ |\varphi|_{\ell + \nu /2} \ |w|_k \right \} \ , \eqno(3.49)$$ 
$$\left | \partial_t|\varphi|_{\ell}^2 - 2 \rho ' |\varphi|_{\ell + \nu /2}^2 \right | \leq C \
t^{-2} |\varphi|_{\ell + \nu /2}^2 \ |\varphi |_{\ell} + C \ t^{-\gamma}|\varphi|_{\ell + \nu /2}
\ |w|_{k + \nu /2}\  |w|_k \ .\eqno(3.50)$$

\vskip 3 truemm      
\noi {\bf Proof.} We use the same regularization as in the proof of Lemma 3.4. From (3.12) we obtain
$$\partial_t |J_{\varepsilon} w|_k^2 - 2 \rho ' | J_{\varepsilon} w|_{k + \nu /2}^2
= 2 {\rm Re} \left ( < j_{\varepsilon}|\xi |^k f \widehat{w}_>, j_{\varepsilon} |\xi |^k f
\partial_t \widehat{w}_> \ > + < j_{\varepsilon} f \widehat{w}_<, j_{\varepsilon} f
\partial_t \widehat{w}_< \ > \right ) \ . \eqno(3.51)$$
\noi We substitute $\partial_t w$ from (2.11) and we estimate the various terms successively. The
term with $\Delta w$ does not contribute. The term $s\cdot \nabla w$ is estimated by (3.22) under
the conditions
$$\ell > n/2 - \nu \qquad , \quad k \geq \nu /2 \qquad , \quad \ell + 1 \geq k - \nu /2 \quad .$$
\noi The term $(\nabla \cdot s)w$ is estimated by (3.24) with $m = k - \nu /2$ or $m = 0$, and
$\ell$ replaced by $\ell + \nu /2$, under the conditions
$$\ell > n/2 - \nu \qquad , \qquad \ell \geq k - \nu \quad .$$

Substituting those two estimates into (3.51), integrating over time and taking the limit
$\varepsilon \to 0$ yields the integrated form of (3.49). The required conditions on $k$, $\ell$
are implied by (3.48). \par

Similarly, we obtain from (3.12)
$$\partial_t |J_{\varepsilon} \ \varphi |_{\ell}^2 - 2 \rho ' | J_{\varepsilon} \ \varphi |_{\ell +
\nu /2}^2 = 2 < j_{\varepsilon} |\xi |^{\ell + 1} \ f \ \widehat{s}_>, j_{\varepsilon}\ |\xi
|^{\ell + 1} \ f \ \partial_t \ \widehat{s}_>\  > + 2 < j_{\varepsilon} |\xi |^{\ell_<}\ 
f \  \widehat{\varphi}_<, j_{\varepsilon} |\xi|^{\ell_<} \ f \ \partial_t \ \widehat{\varphi}_<
\ > \ . \eqno(3.52)$$ 
\noi We substitute $\partial_ts$ and $\partial_t \varphi$ from (2.27)
(2.12)  into (3.52) and we estimate the various terms successively. The term $s\cdot \nabla s$ from
$\partial_ts$ is estimated by (3.25) with $s' = s$, $\ell ' = \ell$ under the conditions
$$\ell > n/2 - \nu \qquad , \quad \ell + 1 \geq \nu /2 \quad .$$
\noi The term $\nabla g_0$ from $\partial_ts$ is estimated by (3.28) with $m = \ell + 2 - \nu /2$,
$w_1 = w_2 = w$, $k_1 (=k'_2) = k + \nu /2$ , $k_2 (=k'_1) = k$ , under the conditions
$$2k > \ell + \lambda - \nu + n/2 \qquad , \qquad k \geq \ell + \lambda - \nu \quad .$$
\noi The terms $|s|^2$ and $g_0$ from $\partial_t \varphi_<$ are estimated by (3.27) and (3.29)
respectively. Using those estimates yields (3.50) in the same way as above. The required conditions
on $k$, $\ell$ are implied by (3.48). In particular the condition $\ell + 1 \geq \nu /2$ follows
from $\ell > n/2 - \nu$. \par
\hfill $\sq$ \par

The next lemma is a regularity result. \\

\noi {\bf Lemma 3.6.} {\it Let $k$, $\ell$ satisfy (3.48) and let $\bar{k}$, $\bar{\ell}$ satisfy
$$\bar{k} - k = \bar{\ell} - \ell \geq 0 \quad . \eqno(3.53)$$
\noi Let $(w, \varphi ) \in {\cal X}_{\rho , loc}^{\bar{k},\bar{\ell}} (I)$ be a solution of the
system (2.11)-(2.12). Then the following estimates hold~:}
$$\left |\partial_t|w|_{\bar{k}}^2 - 2 \rho ' |w|_{\bar{k} + \nu /2}^2 \right | \leq C \ t^{-2} \left
\{ |w|_{\bar{k}+ \nu /2}^2 \ |\varphi|_{\ell} + |w|_{\bar{k}+ \nu /2}\  |\varphi |_{\bar{\ell} + \nu
/2} \ |w|_k \right \} \ , \eqno(3.54)$$
$$\left |\partial_t|\varphi|_{\bar{\ell}}^2 - 2 \rho ' |\varphi|_{\bar{\ell} + \nu /2}^2 \right |
\leq C \ t^{-2} \ |\varphi|_{\bar{\ell}+ \nu /2}^2\  |\varphi|_{\ell} + C \ t^{-\gamma}
|\varphi|_{\bar{\ell}+ \nu /2}\  |w |_{\bar{k} + \nu /2}\  |w|_k  \ . \eqno(3.55)$$

\vskip 3 truemm 
\noi {\bf Proof.} The proof follows the same pattern as that of Lemma 3.5 and we concentrate on
the differences, which bear on the estimates connected with Lemma 3.4, omitting the regularization
for brevity. We estimate only the contribution of the high $|\xi |$ region, since that of the low
$|\xi |$ region is already estimated by Lemma 3.5. \par

We first estimate $\partial_t |w| {2 \over k}$, starting from (3.51) with $k$ replaced by
$\bar{k}$ and with $J_{\varepsilon}$ omitted and we estimate successively the contribution of the
various terms form (2.11). \par

The contribution of $s\cdot \nabla w$, namely
$$2 {\rm Re} \ < |\xi  |^{\bar{k}} \ f \ \widehat{w} , |\xi |^{\bar{k}} \ f(\widehat{s \cdot
\nabla w}) >$$
\noi is estimated in the same way as in the proof of (3.22). The contributions $I_1$ and $I_2$ of
the regions $|\xi - \eta | \leq |\xi | \wedge |\eta |$ and $| \xi | \leq | \xi - \eta | \wedge
|\eta |$ are estimated as in the latter with $k$ replaced by $\bar{k}$ under the conditions
$\bar{k} \geq \nu /2$ , $\ell > n/2 - \nu$. \par

The contribution of the region $|\eta | \leq |\xi | \wedge |\xi - \eta |$ is estimated by (3.42)
with $k$ replaced by $\bar{k}$ and with $\theta = 0$, or equivalently
$$|I_3| \leq C |w|_{\bar{k} + \nu /2} \Big \{ \parallel f\widehat{s}_< \parallel_1 \ \parallel f
\widehat{w} \parallel_2 \ + \parallel (|\xi|^{\bar{k} - \nu /2 + 1} f |\widehat{s}_>|) *
(f|\widehat{w}|) \parallel_2 \big \} \ .$$
\noi The last norm is then estimated by Lemma 3.2, thereby yielding
$$|I_3| \leq C |w|_{\bar{k}+ \nu /2} \ |\varphi |_{\bar{\ell} + \nu /2} \ |w|_k$$
\noi under the conditions
$$k + \bar{\ell} > \bar{k} + n/2 - \nu \qquad , \qquad \bar{\ell} \geq \bar{k} - \nu$$
\noi which for $\bar{k} - k = \bar{\ell} - \ell$ reduce to 
$$\ell > n/2 - \nu \qquad , \qquad \ell \geq k - \nu \ .$$

The contribution of $(\nabla \cdot s)w$ is estimated by 
$$|w|_{\bar{k}+ \nu /2} \parallel |\xi |^{\bar{k} - \nu /2} \ f(\widehat{(\nabla \cdot
s)w})\parallel_2$$
\noi and subsequently in a way similar to the proof of (3.24) with $m = \bar{k} - \nu /2$. Using
the elementary inequality
$$|\xi |^m |\xi - \eta | \leq 2^m \left ( |\xi - \eta |^{m+1} + |\eta |^{m + \theta} |\xi - \eta
|^{1 - \theta } \right ) \eqno(3.56)$$
\noi valid for $m \geq 0$ and $0 \leq \theta \leq 1$, we estimate the last norm by
$$\cdots C \Big \{ \parallel (|\xi |^{m+1} f|\widehat{s}|) * (f|\widehat{w}|)\parallel_2\  +
\parallel (|\xi |^{1 - \theta } f|\widehat{s}|) * (|\xi |^{m + \theta} f|\widehat{w}|)\parallel_2
\Big \} \ .$$
\noi We then estimate the first term by Lemma 3.2 and the second term with $\theta = \nu$ by the
Young inequality and (3.16) (3.17), thereby obtaining 
$$\cdots \leq C \left ( |\varphi |_{\bar{\ell} + \nu /2} \ |w|_k + |\varphi |_{\ell} \
|w|_{\bar{k} + \nu /2} \right )$$
\noi under the same conditions as before, namely
$$k + \bar{\ell} > \bar{k} + n/2 - \nu \quad , \quad \bar{\ell} \geq \bar{k} - \nu \ .$$
\noi This completes the proof of (3.54). \par

We next estimate $\partial_t |\varphi |_{\bar{\ell}}^2$ by substituting similarly the various
terms from (2.12) (2.27) into the right-hand side of (3.52) with $\ell$ replaced by $\bar{\ell}$
and $J_{\varepsilon}$ omitted. The contribution of $s\cdot \nabla s$ to the high $|\xi |$ part of
the norm, namely
$$< |\xi |^{\bar{\ell} + 1} \ f\widehat{s}_> , |\xi |^{\bar{\ell} + 1} \ f(\widehat{s\cdot
\nabla s})_> \ >$$
\noi is estimated by modifying the proof of (3.25) along the same lines as that of (3.22) above.
\par

The contribution of $g_0(w,w)_>$ is estimated by
$$|\varphi |_{\bar{\ell} + \nu /2} \parallel |\xi |^{\bar{\ell} + 1 - \nu /2} \ f\nabla
\widehat{g_0(w,}w)_>\parallel_2$$
\noi followed by (3.28) with $m = \bar{\ell} + 2 - \nu /2$ , $k_1 (=k'_2) = \bar{k} + \nu /2$ , $k_2
(=k'_1) = k$ under the conditions
$$k + \bar{k} > \bar{\ell} + \lambda - \nu + n/2 \qquad , \quad \bar{k} \geq \bar{\ell} + \lambda -
\nu$$
\noi which reduce to the same conditions with $(\bar{k} , \bar{\ell})$ replaced by $(k, \ell)$
under the condition $\bar{k} - k = \bar{\ell} - \ell$. \par

This completes the proof of (3.55). \par
\hfill $\sq$ \par

We next estimate the difference between two solutions of the system (2.11)-(2.12). \\

\noi {\bf Lemma 3.7.} {\it Let $k$, $\ell$ satisfy (3.48) and let $k'$, $\ell '$ satisfy
$$k' \geq \nu /2 \quad , \quad k - k' = \ell - \ell ' \geq 1 - \nu \ . \eqno(3.57)$$
\noi Let $(w_1, \varphi_1)$ and $(w_2, \varphi_2) \in {\cal X}_{\rho , loc}^{k , \ell} (I)$ be two
solutions of the system (2.11)-(2.12) and let $w_{\pm} = w_1 \pm w_2$, $\varphi_{\pm} =
\varphi_1 \pm \varphi_2$. Then the following estimates hold~:}
$$\left | \partial_t |w_-|_{k'}^2 - 2 \rho ' |w_-|_{k' + \nu /2}^2 \right | \leq C \ t^{-2} \Big \{
|w_-|_{k' + \nu /2}^2 \ |\varphi_+|_{\ell}$$
$$+ |w_-|_{k' + \nu /2} \left ( |w_+|_{k + \nu /2} \  |\varphi_-|_{\ell '} + |\varphi_+|_{\ell + \nu
/2} \ |w_-|_{k'} + |\varphi_- |_{\ell ' + \nu /2} \ |w_+|_k \right ) \Big \} \ , \eqno(3.58)$$
$$\left | \partial_t |\varphi_-|_{\ell '}^2 - 2 \rho ' |\varphi_-|_{\ell' + \nu /2}^2 \right | \leq C
\ t^{-2} \Big \{ |\varphi_-|_{\ell ' + \nu /2}^2 \ |\varphi_+|_{\ell} + |\varphi_-|_{\ell ' + \nu
/2} \ |\varphi_-|_{\ell '}\ |\varphi_+|_{\ell + \nu /2} \Big \}$$ 
$$ + C\ t^{-\gamma} |\varphi_-|_{\ell ' + \nu /2} \Big \{ 
|w_+|_{k + \nu /2} \  |w_-|_{k '} + |w_-|_{k' + \nu /2} \ |w_+|_{k} \Big \} \ . \eqno(3.59)$$

\vskip 3 truemm

\noi {\bf Proof.} The proof follows the same pattern as that of Lemma 3.5, using the estimates of
Lemma 3.4. We omit again the regularization for brevity. The equations satisfied by $(w_-,
\varphi_-)$ are
$$\partial_t w_- = i(2t^2)^{-1} \Delta w_- + (2t)^{-2} \left \{ 2 s_+ \cdot \nabla w_- + 2s_- \cdot
\nabla w_+ + (\nabla \cdot s_+) w_- + (\nabla \cdot s_-)w_+ \right \} \ , \eqno(3.60)$$
$$\partial_t \varphi_- = (2t^2)^{-1} (s_+ \cdot s_-) + t^{-\gamma} g_0 (w_+ , w_-) \ , \eqno(3.61)$$
\noi and in the same way as in the proof of Lemma 3.5, we shall also use the equation for $s_-$
obtained by taking the gradient of (3.61), namely 
$$\partial_t s_- = (2t^2)^{-1} \left ( s_+ \cdot \nabla s_- + s_- \cdot \nabla s_+ \right ) +
t^{-\gamma} \nabla g_0(w_+, w_-) \ . \eqno(3.62)$$
\noi From (3.12), we obtain
$$\left | \partial_t |w_-]_{k'}^2 - 2 \rho ' |w_-|_{k'+ \nu /2}^2 \right | = 2 {\rm Re} \Big (
< |\xi |^{k'} f \widehat{w}_{->}, |\xi |^{k'} f \partial_t \widehat{w}_{->} \ >
+ < f \widehat{w}_{-<} , f \ \partial_t \ \widehat{w}_{-<} \ > \Big ) \ . \eqno(3.63)$$
\noi We substitute $\partial_t w_-$ from (3.60) into (3.63) and we estimate the various terms
successively. The term $\Delta w_-$ does not contribute. We consider only the contribution of the
high $|\xi |$ region. The low $|\xi |$ region is treated in a similar and simpler way. \\

\noi \underbar{Term with} $s_+\cdot \nabla w_-$. We apply (3.22) with $k$ replaced by $k'$,
thereby obtaining a contribution
$$|w_-|_{k'+ \nu /2}^2 \ |\varphi _+ |_{\ell}$$
\noi under conditions which follows from (3.48) and from $\nu /2 \leq k' \leq k$. \\

\noi \underbar{Term with} $s_-\cdot \nabla w_+$. We apply (3.23) with $m = k' - \nu /2$ , $\ell$
replaced by $\ell '$ and $k$ replaced by $k + \nu /2$, thereby obtaining a contribution
$$|w_-|_{k' + \nu /2} \ |w_+|_{k + \nu /2} \ |\varphi_-|_{\ell '}$$
\noi under the conditions
$$k + \ell ' > k' + n/2 - \nu \quad , \ k + \nu \geq k' + 1 \quad , \ \ell + 1 \geq k' - \nu /2 \
,$$
\noi which follow from (3.48) (3.57). \\

\noi \underbar{Term with} $(\nabla\cdot  s_+) w_-$. We apply (3.24) with $m = k' - \nu /2$, $k$
replaced by $k'$ and $\ell$ replaced by $\ell + \nu /2$, thereby obtaining a contribution
$$|w_-|_{k' + \nu /2} \ |\varphi_+|_{\ell + \nu /2} \ |w_-|_{k'}$$
\noi under the conditions
$$\ell > n/2 - \nu \qquad , \quad \ell \geq k' - \nu \ ,$$
\noi which follow from (3.48) and from $k' \leq k$. \\

\noi \underbar{Term with} $(\nabla\cdot  s_-) w_+$. We apply (3.24) with $m = k' - \nu /2$, and
$\ell$ replaced by $\ell ' + \nu /2$, thereby obtaining a contribution
$$|w_-|_{k'+ \nu /2} \ |\varphi_-|_{\ell ' + \nu /2} \ |w_+|_k$$
\noi under the conditions
$$k + \ell ' > k' + n/2 - \nu \qquad , \quad \ell ' \geq k' - \nu \qquad , \quad k \geq k' - \nu /2
\ ,$$
\noi which follow from (3.48) (3.57). \par

Collecting the previous four estimates together with the contribution of the low $|\xi |$ region,
yields (3.58). \par

We now turn to the estimate of $\varphi_-$. From (3.12) we obtain
$$\partial_t |\varphi_-|_{\ell '}^2 - 2 \rho ' |\varphi_-|_{\ell ' + \nu /2}^2 = 2 < |\xi
|^{\ell ' + 1} \ f\ \widehat{s}_{->}, |\xi |^{\ell ' + 1} \ f \ \partial_t \ \widehat{s}_{->}
\ > + 2 < |\xi |^{\ell_<} \ f \ \widehat{\varphi}_{-<}, |\xi |^{\ell_<} \ f \
\widehat{\varphi}_{-<} > \ . \eqno(3.64)$$
\noi We substitute $\partial_t s_-$ and $\partial_t \varphi_-$ from (3.62) and (3.61) into (3.64)
and we estimate the various terms successively. \\

\noi \underbar{Term with} $s_+\cdot \nabla s_-$. We apply (3.25) with $s' = s_-$ and obtain a
contribution
$$|\varphi_-|_{\ell ' + \nu /2}^2 \ |\varphi_+|_{\ell}$$
\noi under the conditions
$$\ell > n/2 - \nu \quad , \ \ell ' + 1 \geq \nu /2 \quad , \ \ell \geq \ell ' - \nu /2 \ ,$$
\noi which follow from (3.48), from $\ell ' + 1 \geq k' + 1 - \nu \geq 1 - \nu /2 \geq \nu /2$ and
from $\ell \geq \ell '$. \\

\noi \underbar{Term with} $s_-\cdot \nabla s_+$. We apply (3.26) with $m = \ell ' + 1 - \nu /2$ ,
$\ell$ replaced by $\ell '$ and $\ell '$ replaced by $\ell + \nu /2$, thereby obtaining a
contribution
$$|\varphi_-|_{\ell ' + \nu /2} \ |\varphi_- |_{\ell '} \ |\varphi_+ |_{\ell + \nu /2}$$
\noi under the conditions
$$\ell > n/2 - \nu \qquad , \quad \ell \geq \ell ' + 1 - \nu \ ,$$
\noi which follow from (3.48) (3.57). \\

\noi \underbar{Term with} $\nabla g_0(w_+,w_-)_>$. We apply (3.28) with $m = \ell ' + 1 - \nu /2$,
$w_1 = w_+$, $w_2 = w_-$, $k_1 = k + \nu /2$ , $k_2 = k'$ , $k'_1 = k$ , $k'_2 = k' + \nu /2$,
thereby obtaining a contribution
$$|\varphi_-|_{\ell ' + \nu /2} \left \{ |w_+|_{k + \nu /2} \ |w_-|_{k'} + |w_-|_{k' + \nu /2} \
|w_+|_k \right \}$$
\noi under the conditions
$$k + k' + \nu > \ell ' + \mu - n + 2 + n/2 \quad , \ k' + \nu \geq \ell ' + \mu - n + 2$$
\noi which follow from (3.48) (3.57). \par

The terms with $(s_+, s_-)_<$ and $g_0(w_+w_-)_<$ are treated by the use of (3.27) and (3.29) 	as in
the proof of Lemma 3.5. \par

Collecting the previous estimates yields (3.59). \par
\hfill $\sq$ \par

We conclude this section by introducing a number of estimating functions of time
ge\-ne\-ra\-li\-zing those introduced in Section II.3 and by deriving a number of estimates for
them. Those functions will be defined in terms of the derivative $h'_0$ of a given function $h_0$ on
which we make the following assumptions
$$h_0 \in {\cal C}^1([1, \infty ) , {I\hskip-1truemm R}^+)\ , \ h'_0 \geq 0\ , \ t^{-2}h_0(t) \in
L^1([1, \infty ))\ , \ t^{-1} h'_0(t) \in L^1([1, \infty )) \ . \eqno(3.65)$$
\noi From the relation
$$t_2^{-1} h_0(t_2) - t_1^{-1} \ h_0(t_1) = \int_{t_1}^{t_2} dt \ t^{-1} \ h'_0  (t) -
\int_{t_1}^{t_2} dt \ t^{-2} \ h_0(t) \eqno(3.66)$$
\noi it follows that the last condition on $h_0$ in (3.65) can be replaced by the condition that
$t^{-1} h_0(t)$ tends to zero when $t \to \infty$. A typical example for $h'_0$ is that
considered in Section II.3, namely $h'_0(t) = t^{-\gamma}$. \par

The first and basic estimating function $h$ is defined by
$$h(t) = \int_1^{\infty} dt_1 (t \vee t_1)^{-1} \ h'_0 (t_1) \eqno(3.67)$$
\noi from which it follows that $h(t)$ is decreasing in $t$ and tends to zero when $t \to
\infty$, while $t h(t)$ is increasing in $t$. We next define for any $m \geq 0$
$$N_m(t) = \int_1^t dt_1 \ h'_0(t_1) \ h^m(t_1) \ , \eqno(3.68)$$
$$Q_m (t) = \int_1^{\infty} dt_1 (t \vee t_1)^{-1} \ h'_0(t_1) \ h^m (t_1) \ , \eqno(3.69)$$
\noi where the integral in (3.69) is convergent since $t^{-1}h'_0(t) \in L^1([1, \infty ))$ and
since $h(t)$ is decreasing in $t$. It follows from (3.68) that $N_m(t)$ is increasing in $t$ and
from (3.69) that $Q_m(t)$ is decreasing in $t$ and tends to zero when $t \to \infty$, while $t
Q_m(t)$ is increasing in $t$, so that $Q_m(t) \geq t^{-1} Q_m(1)$. Moreover, for any nonnegative
integers $i$ and $j$ 
$$N_{i+j}(t) \leq h(1)^i\ N_j(t) \leq h(1)^{i+j}\ N_0(t) \ , \eqno(3.70)$$
$$Q_{i+j}(t) \leq h(1)^i\ Q_j(t) \leq h(1)^{i+j}\ h(t) \leq h(1)^{i+j+1} \ . \eqno(3.71)$$ 
\noi Clearly $N_0(t) = h_0(t) - h_0(1)$ and $Q_0 = h$. It will be convenient to introduce the
notation $Q_{-1} = 1$. \par

Finally we set
$$P_m(t) = \int_1^{\infty} dt_1 \ h(t \vee t_1) \ h'_0(t_1) \ h^m (t_1) \eqno(3.72)$$
\noi which is well defined provided
$$P_m(1) = \int_1^{\infty} dt_1 \ h'_0(t_1) \ h^{m+1} (t_1) < \infty \ . \eqno(3.73)$$
\noi It follows from (3.72) that $P_m(t)$ is decreasing in $t$ and tends to zero when $t \to
\infty$ while $h(t)^{-1} P_m(t)$ is increasing in $t$, so that
$$P_m(1) \ h(t) \leq h(1) \ P_m(t) \leq h(1) \ P_m (1) \ . \eqno(3.74)$$

We now collect a number of identities and inequalities satisfied by the previous estimating
functions.\\

\noi {\bf Lemma 3.8} {\it Let $i$, $j$ and $m$ be nonnegative integers, let $1 \leq a \leq b$ and
$t \geq 1$. Then the following identities and inequalities hold~:}
$$\int_t^{\infty} dt_1 \ t_1^{-2} \ N_m(t_1) = Q_m(t_1) \eqno(3.75)$$
$$\int_1^t dt_1 \ t_1^{-2} \ N_0 (t_1) \ N_m (t_1) = N_{m+1}(t) - h(t) \ N_m(t) \leq N_{m+1}(t)
\eqno(3.76)$$
$$\int_t^{\infty} dt_1 \ t_1^{-2} \ N_0(t_1) \ N_m (t_1) = P_m (t) \quad \hbox{\it if (3.73)
holds} \eqno(3.77)$$
$$N_i(t) \ N_j(t) \leq N_0(t) \ N_{i+j}(t) \eqno(3.78)$$
$$N_i(t) \ Q_j(t) \leq h(t) \ N_{i+j}(t) \leq N_{i+j+1} (t) \eqno(3.79)$$
$$Q_i(t) \ Q_j(t) \leq h (t) \ Q_{i+j}(t) \leq 2 Q_{i+j+1} (t) \eqno(3.80)$$
$$\int_t^{\infty} dt_1 \ h'_0(t_1) \ h(t_1) \ Q_{m-1}(t_1) \leq \int_t^{\infty} dt_1 \ h'_0
(t_1) \ Q_m(t_1) \eqno(3.81)$$
$$\int_t^{\infty} dt_1 \ h'_0(t_1) \ Q_m(t_1) \leq P_m(t_1) \quad  \hbox{\it if (3.73)
holds} \eqno(3.82)$$
$$\int_1^t dt_1 \ h'_0(t_1) \ h(t_1) \ Q_{m-1}(t_1) \leq N_{m+1} (t_1) \eqno(3.83)$$
$$\int_1^t dt_1 \ h'_0(t_1) \  Q_m(t_1) \leq N_{m+1}(t) \eqno(3.84)$$
$$\int_a^b dt \ h'_0(t) \ Q_m(t) \leq Q_m(a) \left ( h_0(b) - h_0(a) \right ) \eqno(3.85)$$
$$\int_a^b dt \ h'_0(t) \ h(t) \ Q_{m-1}(t) \leq 2Q_m(a) \left ( h_0(b) - h_0(a) \right ) \ . 
\eqno(3.86)$$
\vskip 3 truemm

\noi {\bf Proof.} This lemma is a generalization of Lemma II.3.6 and most of the proofs are
obtained by manipulations of integrals similar to those of the corresponding integrals in Lemma
II.3.6, after the replacement of $t^{-\gamma}$ by $h'_0(t)$. This applies to (3.75) (3.76) (3.77)
(3.81) (3.82) (3.83) (3.84) and (3.80b). (The latter is the generalization of (II.3.49)). The
estimates (3.78) and (3.80a) follow from the H\"older inequality. The estimate (3.79a) is the
pointwise counterpart of (II.3.39) and is proved in the same way. The estimate (3.79b) follows from
the decrease of $h$, (3.85) follows from the decrease of $Q_m$, and (3.86) follows from (3.80b) and
from (3.85). \par
\hfill $\sq$\par

\section{Cauchy problem and preliminary asymptotics for the auxiliary system}
\hspace*{\parindent} In this section, we study the existence of solutions in a neighborhood of
infinity in time for the auxiliary system (2.11)-(2.12) and we derive some preliminary
results on the asymptotic behaviour in time of the solutions thereby obtained. We first
introduce some notation.  \par

We choose once for all a strictly positive ${\cal C}^1$ function defined in $[1, \infty )$,
which we call $|\rho '|$ for reasons that will soon become obvious. We assume that $|\rho '|$
satisfies the following properties \par

(i) $|\rho ' | \in L^1([1, \infty ))$, \par
(ii) the function $t^{-\gamma} |\rho '|^{-1}$ is nondecreasing (and therefore $|\rho '|$ is
nonincreasing), and the function $t^{-2}|\rho '|^{-1}$ is nonincreasing and tends to zero at
infinity. \par

Typical examples of suitable functions $|\rho '|$ are $t^{-1 - \varepsilon}$ for $\varepsilon$
sufficiently small, depending on $\gamma$, or $t^{-1}(\alpha + \ell n t)^{-\alpha}$ for $\alpha >
1$. It will be useful to keep those examples in mind in order to understand the time decay
implied by the subsequent estimates. \par

Let now $I \subset [1, \infty )$ be an interval, possibly unbounded, and let $t_0 \in I$ (or
$t_0 = \infty $ if $I$ is unbounded). We define a function $\rho$ in $I$ by
$$\rho (t) = \rho (t_0) - \left | \int_{t_0}^t dt_1 |\rho ' (t_1)| \right | \eqno(4.1)$$
\noi so that $\rho$ is increasing (resp. decreasing) for $t \leq t_0$ (resp. $t \geq t_0$) and has
$|\rho '|$ as the absolute value of its derivative, which justifies the notation. We take $\rho
(t_0)$ sufficiently large so that $\rho$ is nonnegative in $I$. All subsequent estimates will be
independent of $\rho (t_0)$ (they will however depend on $|\rho '|$). The previous choice of
$\rho$ will be used in this section without further comment unless otherwise stated. We now define
the global version of the fundamental spaces, corresponding to the local version (3.13). We define 
$${\cal X}_{\rho}^{k,\ell} (I) \equiv \left ({\cal C} \cap L^{\infty}\right ) \left (I, K_{\rho}^k
\oplus Y_{\rho}^{\ell}\right ) \cap L_{|\rho '|}^2 \left ( I, K_{\rho}^{k + \nu /2} \oplus
Y_{\rho}^{\ell + \nu /2} \right ) \eqno(4.2)$$
\noi where $L^2_{|\rho '|}$ denotes weighted $L^2$ in time with weight $|\rho '|$. More precisely,
especially as regards continuity, $(w, \varphi ) \in {\cal X}_{\rho}^{k, \ell} (I)$ is understood
to mean that    
$$(F^{-1} f \widehat{w}, F^{-1}f\widehat{\varphi}) \in {\cal X}_{0}^{k,\ell} (I) \equiv \left ({\cal
C} \cap L^{\infty}\right ) \left (I, K_{0}^k \oplus Y_{0}^{\ell}\right ) \cap L_{|\rho '|}^2 \left (
I, K_{0}^{k + \nu /2} \oplus Y_{0}^{\ell + \nu /2} \right ) \ .\eqno(4.3)$$
\noi Note that in (4.3) we keep the weight $|\rho '|$ in the $L^2$ part. \par

The norm of $(w, \varphi )$ in ${\cal X}_{\rho}^{k,\ell}$ is made up of several pieces for which we
introduce additional notation which will be most helpful in the derivation of the estimates. Let
$H$ be a continuous strictly positive function defined in $I$. We define 
$$y(w;I,H,k) = \ \mathrel{\mathop {\rm Sup}_{t \in I}} H^{-1}(t) \ |w(t)|_k \ , \eqno(4.4)$$
$$y_1(w;I,H,k) =  \left \{ \int_I dt | \rho ' (t) | H^{-2}(t)| w(t) |_{k + \nu /2}^2
\right \}^{1/2} \ ,  \eqno(4.5)$$
$$Y(w;I, H, k) = (y \vee y_1) (w;I,H,k) \ , \eqno(4.6)$$   
$$z(\varphi ;I,H,\ell) = \ \mathrel{\mathop {\rm Sup}_{t \in I}} H^{-1}(t) \ |\varphi (t)|_{\ell} \ ,
\eqno(4.7)$$
$$z_1(w;I,H,\ell ) = \left \{ \int_I dt | \rho ' (t) | H^{-2}(t)| \varphi (t) |_{\ell + \nu /2}^2
\right \}^{1/2} \ ,  \eqno(4.8)$$
$$Z(\varphi ;I, H, \ell ) = (z \vee z_1) (\varphi ,I,H,\ell ) \ . \eqno(4.9)$$   
\noi We then take
$$\parallel (w, \varphi ); {\cal X}_{\rho}^{k,\ell}(I)\parallel = Y(w;I,1,k) + Z(\varphi ; I, 1,
\ell ) \ . \eqno(4.10)$$

In order to study the Cauchy problem for the auxiliary system (2.11)-(2.12) we need
additional a priori estimates of solutions of that system, which are a continuation of those
of Lemmas 3.5, 3.6 and 3.7 in which we now take into account the dependence on time in the
framework of the spaces ${\cal X}_{\rho}^{k,\ell}$ just defined. \\

\noi {\bf Lemma 4.1.} {\it Let $k$, $\ell$ satisfy (3.48). Let $ 1 \leq T \leq t_0 < \infty$ and
let $(w, \varphi ) \in {\cal X}_{\rho , loc}^{k, \ell} ([T, \infty ))$ be a solution of the system
(2.11)-(2.12). Let $h_0$ and $h_1$ be ${\cal C}^1$ positive functions defined in $[T, \infty
)$, with $h_0$ nondecreasing, $h_1$ nonincreasing, $h_0 \geq t^{-\gamma} |\rho '|^{-1}$ and
$h_1 \geq t^{-2}|\rho '|^{-1}h_0$. Let $a_0 = |w(t_0)|_k$ and $b_0 = h_0 (t_0)^{-1} |\varphi
(t_0)|_{\ell}$. \par
(1) There exist constants $c$ and $C$ such that if
$$\left ( b_0 + a_0^2 \right ) \ h_1 (t_0) \leq c \ , \eqno(4.11)$$
\noi then $(w, h_0^{-1} \varphi ) \in {\cal X}_{\rho}^{k, \ell} ([t_0, \infty ))$ and $(w, \varphi
)$ satisfies the estimates} 
$$\left \{ \begin{array}{l} Y(w;[t_0 , \infty ), 1, k) \leq C \ a_0 \ , \\ \\ Z(\varphi ; [t_0,
\infty ], h_0 , \ell ) \leq C \left ( b_0 + a_0^2 \right ) \ . \end{array} \right . \eqno(4.12)$$

{\it (2) There exist constants $c$ and $C$ such that if
$$\left ( b_0 + a_0^2 \right ) \ T^{-2}|\rho ' (T) |^{-1} h_0 (t_0)  \leq c \ , \eqno(4.13)$$
\noi then $(w, \varphi )$ satisfies the estimates}
$$\left \{ \begin{array}{l} Y(w;[T , t_0 ], 1, k) \leq C \ a_0 \ , \\ \\ Z(\varphi ; [T, t_0
], 1 , \ell ) \leq C \left ( b_0 + a_0^2 \right ) \ h_0(t_0) \ . \end{array} \right . \eqno(4.14)$$

\vskip 3 truemm

\noi {\bf Proof.} The proof requires the same regularization procedure as that of (the integrated
form of) the estimates of Lemma 3.5, but we omit it for brevity. We begin the proof by treating
simultaneously the cases $t \geq t_0$ and $t \leq t_0$. Let $H$ be a ${\cal C}^1$ positive function
of time, increasing for $t \geq t_0$ and decreasing for $t \leq t_0$ and let $\widetilde{\varphi} =
H^{-1} \varphi$. From Lemma 3.5 we obtain 
$$\partial_t |w|_k^2 \mathrel{\mathop <_{>}} 2 \rho ' |w|_{k + \nu /2}^2 \pm C \ t^{-2} \left \{
|w|_{k + \nu /2}^2 \ |\varphi|_{\ell} + |w|_{k + \nu /2} \ |\varphi |_{\ell + \nu /2} \ |w|_k
\right \} \eqno(4.15)$$
$$\partial_t |\widetilde{\varphi}|_{\ell}^2  \mathrel{\mathop <_{>}} 2 \rho '
|\widetilde{\varphi}|_{\ell + \nu /2}^2 \pm C \ t^{-2} |\widetilde{\varphi}|_{\ell + \nu /2}^2 \
|\varphi |_{\ell} \pm C \ t^{-\gamma} \ H^{-1} |\widetilde{\varphi}|_{\ell + \nu /2} \ |w|_{k +
\nu /2} \ |w|_k \eqno(4.16)$$
\noi for $t \displaystyle{\mathrel{\mathop >_{<}}} t_0$. In (4.16), we have dropped the term $-
2H'H^{-1}|\widetilde{\varphi } |_{\ell}^2$ coming from the derivative of $H$. Let now $y_{(1)} =
y_{(1)}(w;I,1,k)$ and $z_{(1)} = z_{(1)}(\varphi ; I, H, \ell )$ where $I = [t_0, t]$ for $t \geq
t_0$ and $I = [t, t_0]$ for $t \leq t_0$. Integrating (4.15) (4.16) over time and using the
Schwarz inequality, we obtain  
$$|w(t)|_k^2 + y_1^2 \leq y_0^2 + C \left \{  \mathrel{\mathop {\rm Sup}_{t \in I}}
t^{-2} |\rho '|^{-1} H \right \} \left ( y_1^2 \ z + y_1 \ z_1 \ y \right )\eqno(4.17)$$
$$|\widetilde{\varphi}(t)|_{\ell}^2 + z_1^2 \leq z_0^2 + C \left \{  \mathrel{\mathop {\rm Sup}_{t
\in I}} \left ( t^{-2} |\rho '|^{-1} H \right ) \right \}  z_1^2 \ z + C \left \{ \mathrel{\mathop
{\rm Sup}_{t \in I}}  t^{-\gamma} |\rho '|^{-1} H^{-1} \right \}z_1 \ y_1 \ y \eqno(4.18)$$
\noi where $y_0 = |w(t_0)|_k = a_0$ and $z_0 = |\widetilde{\varphi}(t_0)|_{\ell}$, and
since the RHS of (4.17) (4.18) is increasing in $|t - t_0|$,
$$y^2 \vee y_1^2 \leq \hbox{RHS of (4.17)}$$
$$z^2 \vee z_1^2 \leq \hbox{RHS of (4.18)}$$
\noi so that $Y = y \vee y_1$ and $Z = z \vee z_1$ satisfy 
$$\left \{ \begin{array}{l} Y^2 \leq y_0^2 + m_1 \ Y^2\ Z \\ \\ Z^2 \leq
z_0^2 + m_1 \ Z^3 + m_2 \ Z\  Y^2\end{array}\right .$$
\noi with
$$m_1 = C \mathrel{\mathop {\rm Sup}_{t \in I}} t^{-2}|\rho '|^{-1} H \quad , \ m_2 =
C \mathrel{\mathop {\rm Sup}_{t \in I}} \left ( t^{-\gamma}|\rho '|^{-1} H^{-1} \right )
\ .$$
\noi If we can arrange that $m_1 Z \leq 1/2$, then we obtain  
$$\left \{ \begin{array}{l} Y^2 \leq 2y_0^2  \\ \\ Z^2 \leq
2z_0^2 + 2m_2 \ Z \ Y^2 \leq 2z_0^2 + 4m_2 \ y_0^2\  Z\end{array}\right .$$
\noi so that 
$$\left \{ \begin{array}{l} Y \leq 2a_0\\ \\ Z \leq
2z_0 + 4m_2 \ a_0^2 \end{array}\right . \eqno(4.19)$$
\noi which will yield the final estimates, while the condition $m_1Z \leq 1/2$ is
implied by
$$4m_1 \left ( z_0 + 2m_2 \ a_0^2 \right ) \leq 1 \ . \eqno(4.20)$$

We now consider separately the cases $t \geq t_0$ and $t \leq t_0$. \par

For $t \geq t_0$, we choose $H = h_0$, so that $m_2 \leq 1$ and $m_1 \leq h_1(t_0)$.
Furthermore $z_0 = b_0$. Then (4.19) reduces to (4.12) and (4.20) reduces to (4.11).
\par

For $t \leq t_0$, we choose $H = 1$, so that $m_2 = t_0^{-\gamma} |\rho ' (t_0)|^{-1}
\leq h_0(t_0)$, $z_0 = b_0 \ h_0 (t_0)$ and $m_1 = t^{-2} |\rho ' (t) |^{-1}$. Then
(4.19) reduces to (4.14) and (4.20) reduces to (4.13) with $t$ replaced by $T$. \par
\nobreak
\hfill $\sq$ \par

\noi {\bf Remark 4.1.} The optimal time decay results in Lemma 4.1 are obtained by
saturating the conditions on $h_0$ and $h_1$, namely by taking $h_0 = t^{-\gamma} |\rho
'|^{-1}$, which we have already assumed to be nondecreasing, and $h_1 = t^{-2} |\rho
'|^{-1} h_0 = t^{-2-\gamma} |\rho '|^{-2}$, in so far as that function is
nonincreasing, a property that we could have (but have not) assumed. We have stated
Lemma 4.1 with inequalities instead of the previous special choices, in order not to
hide the flexibility allowed by the proof, and we shall proceed in the same way in all
subsequent similar estimates. In the special case $|\rho '| = t^{-1 -\varepsilon}$,
the optimal decays obtained with the previous special choices are $h_0 = t^{1 - \gamma
+ \varepsilon}$ and $h_1 = t^{-\gamma + 2 \varepsilon}$. \par

We now turn to the extension of the regularity result of Lemma 3.6. \\

\noi {\bf Lemma 4.2.} {\it Let $k$, $\ell$ satisfy (3.48) and let $\bar{k}$,
$\bar{\ell}$ satisfy $\bar{k} - k = \bar{\ell} - \ell \geq 0$. Let $1 \leq T \leq t_0
< \infty$ and let $(w, \varphi ) \in {\cal X}_{\rho , loc}^{\bar{k}, \bar{\ell}} ([T,
\infty ))$ be a solution of the system (2.11)-(2.12). Let $h_0$ and $h_1$ be as in
Lemma 4.1 and assume that $(w, \varphi )$ satisfies}
$$|w(t)|_k \leq a \quad , \ |\varphi (t)|_{\ell} \leq b \ h_0(t) \quad \hbox{\it or} \
|\varphi (t)|_{\ell} \leq b \ h_0 (t \vee t_0) \ . \eqno(4.21)$$ 

{\it (1) There exist constants $c$ and $C$ such that if
$$(b + a^2) \ h_1 (t_0) \leq c \ , \eqno(4.22)$$
\noi then $(w, h_0^{-1}\varphi ) \in {\cal X}_{\rho}^{\bar{k}, \bar{\ell}} ([t_0 ,
\infty ))$ and $(w, \varphi )$ satisfies the estimates}

$$\left \{ \begin{array}{l} Y(w;[t_0, \infty ), 1, \bar{k}) \leq C \left \{
|w(t_0)|_{\bar{k}} + a \ h_1 (t_0) \ h_0 (t_0)^{-1}\  | \varphi (t_0)|_{\bar{\ell}} \right \}
\ , \\ \\ Z(\varphi ;[t_0, \infty ), h_0, \bar{\ell }) \leq C \left \{
h_0(t_0)^{-1} \ |\varphi (t_0)|_{\bar{\ell }} + a |w (t_0)|_{\bar{k}}\right \} \ .
\end{array}\right . \eqno(4.23)$$

{\it (2) In the case where $|\varphi (t)|_{\ell} \leq b \ h_0 (t \vee t_0)$, there
exist constants $c$ and $C$ such that if
$$(b + a^2) \ T^{-2} |\rho ' (T)|^{-1} \ h_0(t_0) \leq c \ , \eqno(4.24)$$
\noi then $(w, \varphi )$ satisfies the estimates}
 $$\left \{ \begin{array}{l} Y(w;[t, t_0], 1, \bar{k}) \leq C \left \{
|w(t_0)|_{\bar{k}} + a \ t^{-2}\ |\rho '|^{-1} \ | \varphi (t_0)|_{\bar{\ell}}
\right \}  \quad \hbox{\it for all} \ t \in [T, t_0]\ , \\ \\ Z(\varphi ;[T, t_0],
1, \bar{\ell }) \leq C \left \{ |\varphi (t_0)|_{\bar{\ell }} + a \ h_0 (t_0)\  |w
(t_0)|_{\bar{k}}\right \} \ . \end{array}\right . \eqno(4.25)$$

{\it (3) In the case where $|\varphi (t)|_{\ell} \leq b\ h_0 (t)$, there exist
constants $c$ and $C$ such that if
$$(b + a^2) \ h_1 (T) \leq c \ , \eqno(4.26)$$
\noi then $(w, \varphi )$ satisfies the estimates}
 $$\left \{ \begin{array}{l} Y(w;[t, t_0], h_0^{-1}, \bar{k}) \leq C \left \{
h_0(t_0) \ |w(t_0)|_{\bar{k}} + a \ h_1(t)\ | \varphi (t_0)|_{\bar{\ell}}
\right \}  \quad \hbox{\it for all} \ t \in [T, t_0]\ , \\ \\ Z(\varphi ;[T, t_0],
1, \bar{\ell }) \leq C \left \{ |\varphi (t_0)|_{\bar{\ell }} + a \ h_0 (t_0)\  |w
(t_0)|_{\bar{k}}\right \} \ . \end{array}\right . \eqno(4.27)$$

\vskip 3 truemm

\noi {\bf Proof.} The proof follows closely that of Lemma 4.1. Let $h_2$ and $h_3$ be
${\cal C}^1$ positive functions of time, increasing for $t \geq t_0$ and decreasing for
$t \leq t_0$, and let $\widetilde{w} = h_2^{-1} w$, $\widetilde{\varphi} = h_3^{-1}
\varphi$. From Lemma 3.6 we obtain for $t \displaystyle{\mathrel{\mathop >_{<}}} t_0$
$$\partial_t |\widetilde{w}|_{\bar{k}}^2 \mathrel{\mathop <_{>}} 2 \rho '
|\widetilde{w}|_{\bar{k} + \nu /2}^2 \pm C \ t^{-2} \left \{ |\widetilde{w}|_{\bar{k}+
\nu /2}^2 \ |\varphi |_{\ell} + h_3 \ h_2^{-1} \ |\widetilde{w}|_{\bar{k} + \nu /2} \
|\widetilde{\varphi}|_{\bar{\ell} + \nu /2} \ |w|_k \right \} \ , \eqno(4.28)$$
 $$\partial_t |\widetilde{\varphi}|_{\bar{\ell}}^2 \mathrel{\mathop <_{>}} 2 \rho '
|\widetilde{\varphi}|_{\bar{\ell} + \nu /2}^2 \pm C \ t^{-2} \  
|\widetilde{\varphi}|_{\bar{\ell}+ \nu /2}^2 \ |\varphi |_{\ell} \pm C \ t^{-\gamma}
\ h_3^{-1}\ h_2   \ |\widetilde{\varphi}|_{\bar{\ell} + \nu /2} \
|\widetilde{w}|_{\bar{k} + \nu /2} \ |w|_k  \eqno(4.29)$$ 
\noi where we have omitted the terms containing $h'_2$ and $h'_3$. We define $y_{(1)} =
y_{(1)} (w; I, h_2, \bar{k})$ and $z_{(1)} = z_{(1)} (\varphi ; I, h_3, \bar{\ell})$
where $I = [t_0, t]$ for $t \geq t_0$ and $I = [t, t_0]$ for $t \leq t_0$. Proceeding
as in the proof of Lemma 4.1, we obtain from (4.28) (4.29)
$$\left \{ \begin{array}{l} y^2 \vee y_1^2 \leq y_0^2 + m_0 \ y_1^2 + m_1 \ a \ y_1 \
z_1 \\ \\  z^2 \vee z_1^2 \leq z_0^2 + m_0 \ z_1^2 + m_2 \ a \ y_1 \
z_1 \end{array}\right . \eqno(4.30)$$
\noi where $y_0 = |\widetilde{w}(t_0)|_{\bar{k}}$, $z_0 =
|\widetilde{\varphi}(t_0)|_{\bar{\ell}}$,
$$m_0 = C \mathrel{\mathop {\rm Sup}_{t \in I}} t^{-2} |\rho '|^{-1} \ |\varphi
(t)|_{\ell} \ ,$$
$$m_1 = C \mathrel{\mathop {\rm Sup}_{t \in I}} t^{-2} |\rho '|^{-1} \ h_3 \ h_2^{-1}
\quad , \quad m_2 = C \mathrel{\mathop {\rm Sup}_{t \in I}} t^{-\gamma} \ |\rho '|^{-1} \
h_2 \ h_3^{-1} \ . \eqno(4.31)$$
\noi If we can arrange that $m_0 \leq 1/2$, then $Y = y \vee y_1$ and $Z = z \vee z_1$
satisfy 
$$\left \{ \begin{array}{l} Y^2 \leq 2y_0^2 + 2m_1 \ a\ Y \ Z \\ \\  Z^2 \leq
2z_0^2 + 2 m_2  \ a \ Y \ Z \ .\end{array}\right . \eqno(4.32)$$
\noi By an elementary computation, one obtains from (4.32) the estimates 
$$\left \{ \begin{array}{l} Y \leq 4 \left (y_0 + 2a\ m_1 \ z_0 \right ) \\ \\  Z \leq
4\left (z_0 + 2 a \ m_2 \ y_0 \right ) \end{array}\right . \eqno(4.33)$$
\noi under the condition
$$8 a^2 \ m_1 \ m_2 \leq 1 \ . \eqno(4.34)$$

We now consider separately the various cases at hand. \par

For $t \geq t_0$, we take $h_2 = 1$ and $h_3 = h_0$, so that
$$y_0 = |w(t_0)|_{\bar{k}} \quad , \quad z_0 = h_0 (t_0)^{-1} \ |\varphi
(t_0)|_{\bar{\ell}} \ ,$$
$$m_0 \leq C \ b\ h_1(t_0) \quad , \ m_1 \leq C \ h_1(t_0) \quad , \ m_2 \leq C \ .$$
\noi The estimate (4.33) then reduces to (4.23), while the conditions $m_0 \leq 1/2$
and (4.34) recombine to yield (4.22). \par

For $t \leq t_0$, in the case where $|\varphi (t)|_{\ell} \leq b\ h_0 (t_0)$, we take
$h_2 = h_3$ = 1, so that
$$y_0 = |w(t_0)|_{\bar{k}} \quad , \quad z_0 = |\varphi (t_0)|_{\bar{\ell}} \ ,$$
$$m_0 \leq C \ b \ t^{-2} \ \rho '^{-1} \ h_0(t_0)\ , \quad m_1 \leq C \ t^{-2} \ \rho
'^{-1} \quad \hbox{and} \ m_2 \leq C \ h_0 (t_0) \ ,$$
\noi thereby obtaining (4.25) under the condition (4.24). \par

For $t \leq t_0$, in the case where $|\varphi (t)|_{\ell} \leq b \ h_0
(t)$, we take $h_2 = h_0^{-1}$, $h_3 = 1$, so that 
$$y_0 = h_0 (t_0) \  |w(t_0)|_{\bar{k}} \quad , \quad z_0 = |\varphi
(t_0)|_{\bar{\ell}} \ ,$$
$$m_0 \leq C \ b\ h_1(t)\ , \quad m_1 \leq C \ h_1(t) \quad \hbox{and} \
m_2 \leq C \ ,$$
\noi thereby obtaining (4.27) under the condition (4.26). \par \nobreak
\hfill $\sq$ \par

We now turn to the estimate of the difference of two solutions of the system
(2.11)-(2.12).\\

\noi {\bf Lemma 4.3.} {\it Let $k$, $\ell$ satisfy (3.48) and let $k'$, $\ell '$
satisfy (3.57). Let $1 \leq T \leq t_0 < \infty$. Let $h_0$ and $h_1$ be as in Lemma
4.1 and let $(w_i, \varphi_i)$, $i = 1,2$, be two solutions of the system (2.11)-(2.12) such
that $(w_i, h_0^{-1} \varphi_i) \in {\cal X}_{\rho}^{k,\ell} ([T, \infty ))$ and such that
$(w_i, \varphi_i)$ satisfy the estimates $$Y \left ( w_i , [T, \infty ), 1, k \right ) \leq
a \ , \eqno(4.35)$$   $$Z \left ( \varphi_i , [t_0, \infty ), h_0, \ell  \right ) \leq b \ ,
\eqno(4.36)$$ \noi and either
$$Z \left ( \varphi_i , [T, t_0 ), 1, \ell  \right ) \leq b\ h_0(t_0)  \eqno(4.37)$$
\noi or 
$$Z \left ( \varphi_i , [T, t_0 ), h_0, \ell  \right ) \leq b \ . \eqno(4.38)$$
\noi Let $w_{\pm} = w_1 \pm w_2$ and $\varphi_{\pm} = \varphi_1 \pm \varphi_2$. \par

(1) There exist constants $c$ and $C$ such that under the condition (4.22), $(w_-,
\varphi_-)$ satisfies the estimates
$$\left \{ \begin{array}{l} Y(w_-;[t_0, \infty ), 1, k ') \leq C \left \{
|w_-(t_0)|_{k'}+ a \ \ h_1 (t_0) \ h_0(t_0)^{-1} \ |\varphi_-(t_0)|_{\ell '} \right \} \
,  \\ \\  Z(\varphi_-;[t_0, \infty ), h_0, \ell ') \leq C
\left \{ h_0(t_0)^{-1}\ |\varphi _-(t_0)|_{\ell '} + a \
|w_-(t_0)|_{k '} \right \} 
  \ .\end{array}\right . \eqno(4.39)$$

(2) In the case where $\varphi_i$ satisfy (4.37), there exist constants $c$ and $C$
such that under the condition (4.24), $(w_-, \varphi_-)$ satisfies the estimates
$$\left \{ \begin{array}{l} Y(w_-;[t, t_0, ], 1, k ') \leq C \left \{
|w_-(t_0)|_{k'}+ a \ t^{-2} \ |\rho '|^{-1} \ |\varphi_-(t_0)|_{\ell '} \right \} 
  \quad  \hbox{for all} \ t \in [T, t_0]\ , \\ \\ Z(\varphi_-;[T,t_0], 1, \ell ') \leq C
\left \{ |\varphi _-(t_0)|_{\ell '} + a \ h_0(t_0)\  |w_-(t_0)|_{k '} \right \} 
  \ .\end{array}\right . \eqno(4.40)$$

(3) In the case where $\varphi_i$ satisfy (4.38), there exist constants $c$ and $C$
such that under the condition (4.26) $(w_-, \varphi_-)$ satisfies the estimates
$$\left \{ \begin{array}{l} Y(w_-,[t, t_0, ], h_0^{-1}, k ') \leq C \left \{
h_0(t_0)\ |w_-(t_0)|_{k'}+ a \ h_1(t) \  |\varphi_-(t_0)|_{\ell
'} \right \} \quad  \hbox{for all} \ t \in [T, t_0]\ , \\ \\ Z(\varphi_-,[T,t_0], 1,
\ell ') \leq C \left \{ |\varphi _-(t_0)|_{\ell '} + a \ h_0(t_0)\  |w_-(t_0)|_{k '}
\right \} 
  \ .\end{array}\right . \eqno(4.41)$$

(The estimates (4.39) (4.40) and (4.41) are obtained from (4.23) (4.25) and (4.27) by
replacing $(w, \varphi )$ by $(w_- , \varphi_-)$ and $(\bar{k}, \bar{\ell})$ by $(k',
\ell '$).}\\

\noi {\bf Proof.} The proof follows closely that of Lemma 4.2. Let $h_2$ and $h_3$ be
${\cal C}^1$ positive functions of time, increasing for $t \geq t_0$ and decreasing for
$t \leq t_0$ and let $\widetilde{w}_- = h_2^{-1} w_-$, $\widetilde{\varphi}_- =
h_3^{-1} \varphi_-$. Let $H(t) = h_0(t)$ in cases where (4.36) (4.38) are relevant,
$H(t) = h_0(t_0)$ in the case where (4.37) is relevant. From Lemma 3.7 we obtain for $t
\displaystyle{\mathrel{\mathop >_{<}}} t_0$
$$\partial_t |\widetilde{w}_-|_{k'}^2 \mathrel{\mathop <_{>}} 2 \rho ' \ 
|\widetilde{w}_-|_{k' + \nu /2}^2 \pm C \ t^{-2} \Big \{ |\widetilde{w}_-|_{k'+
\nu /2}^2 \ |\varphi_+ |_{\ell} + |\widetilde{w}_-|_{k' + \nu /2} \
|\widetilde{w}_-|_{k'} \ |\varphi_+|_{\ell + \nu /2} \Big \}$$ 
$$ \pm C \ t^{-2} \ h_3 \
h_2^{-1} \  |\widetilde{w}_-|_{k'+ \nu /2} 
 \Big \{ |\widetilde{\varphi}_-|_{\ell '} \ |w_+|_{k + \nu /2} + 
|\widetilde{\varphi}_-|_{\ell '+ \nu /2} \ |w_+ |_{k}  \Big \} \ ,
\eqno(4.42)$$

$$\partial_t |\widetilde{\varphi }_-|_{\ell '}^2 \mathrel{\mathop <_{>}} 2 \rho ' \ 
|\widetilde{\varphi}_-|_{\ell ' + \nu /2}^2 \pm C \ t^{-2} \Big \{
|\widetilde{\varphi}_-|_{\ell '+ \nu /2}^2 \ |\varphi_+ |_{\ell} +
|\widetilde{\varphi}_-|_{\ell ' + \nu /2} \ |\widetilde{\varphi}_-|_{\ell'} \
|\varphi_+|_{\ell + \nu /2} \Big \}$$  $$ + C \ t^{-\gamma} \ h_2 \
h_3^{-1} \  |\widetilde{\varphi}_-|_{\ell '+ \nu /2} 
 \Big \{ |\widetilde{w}_-|_{k '} \ |w_+|_{k + \nu /2} + 
|\widetilde{w}_-|_{k '+ \nu /2} \ |w_+ |_{k}  \Big \} \ ,
\eqno(4.43)$$
\noi where we have omitted the terms containing $h'_2$ and $h'_3$. We define $y_{(1)} =
y_{(1)} (w_-;I,h_2, k')$ and $z_{(1)} = z_{(1)} (\varphi_-; I, h_3, \ell ')$ where $I =
[t_0, t]$ for $t \geq t_0$ and $I = [t, t_0]$ for $t \leq t_0$. Proceding as in the
proof of Lemmas 4.1 and 4.2, we obtain from (4.42) (4.43) supplemented by (4.35)-(4.38)
applied to $(w_+, \varphi_+)$ 
$$\left \{ \begin{array}{l} y^2 \vee y_1^2 \leq y_0^2 + m_0 \ b\ y_1 (y +
y_1) + m_1 \ a \ y_1 (z + z_1) \\ \\ z^2 \vee z_1^2 \leq z_0^2 + m_0 \ b\ z_1 (z +
z_1) + m_2 \ a \ z_1 (y + y_1)\end{array} \right . \eqno(4.44)$$ \noi where 
$$y_0 = |\widetilde{w}_-(t_0)|_{k'} \quad , \quad z_0 =
|\widetilde{\varphi}_-(t_0)|_{\ell'}\ , $$ 
$$m_0 = C \mathrel{\mathop {\rm Sup}_{t \in
I}} \left ( t^{-2} \ |\rho '|^{-1} \ H(t) \right )$$ \noi and $m_1$, $m_2$ are defined by
(4.31). From there on, the proof is identical with that of Lemma 4.2, with (4.30)
replaced by (4.44). \par \nobreak \hfill $\sq$ \par

\noi {\bf Remark 4.2.} It is an unfortunate feature of Lemmas 4.2 and 4.3 that the
derivation of regularity and of difference estimates requires a large time restriction
(see (4.22) (4.24) (4.26)) whereas one would expect those estimates to hold for all
times where the solution is a priori defined, since those estimates are linear in the
higher or difference norm and are expected to follow from some kind of Gronwall's
inequality. The reason for that fact is the occurrence of integral norms in the
definition of the spaces, for which we obtain algebraic inequalities which require
some kind of smallness condition in order to enable us to conclude. \par

In practice, the conditions (4.22) (4.24) (4.26) required for those estimates to hold
have the same form and the same dependence on basic parameters such as $a_0$, $b_0$ as
the conditions that will be needed anyway in order to derive the a priori estimates on
one single solution that are needed to solve the Cauchy problem. We shall impose all
such conditions together, without any significant limitation on the range of validity
of the results (see the proof of Proposition 4.1 below). \par

We now turn to the Cauchy problem for large time for the system (2.11)-(2.12). \\

\noi {\bf Proposition 4.1.} {\it Let $(k , \ell )$ satisfy (3.48) and in addition $k
\geq 1 - \nu /2$. Let $h_0$ and $h_1$ be ${\cal C}^1$ positive functions defined in
$[1, \infty )$ with $h_0$ nondecreasing, $h_1$ nonincreasing and tending to zero at
infinity, $h_0 \geq t^{-\gamma} |\rho '|^{-1}$ and $h_1 \geq t^{-2} |\rho '|^{-1}
h_0$. Let $a_0 > 0$, $b_0 > 0$. \par

(1) There exists $T_0 < \infty $, depending on $a_0$, $b_0$, such that for all $t_0
\geq T_0$, there exists $T \leq t_0$, depending on $a_0$, $b_0$ and $t_0$, such that
for any $(w_0, \varphi_0) \in K_{\rho_{_0}}^k \oplus Y_{\rho_{_0}}^{\ell}$, where $\rho_0 =
\rho (t_0)$, satisfying $|w_0|_k \leq a_0$, $|\varphi_0|_{\ell} \leq h_0(t_0)^{-1}
b_0$, the system (2.11)-(2.12) has a unique solution in the interval $[T, \infty )$
with $w(t_0) = w_0$, $\varphi (t_0) = \varphi_0$, such that $(w, h_0^{-1} \varphi )
\in {\cal X}_{\rho}^{k, \ell} ([T, \infty ))$. One can define $T_0$ and $T$ by
$$\left ( b_0 + a_0^2 \right ) \ h_1 (T_0) = c \eqno(4.45)$$
$$T^{-2} \ |\rho ' (T)|^{-1} \ h_0(t_0) \ h_1 (T_0)^{-1} = 1 \eqno(4.46)$$
\noi and the solution $(w, \varphi )$ satisfies the estimates
$$Y(w;[T, \infty ), 1, k ) \leq C \ a_0 \ , \eqno(4.47)$$
$$Z (\varphi ; [t_0, \infty ), h_0, \ell ) \leq C \left ( b_0 + a_0^2 \right ) \ ,
\eqno(4.48)$$ 
$$Z (\varphi ; [T, t_0 ], 1, \ell ) \leq C \left ( b_0 + a_0^2 \right ) h_0(t_0) \ .
\eqno(4.49)$$ 

(2) If $(w_0, \varphi_0) \in K_{\rho_0}^{\bar{k}} \oplus Y_{\rho_0}^{\bar{\ell}}$ for
some $\bar{k}$, $\bar{\ell}$ with $\bar{k} - k = \bar{\ell} - \ell > 0$, then $(w,
h_0^{-1} \varphi ) \in {\cal X}_{\rho}^{\bar{k}, \bar{\ell}} ([T, \infty ))$, possibly
after changing the constant $c$ in (4.45) (see Remark 4.2). \par

(3) The map $(w_0, \varphi_0) \to (w, \varphi )$ is norm continuous on the bounded sets
of $K_{\rho_{_0}}^k \oplus Y_{\rho_{_0}}^{\ell}$ from the norm of $(w_0, \varphi_0)$ in
$K_{\rho_{_0}}^{k'} \oplus Y_{\rho_{_0}}^{\ell '}$ to the norm of $(w, h_0^{-1} \varphi )$ in
${\cal X}_{\rho}^{k'' , \ell ''}([T, \infty ))$ for $k' \geq \nu /2$, $k - k' = \ell -
\ell ' \geq 1 - \nu$, $k - k'' = \ell - \ell '' > 0$. Furthermore the same map is
continuous from the same topology on $(w_0, \varphi_0)$ to the weak-$*$ topology of
$(w, h_0^{-1} \varphi )$ in ${\cal X}_{\rho}^{k, \ell } ([T, \infty ))$.} \\

\noi {\bf Proof.} \underbar{Part (1)}. The proof proceeds in several steps using a
parabolic regularization and a limiting procedure. We consider first the case $t \geq
t_0$. We shall then indicate briefly the modifications needed in the case $t \leq
t_0$. We shall need the function $\widetilde{w}$ defined by $\widetilde{w}(t) =
U(1/t) w(t)$. Since the operator $U(1/t)$ is unitary in $K_{\rho}^k$ for all $k$ and
$\rho$, all subsequent norm estimates for $\widetilde{w}$ in $K_{\rho}^k$ are
identical with the same estimates for $w$.\par

\noi \underbar{Step 1}. Parabolic regularization and local resolution. \par

We rewrite the system (2.11)-(2.12) in the equivalent form
$$\left \{ \begin{array}{l} \partial_t \widetilde{w} = \left ( 2t^2
\right )^{-1} \ U(1/t) \left ( 2s \cdot \nabla + (\nabla \cdot s) \right )
U(1/t)^* \widetilde{w} \equiv G_1(\widetilde{w} , s) \ , \\ \\ \partial_t \varphi = \left ( 2t^2
\right )^{-1}\ |s|^2 + t^{-\gamma} \ g_0 \left ( U(1/t)^* \widetilde{w} \right )
\equiv G_2(\widetilde{w}, s) \ .  \end{array} \right .\eqno(4.50)$$
   \noi We introduce a parabolic regularization and consider the regularized system 
$$\left \{ \begin{array}{l} \partial_t \widetilde{w} = \theta\ \Delta\ \widetilde{w} +
G_1(\widetilde{w} , s)  \\ \\ \partial_t \varphi = \theta\ \Delta\ \varphi +
G_2(\widetilde{w}, s)   \end{array} \right . \eqno(4.51)$$
\noi with $\theta > 0$. We also regularize the initial data $(\widetilde{w}_0,
\varphi_0)$ at time $t_0$ to $(\bar{\widetilde{w}}_0, \bar{\varphi}_0) \in
X_{\rho_{_0}}^{\bar{k}, \bar{\ell}} \equiv K_{\rho_{_0}}^{\bar{k}} \oplus
Y_{\rho_0}^{\bar{\ell}}$ with $\rho_0 = \rho (t_0)$ , $\bar{k} \geq k \vee (1 + \nu
/2)$ , $\bar{\ell} - \ell = \bar{k} - k$. For the purpose of proving Part (1), one can
take equality in the previous condition, namely $\bar{k} = k \vee (1 + \nu /2)$ and
in particular that second regularization is unnecessary if $k \geq 1 + \nu /2$. We
continue the argument with general $(\bar{k}, \bar{\ell})$ because it will be useful
for the proof of Part (2). \par

We rewrite the Cauchy problem for the system (4.51) in the integral form
$${\widetilde{w} \choose \varphi} (t) = V_{\theta} (t- t_0) {\bar{\widetilde{w}}_0
\choose \bar{\varphi}_0} + \int_{t_0}^t dt' \ V_{\theta} (t - t') {G_1 (\widetilde{w},s)
\choose G_2(\widetilde{w},s)} (t')  \eqno(4.52)$$
\noi where $V_{\theta}(t) \equiv \exp (\theta t \Delta )$ is a contraction in
$X_{\rho_{_0}}^{\bar{k}, \bar{\ell}}$ and satisfies the bound
$$\parallel \nabla V_{\theta} (t) ; {\cal L} (X_{\rho_{_0}}^{\bar{k},\bar{\ell}})\parallel \
\leq C (\theta t)^{-1/2} \ . \eqno(4.53)$$
\noi By (3.23) (3.24) (3.26) (3.27) (3.28) (3.29) of Lemma 3.4, we estimate 
$$\left \{ \begin{array}{l} |G_1(\widetilde{w}, s)|_{\bar{k}-1} \leq C \ t^{-2}
\left (|\widetilde{w}|_{\bar{k}} \ |\varphi |_{\bar{\ell} - 1} +
|\widetilde{w}|_{\bar{k} - 1} \ |\varphi |_{\bar{\ell}} \right ) \\ \\ 
|G_2(\widetilde{w}, s)|_{\bar{\ell}-1} \leq C \ t^{-2} \ |\varphi |_{\bar{\ell}} \
|\varphi |_{\bar{\ell} - 1} + C\ t^{-\gamma} \  |\widetilde{w}|_{\bar{k}} \
|\widetilde{w}|_{\bar{k} - 1} 
\end{array} \right
. \eqno(4.54)$$    
\noi under the conditions (which follow from (3.48))
$$\bar{\ell} > n/2 \ , \ \bar{\ell} \geq \bar{k} - 1 \geq 0 \ , \ \bar{k} \geq
\bar{\ell} + \lambda - 1 \ , \ 2\bar{k} > \bar{\ell} + \lambda + n/2 \ , $$
\noi where $\lambda = \mu - n +2$. In (4.54) the various norms are taken with constant
$\rho (t) \equiv \rho_0$. By a standard contraction argument, the system (4.52) has a
unique solution
$$(\widetilde{w}, \varphi ) \in {\cal C} \left ( [t_0, t_0 + \tau )],
X_{\rho_{_0}}^{\bar{k}, \bar{\ell}} \right ) \eqno(4.55)$$
\noi where one can take
$$\tau = C \ \theta \left ( t_0^{-2} \ h_0(t_0) + t_0^{-\gamma} \ h_0(t_0)^{-1} \right
)^{-2} \left ( \bar{a}_0 + \bar{b}_0 \right )^{-2} \eqno(4.56)$$
\noi with
$$\bar{a}_0 = |\bar{\widetilde{w}}_0|_{\bar{k}} \quad , \quad \bar{b}_0 = h_0(t_0)^{-1}
\ |\bar{\varphi}_0|_{\bar{\ell}} \ . $$
\noi Furthermore, from the estimates
$$\partial_t |\widetilde{w}|_{\bar{k}}^2 + 2 \theta |\nabla \widetilde{w}|_{\bar{k}}^2
= 2 {\rm Re} \ <|\xi |^{\bar{k} + 1} \ f\ \widehat{\widetilde{w}}, |\xi |^{\bar{k}-1} \
f\ \widehat{G_1} (\widetilde{w} , \varphi )>  + \ \hbox{lower order terms}$$
$$ \leq |\widetilde{w}|_{\bar{k} + 1} \ A_1 \left ( t,\ 
|\widetilde{w} |_{\bar{k}} , \ |\varphi |_{\bar{\ell}} \right )  \ , \eqno(4.57)$$
$$\partial_t |\varphi|_{\bar{\ell}}^2 + 2 \theta |\nabla \varphi|_{\bar{\ell}}^2
= 2 <|\xi |^{\bar{\ell} + 3} \ f\ \widehat{\varphi}, |\xi |^{\bar{\ell}+1} \
f\ \widehat{G_2} (\widetilde{w} , \varphi )>  + \ \hbox{lower order terms}$$
$$ \leq |\varphi|_{\bar{\ell} + 1} \ A_2 \left ( t,\ 
|\widetilde{w} |_{\bar{k}} , \ |\varphi |_{\bar{\ell}} \right ) \eqno(4.58)$$
\noi for some estimating functions $A_1$ and $A_2$, it follows that
$$( \widetilde{w} , \varphi ) \in L^2 \left ( [t_0 , t_0 + \tau ], X_{\rho_{_0}}^{\bar{k} +
1, \bar{\ell} + 1} \right ) \ . \eqno(4.59)$$
\noi The estimates (4.57) (4.58) are derived with the help of a regularization
$j_{\varepsilon}$, which we have omitted for brevity, and of the estimate (4.54). \par

\noi \underbar{Step 2}. Uniform estimates and globalisation. \par

>From the regularity conditions (4.55) (4.59), from the fact that $\rho$ is decreasing
and from Lemmas 4.1 and 4.2, especially (4.12) and (4.23), it follows that
$(\widetilde{w}, \varphi ) \in {\cal X}_{\rho}^{\bar{k}, \bar{\ell}} ([t_0 , t_0 + \tau
])$ and that $(\widetilde{w}, \varphi )$ satisfies the estimates 
$$\left \{ \begin{array}{l} Y(\widetilde{w};I,1,k) \leq C\ a_0 \\ \\ Z(\varphi ; I, h_0,
\ell ) \leq C \left ( b_0 + a_0^2 \right ) \ , \end{array} \right . \eqno(4.60)$$
        $$\left \{ \begin{array}{l} Y(\widetilde{w};I,1,\bar{k}) \leq C\left ( 
\bar{a}_0 + h_1(t_0) a_0 \ \bar{b}_0 \right )\\ \\ Z(\varphi ; I, h_0, \bar{\ell}
) \leq C \left ( \bar{b}_0 + a_0\ \bar{a}_0 \right ) \ , \end{array} \right .
\eqno(4.61)$$
\noi for $I = [t_0, t_0 + \tau ]$ and $t_0 \geq T_0$, under the condition (4.45) which ensures
(4.11) (4.22). The estimates (4.60) (4.61) are uniform in $\theta$. From (4.56), from (4.60) (4.61)
for general $I$ and from the fact that $\rho$ is decreasing, it follows by a minor
modification of a standard globalisation argument that $(\widetilde{w}, \varphi )$ can
be continued to a solution of the system (4.51) such that $(\widetilde{w}, h_0^{-1}
\varphi )$ belongs to ${\cal X}_{\rho}^{\bar{k}, \bar{\ell}} ([t_0, \infty ))$ and that
$(\widetilde{w}, \varphi )$ satisfies (4.60) (4.61) with $I = [t_0, \infty )$. \par

\noi \underbar{Step 3}. Limiting procedure. \par

We now take the limits $\theta \to 0$ and $(\bar{\widetilde{w}}_0, \bar{\varphi}_0) \to
(\widetilde{w}_0, \varphi_0)$ in that order. We first keep $(\bar{\widetilde{w}}_0,
\bar{\varphi}_0)$ fixed and consider two solutions $(\widetilde{w}_1, \varphi_1)$ and
$(\widetilde{w}_2, \varphi_2)$ with $(\widetilde{w}_i , h_0^{-1} \varphi_i) \in {\cal
X}_{\rho_{_0}}^{\bar{k},\bar{\ell}}([t_0, \infty ))$ as obtained in Step 2, corresponding
to two values $\theta_1$ and $\theta_2$. We estimate the difference $(\widetilde{w}_-,
\varphi_-) = (\widetilde{w}_1 - \widetilde{w}_2, \varphi_1 - \varphi_2)$ by a minor
variation of Lemma 4.3 with $(k, \ell )$ replaced by $(\bar{k}, \bar{\ell})$ and $(k',
\ell ') = (\bar{k} - 1, \bar{\ell} - 1)$, under the condition (4.22) which follows
from $t_0 \geq T_0$ and from (4.45) (4.60), possibly after changing the constant $c$.
More precisely in the proof of (4.39), we take the initial condition
$\widetilde{w}_-(t_0) = 0$, $\varphi_-(t_0) = 0$, but we have an additional term
coming from the parabolic regularization in the analogue of (4.42) (4.43), namely 
$$\left \{ \begin{array}{l} \partial_t |\widetilde{w}_-|_{\bar{k}-1}^2 \leq 2
|\widetilde{w}_-|_{\bar{k}} \left \{ \theta_1 |\widetilde{w}_1|_{\bar{k}} +
\theta_2 |\widetilde{w}_2|_{\bar{k}} \right ) + \hbox{previous terms} \ ,\\
\\ \partial_t |\varphi_-|_{\bar{\ell}-1}^2 \leq 2
|\varphi_-|_{\bar{\ell}} \left ( \theta_1 |\varphi_1|_{\bar{\ell }} +
\theta_2 |\varphi_2|_{\bar{\ell }} \right ) + \hbox{previous terms} \ .\end{array} \right
. \eqno(4.62)$$ 
\noi (The tildas in (4.62) and in (4.42) (4.43) have different meanings, but this has
no implication on the argument). From (4.61) and (4.62), for any $t_1$ with $t_0 < t_1
< \infty$, we obtain estimates of the type (4.44) for the quantities $y_{(1)} =
y_{(1)}(w_-;[t_0,t_1],1,\bar{k} - 1)$ and  $z_{(1)} = z_{(1)} (\varphi_-; [t_0, t_1], h_0, 
\bar{\ell} - 1)$ where now  $$y_0^2 \vee z_0^2 \leq |t_1 - t_0|(\theta_1 + \theta_2)
\ A(t_0, \bar{a}_0, \bar{b}_0) \ . \eqno(4.63)$$

Therefore, by the same argument as in the proof of (4.39)
$$\left \{ \begin{array}{l} Y(\widetilde{w}_-;[t_0,t_1], 1, \bar{k} - 1) \leq
C \left ( y_0 + a_0 \ h_1(t_0) z_0 \right ) \\ \\ Z(\varphi _-;[t_0,t_1],
h_0, \bar{\ell } - 1) \leq C \left ( z_0 + a_0 \ y_0 \right )\ , \end{array} \right
. \eqno(4.64)$$
\noi which implies that the solution $(\widetilde{w}_{\theta} , \varphi_{\theta})$
associated with $\theta$ converges in norm in ${\cal X}_{\rho}^{\bar{k}-1,
\bar{\ell}-1}([t_0, t_1])$ when $\theta \to 0$ for all $t_1\geq t_0$. Furthermore, the
limit $(\widetilde{w}, \varphi )$ is such that $(\widetilde{w}, h_0^{-1} \varphi ) \in
{\cal X}_{\rho}^{\bar{k}, \bar{\ell}} ([t_0, \infty ))$ and $(\widetilde{w}, \varphi )$ satisfies
(4.47) (4.48). This follows from the bound (4.61) with $I = [t_0, \infty )$, which is uniform in
$\theta$, and from the previous convergence by standard compactness arguments, except
for the strong continuity in time. The latter follows from the weak continuity, which
also follows from a compactness argument, and from the fact that the
$K_{\rho}^{\bar{k}} \oplus Y_{\rho}^{\bar{\ell}}$ norm of $(\widetilde{w}, \varphi)$ is (absolutely)
continuous in $t$ by Lemma 3.5 with $k$, $\ell$ replaced by $\bar{k}$, $\bar{\ell}$. The limit
obviously satisfies the system (4.50). \par

We let now $(\bar{\widetilde{w}}_0, \bar{\varphi}_0)$ tend to $(\widetilde{w}_0,
\varphi_0)$ in $K_{\rho_{_0}}^k \oplus Y_{\rho_{_0}}^{\ell}$ (that step is not needed if
$\bar{k} = k \geq 1 + \nu /2$). Let $(\bar{\widetilde{w}}_{0i} , \bar{\varphi}_{0i})
\in K_{\rho_{_0}}^{\bar{k}} \oplus Y_{\rho_{_0}}^{\bar{\ell}}$, $i = 1,2
$, be two sets of regularized initial conditions and let $(\widetilde{w}_i,
\varphi_i)$ be the solutions of the system (4.50) obtained previously. The difference
$(\widetilde{w}_-, \varphi_-) = (\widetilde{w}_1 - \widetilde{w}_2, \varphi_1 -
\varphi_2)$ is then estimated by (4.39) with $a$ replaced by $a_0$ as follows from
(4.60), under the condition $t_0 \geq T_0$ and (4.45) as previously. One can (but need
not) take $k' = k - 1 + \nu$, $\ell ' = \ell - 1 + \nu$. This implies that the
solution $(\bar{\widetilde{w}}, \bar{\varphi})$ associated with $(\bar{\widetilde{w}}_0,
\bar{\varphi}_0)$ converges in the norm of $(\bar{\widetilde{w}}, h^{-1}
\bar{\varphi})$ in ${\cal X}_{\rho}^{k', \ell '} ([t_0 , \infty ))$ when
$(\bar{\widetilde{w}}_0, \bar{\varphi}_0)$ converges to $(\widetilde{w}_0,
\varphi_0)$ in the norm of $K_{\rho_{_0}}^{k'} \oplus Y_{\rho_{_0}}^{\ell '}$ on the
bounded sets of $K_{\rho_{_0}}^k \oplus Y_{\rho_{_0}}^{\ell}$ (and a fortiori in the norm
of $K_{\rho_{_0}}^k \oplus Y_{\rho_{_0}}^{\ell}$). Let $(\widetilde{w}, \varphi )$ be the
limit. By the same arguments as above, $(\widetilde{w}, h_0^{-1} \varphi ) \in
{\cal X}_{\rho}^{k,\ell}([t_0, \infty ))$ and $(\widetilde{w}, \varphi )$
satisfies the system (4.50) and the estimates (4.47) (4.48). \par
     
\noi \underbar{Step 4}. Uniqueness follows immediately from Lemma 4.3, part (1). \par

We now turn to the case $t \leq t_0$. The proof proceeds exactly in the same way, with
Parts 1 of Lemmas 4.1, 4.2 and 4.3 replaced by Parts 2 of the same Lemmas. In the same
way as before, (4.14) implies that (4.24) follows from (4.13), possibly after a change
of constant $c$. With $T_0$ defined by (4.45) and possibly with another change of
constant $c$, the condition (4.13) follows from (4.46). \par

\noi \underbar{Part (2)}. If $\bar{k} \geq 1 + \nu /2$, the result follows from the
proof of Part (1) with the second limiting procedure omitted. If $\bar{k} < 1 + \nu
/2$, the result follows from the proof of Part (1) with $\bar{k}$ replaced by $1 + \nu
/2$ and with the second limiting procedure going down from $1 + \nu /2$ to $\bar{k}$
instead of going down from $1 + \nu /2$ to $k$. \par

Note that in the previous proof, the constant $c$ in (4.45) comes from a successive
application of Lemma 4.1, especially from (4.11) (4.13), and of Lemmas 4.2 and 4.3,
especially from (4.22) (4.24). In particular that constant depends on $(k , \ell )$
and on $(\bar{k}, \bar{\ell})$ (the pair $(k', \ell ')$ in the applications of Lemma
4.3 is chosen as a function of $(k, \ell )$ or $(\bar{k}, \bar{\ell}))$. In the proof
of Part (1), one can in addition choose $\bar{k} = k \vee (1 + \nu /2)$, $\bar{k} - k
= \bar{\ell} - \ell$, and the constant $c$ therefore depends only on $(k, \ell )$. In
contrast with that, in the proof of Part (2), the pairs $(k , \ell )$ and $(\bar{k},
\bar{\ell})$ are independent, and the constant $c$ in (4.45) needed for Part (2) to
hold may (and is expected to) depend on both $(k , \ell )$ and $(\bar{k},
\bar{\ell})$. The crucial information contained in Part (2) is the fact that (4.45)
does not involve $|w_0|_{\bar{k}}$, $|\varphi_0|_{\bar{\ell}}$, but only $|w_0|_k$ and
$|\varphi_0|_{\ell}$. \par

\noi \underbar{Part (3)}. Continuity with respect to initial data follows from Lemma
4.3 parts 1 and 2, and from the a priori estimates (4.47) (4.48) (4.49) by
interpolation and compactness arguments. \par \nobreak
\hfill $\sq$ \par

We conclude this section with two properties of the behaviour at infinity of solutions
of the system (2.11)-(2.12) in spaces ${\cal X}_{\rho}^{k, \ell}$. There is no initial
time involved in those properties and $\rho$ is only required to satisfy a suitable
monotony condition at infinity. The first property is the existence of a limit for
$w(t)$ as $t \to \infty$. It applies in particular to the solutions constructed in
Proposition 4.1. There is a large flexibility on the assumptions under which such a
limit exist. We shall give another example in the next section in a different context.\\

\noi {\bf Proposition 4.2.} {\it Let $T \geq 1$, $\rho_{\infty} \geq 0$ and
$$\rho (t) = \rho_{\infty} + \int_t^{\infty} dt_1 |\rho ' (t_1)|$$
\noi for $t \geq T$. Let $k$, $\ell$ satisfy
$$\ell > n/2 - 1 \quad , \quad 1 - \nu /2 \leq k \leq \ell + 1 \ . \eqno(4.65)$$
\noi Let $h_0$ and $h_1$ be as in Proposition 4.1. Let $(w, \varphi )$ be a solution of
the system (2.11)-(2.12) such that $(w, h_0^{-1} \varphi ) \in {\cal
X}_{\rho}^{k,\ell}([T , \infty ))$. Let 
$$a = Y (w; [T, \infty ), 1, k) \quad , \quad b = Z(\varphi ; [T , \infty ), h_0, \ell
) \ .$$
\noi Then there exists $w_+ \in K_{\rho_{\infty}}^k$ such that $w(t)$ tends to $w_+$
strongly in $K_{\rho_{\infty}}^{k'}$ for $k' < k$ and weakly in $K_{\rho_{\infty}}^k$
when $t \to \infty$. Furthermore, the following estimates hold~:
$$|w_+|_k \equiv \parallel w_+;K_{\rho_{\infty}}^k\parallel \ \leq a \ , \eqno(4.66)$$
$$Y(\widetilde{w} - \widetilde{w}(t_1);[t_1, \infty ), 1, k') \leq C \ a\ b\ h_1(t_1)
\ , \eqno(4.67)$$
$$\parallel \widetilde{w}(t) - w_+ ; K_{\rho_{\infty}}^{k'} \parallel \leq C \ a\ b\
h_1(t) \ , \eqno(4.68)$$
\noi for $k' = k - 1 + \nu$, for all $t$, $t_1 \geq T$, and with $ \widetilde{w}(t) =
U(1/t) w(t)$.} \\

\noi {\bf Proof.} We first prove (4.67). Let $t \geq t_1$ and $w_1(t) = 
\widetilde{w}(t) -  \widetilde{w}(t_1)$ so that by (4.50)
$$\partial_t  w_1 =
(2t^2)^{-1} U(1/t)(2s \cdot \nabla + (\nabla \cdot s)) U(1/t)^*  \widetilde{w}\ .$$ 

Using (3.12) and (3.23) (3.24) from Lemma 3.4 with $k' = k - 1 + \nu$ , $m = k' - \nu /2$, we obtain
$$\partial_t |w_1|_{k'}^2 - 2 \rho ' |w_1|_{k' + \nu /2}^2 \leq C \ t^{-2} |w_1|_{k' +
\nu /2} \left \{ |w|_{k + \nu /2} \ |\varphi|_{\ell} + |\varphi|_{\ell + \nu /2} \ |w|_k \right \} 
\eqno(4.69)$$  
\noi under the conditions $k' \geq \nu /2$ , $\ell > n/2 - 1$ , $\ell \geq k - 1$, which
reduce to (4.65). Let $y_{(1)} = y_{(1)} (w_1; [t_1, t],1,k')$. Integrating (4.69)
over time in the same way as in Lemma 4.1, we obtain 
$$y^2 \vee y_1^2 \leq C \ ab\ y_1 \left \{ \mathrel{\mathop {\rm Sup}_{t' \in
[t_1,t]}}t'^{-2} |\rho '(t')|^{-1} \ h_0(t') \right \} \leq C\ ab\ y_1\ h_1(t_1)$$
\noi from which (4.67) follows. From (4.67) and from the fact that $\rho$ is
decreasing, it follows that $w(t)$ has a strong limit $w_+ \in K_{\rho_{_\infty}}^{k'}$ and
that (4.68) holds. By standard compactness and interpolation arguments, $w_+ \in
K_{\rho_{\infty}}^k$, $w_+$ satisfies (4.66) and $w(t)$ converges to $w_+$ weakly in
$K_{\rho_{_\infty}}^k$ and strongly in $K_{\rho_{_\infty}}^{k'}$ for $k' < k$. \par
\nobreak
\hfill $\sq$ \par

The second property is a uniqueness property for solutions with suitable restrictions
on the behaviour at infinity. It will be used in Section 6. Since it corresponds to a
situation with infinite initial time, it requires $\rho$ to be increasing. \\

\noi {\bf Proposition 4.3.} {\it Let $T \geq 1$, $\rho_{\infty} > 0$ and
$$\rho (t) = \rho_{\infty} - \int_t^{\infty} dt_1 |\rho ' (t_1)|$$
\noi with $\rho (t) \geq 0$ for $t \geq T$. Let $(k, \ell )$ satisfy (3.48) and $k \geq
1 - \nu /2$. Let $h_0$ and $h_1$ be as in Proposition 4.1. Let $(w_i, \varphi_i)$, $i =
1,2$ be two solutions of the system (2.11)-(2.12) such that $(w_i, h_0^{-1} \varphi_i)
\in {\cal X}_{\rho}^{k, \ell} ([T, \infty ))$ and such that
$$|w_1(t) - w_2(t)|_{k'}  \ h_0(t) \to 0 \quad , \quad |\varphi_1(t) -
\varphi_2(t)|_{\ell '} \to 0$$
\noi when $t \to \infty$, for some $(k', \ell ')$ such that $\nu /2 \leq k' \leq k - 1
+ \nu$, $\ell - \ell ' = k - k'$. Let 
$$a = \ \mathrel{\mathop {\rm Max}_{i}} Y\left ( w_i;[T, \infty ), 1,k\right ) \quad ,
\quad   b = \ \mathrel{\mathop {\rm Max}_{i}} Z\left (\varphi_i;[T, \infty ), h_0,\ell
\right )\ .$$

Then there exists a constant $c$ such that if (4.26) holds, then $(w_1, \varphi_1) =
(w_2 , \varphi_2)$.} \\

\noi {\bf Proof.} The result follows immediately from Lemma 4.3, part (3) by taking the
limit $t_0 \to \infty$ in (4.41).\par
\nobreak
\hfill $\sq$ \par

\section{Asymptotics of solutions of the auxiliary system}
\hspace*{\parindent} In this section, we continue the study of the asymptotic
properties of the solutions of the auxiliary system (2.11)-(2.12) obtained in
Proposition 4.1. We have already proved the existence of a limit $w_+$ for $w(t)$ as $t
\to \infty$ for such solutions. Under suitable additional regularity assumptions in the
form of stronger lower bounds on $(k, \ell)$, we shall obtain estimates on the
asymptotic behaviour in time of the asymptotic functions $(w_m, \varphi_m)$ and $(W_m
, \phi_m)$ defined by (2.14)-(2.15) and estimates on the remainders $(q_{p+1},
\psi_{p+1}) = (w - W_p, \varphi - \phi_p)$, also defined by (2.14)-(2.15), eventually
leading to the existence of an asymptotic state for the phase $\varphi$ in the form of
a limit $\psi_+$ for $\psi_{p+1}$ as $t \to \infty$ for sufficiently large $t$. \par

At a technical level however, the situation here differs significantly from that in
Section 4. We are no longer trying to solve the system (2.11)-(2.12), but instead we
assume a solution of that system to be given, and we estimate successively the
asymptotic functions $(w_m, \varphi_m)$ and the remainders $(q_m, \psi_m)$, which are
defined by triangular systems of equations. As a consequence, there is no need to
control the loss of derivatives as in (2.11)-(2.12), and we can simply let that loss
accumulate in the solution of the triangular system. Therefore we no longer need a
time dependent $\rho$ (see (3.12)), integral norms in ${\cal X}_{\rho}^{k,\ell}$,
and integration by parts (see (3.22) (3.25)). We do not even need Gevrey spaces, and
we could use instead ordinary Sobolev spaces, in the same way as in II. We shall of
course nevertheless keep Gevrey spaces in order to make contact with Section 4, but
in all this section we shall take $\rho$ to be constant. Instead of the spaces
${\cal X}_{\rho}^{k,\ell}$ defined by (4.2), we shall use the simpler spaces
$${\cal Y}_{\rho}^{k, \ell}(I) = ({\cal C} \cap L^{\infty}) (I, K_{\rho}^k \oplus
Y_{\rho}^{\ell}) \ . \eqno(5.1)$$

The contact with Section 4, in particular the applicability of the results of this
section to the solutions obtained in Proposition 4.1, will be achieved through the
fact that if $\rho$ is defined in $[t_0, \infty )$ by (4.1) or equivalently by
$$\rho (t) = \rho_{\infty} + \int_t^{\infty} dt_1 \ |\rho '(t_1)| \ , \eqno(5.2)$$
     \noi then
$${\cal X}_{\rho}^{k,\ell} ([t_0 , \infty )) \hookrightarrow {\cal
Y}_{\rho_{\infty}}^{k,\ell} ([t_0, \infty )) \ . \eqno(5.3)$$

In all the estimates of this section, the function $f$ occurring in the definition
(3.8) (3.9) of the spaces plays no role whatsoever, and is consistently eliminated from
the proofs by using the submultiplicativity property (3.3). As a consequence all the
estimates are uniform in (actually independent of) $\rho$ and $\nu$, and no assumption
is made connecting $\nu$ to other parameters such as $\mu$ or $(k, \ell)$~: we only
assume $0 \leq \nu \leq 1$ and $\rho \geq 0$. \par

As a preliminary result, we give an existence result of the limit $w_+$ of $w(t)$ as
$t \to \infty$ for a solution of the system (2.11)-(2.12). That result is a
simplified version of Proposition 4.2 appropriate to the present context. \\

\noi {\bf Proposition 5.1.} {\it Let $(k, \ell )$ satisfy}
$$k \geq 1 \quad , \quad \ell \geq 0 \quad , \quad k + \ell > n/2 \ . \eqno(5.4)$$
\noi {\it Let $1 \leq T < \infty$ and let $h_0$ be a ${\cal C}^1$ positive
nondecreasing function defined in $[T, \infty )$ with $t^{-2} h_0 \in L^1([T, \infty
))$. Let $(w, \varphi )$ be a solution of the system (2.11)-(2.12) such that $(w,
h_0^{-1} \varphi ) \in {\cal Y}_{\rho}^{k, \ell} ([T, \infty ))$. Then there exists
$w_+ \in K_{\rho}^k$ such that $w(t)$ tends to $w_+$ strongly in $K_{\rho}^{k'}$ for
$k' < k$ and weakly in $K_{\rho}^k$. Furthermore, the following estimate holds~: 
$$|\widetilde{w}(t) - w_+|_{k'} \leq C \ a\ b\ h(t) \eqno(5.5)$$
\noi for $0 \leq k' \leq k - 1$, $k' \leq \ell$, $k' < k + \ell - n/2$, where
$$a = \ \mathrel{\mathop {\rm Sup}_{t}}|w(t)|_k \quad , \quad b = \ \mathrel{\mathop {\rm
Sup}_{t}}h_0^{-1} (t) |\varphi (t) |_{\ell} \ , \eqno(5.6)$$
\noi $ h(t) = \int_t^{\infty} dt_1\ t_1^{-2} h_0(t_1)$ and $\widetilde{w}(t) = U(1/t)
w(t)$.} \\

\noi {\bf Proof.} By (3.23) (3.24), we estimate
$$|\partial_t (\widetilde{w}(t) - \widetilde{w}(t_1))|_{k'} \leq C \ t^{-2} \
|w(t)|_k \ |\varphi (t) |_{\ell} \leq C \ t^{-2} \ h_0(t) \ a \ b $$
\noi under the conditions stated on $k'$, and therefore by integration
$$|\widetilde{w}(t)  - \widetilde{w}(t_1)|_{k'} \leq C \ a\ b \ h(t \wedge t_1)$$     
\noi which implies the existence of $w_+ \in K_{\rho}^{k'}$ satisfying the estimate
(5.5). By standard compactness arguments, $w_+ \in K_{\rho}^k$, $|w_+|_k \leq a$ and
$w(t)$ tends weakly to $w_+$ in $K_{\rho}^k$ when $t \to \infty$. The other
convergences follow by interpolation. \par \nobreak
\hfill $\sq$ \par

We now estimate the asymptotic functions $(w_m, \varphi_m)$ defined by successive
integrations from the system (2.18)-(2.19) with initial condition (2.20)-(2.21). For
that purpose, we need the function
$$\bar{h}_0 (t) = \int_1^t dt_1 \ t_1^{-\gamma} \eqno(5.7)$$
\noi and the associated estimating functions defined by (3.68) (3.69) (3.72), which we
denote by $\bar{N}_m$, $\bar{Q}_m$ and $\bar{P}_m$ for that special choice of $h_0$.
(Those functions appeared already in II where they were called $h_0$, $N_m$, $Q_m$ and
$P_m$, see (II.3.19) (II.3.25) (II.3.26) (II.3.31)). We recall that $\lambda = \mu - n
+ 2$ and we define $\bar{\lambda} = \lambda \vee 1$. The following proposition is the
extension to the present context of Proposition II.5.1. \\ 
 
\noi {\bf Proposition 5.2.} {\it Let $p \geq 0$ be an integer, let $k$ satisfy
$$k > n/2 \quad , \quad k \geq (p+2) \bar{\lambda} - 1 \ , $$
\noi and let
$$k_m = k - m \bar{\lambda} \quad , \quad \ell_m = k - m\bar{\lambda} - \lambda \quad
, \quad 0 \leq m \leq p + 1 \ . \eqno(5.8)$$ 
\noi Let $w_+ \in K_{\rho}^k$ and let $a = |w_+|_k$. Let $\{w_0 = w_+, w_{m+1}\}$ and
$\{\varphi_m\}$, $0 \leq m \leq p$ be the solution of the system (2.18)-(2.19) with
initial conditions (2.20)-(2.21). Then\par
(1) $w_{m+1} \in {\cal C}([1, \infty )$, $K_{\rho}^{k_{m+1}})$, $\varphi_m \in {\cal
C}([1, \infty ), Y_{\rho}^{\ell_m})$ and the following estimates hold for all $t \geq
1$~: 
$$|w_{m+1}(t)|_{k_{m+1}} \leq A(a) \ \bar{Q}_m(t) \ , \eqno(5.9)$$
$$|\varphi_m(t) |_{\ell_m} \leq A(a) \ \bar{N}_m(t) \ , \eqno(5.10)$$
\noi for some estimating function $A(a)$. \par

If in addition $(p+2)\gamma > 1$ and if we define $\varphi_{p+1}$ by (2.19) with
initial condition $\varphi_{p+1}(\infty ) = 0$, then $\varphi_{p+1} \in {\cal C}([1,
\infty ), Y^{\ell_{p+1}})$ and the following estimate holds~:
$$|\varphi_{p+1}(t)|_{\ell_{p+1}} \leq A(a) \ \bar{P}_p(t) \ . \eqno(5.11)$$

(2) The functions $\{\varphi_m\}$ are gauge invariant, namely if $w'_+ = w_+ \exp (i
\omega )$ for some real valued function $\omega$ and if $w'_+$ gives rise to $\{\varphi
'_m\}$, then $\varphi '_m = \varphi_m$ for $0 \leq m \leq p + 1$.\par

(3) The map $w_+ \to \{w_{m+1}, \varphi_m\}$ is uniformly Lipschitz continuous on the
bounded sets from the norm topology of $w_+$ in $K_{\rho}^k$ to the norms $\parallel
\bar{Q}_m^{-1} w_{m+1}; L^{\infty}([1, \infty ), K_{\rho}^{k_{m+1}})\parallel$ and
$\parallel \bar{N}_m^{-1} \varphi_{m}; L^{\infty}([1, \infty ),
Y_{\rho}^{\ell_{m}})\parallel$, $0 \leq m \leq p$. A similar continuity holds for
$\varphi_{p+1}$.} \\

\noi {\bf Proof.} The proof is essentially the same as that of Proposition II.5.1.
\par

\noi \underbar{Part (1)}. We proceed by induction on $m$, starting from the assumption on $w_+$. We
assume the results to hold for $(w_j, \varphi_j)$ for $j \leq m$ and we prove them for $w_{m+1}$ and
$\varphi_{m+1}$. We first consider $w_{m+1}$. Letting the exponents $(k_m, \ell_m)$ be
undefined for the moment except for being nonincreasing in $m$, we obtain from (2.18)
and from (3.23) (3.24)
$$|\partial_t \ w_{m+1}|_{k_{m+1}} \leq A(a) \ t^{-2} \left \{ \sum_{0\leq j \leq m -
1} \ \bar{N}_j (t) \ \bar{Q}_{m-j-1}(t) + \bar{N}_m(t) \right \} \eqno(5.12)$$
\noi under the conditions
$$\left \{ \begin{array}{l} k_{m+1} \geq 0 \quad , \quad k_{m+1} \leq \left (
k_m - 1 \right ) \wedge \ell_m \ , \\ \\ k_{m+1} + n/2 < k_j + \ell_{m-j}
\quad , \quad 0 \leq j \leq m \ .\end{array} \right . \eqno(5.13)$$
\noi Integrating (5.12) between $t$ and $\infty$ with $w_{m+1}(\infty ) = 0$ and using
(3.79) (3.75) yields the result for $w_{m+1}$.  \par

We next consider $\varphi_{m+1}$. From (2.19) and from (3.26) (3.27) (3.28) (3.29) we
obtain
$$|\partial_t \ \varphi_{m+1}|_{\ell_{m+1}} \leq A(a) \left \{  t^{-2} 
\sum_{0\leq j \leq m} \ \bar{N}_j (t) \ \bar{N}_{m-j}(t) + t^{-\gamma} \left ( 
\sum_{0 \leq j \leq m - 1} \bar{Q}_j(t) \ \bar{Q}_{m-1-j}(t) + \bar{Q}_m(t) \right )\right \}
\eqno(5.14)$$  
\noi under the conditions  
$$\left \{ \begin{array}{l} \ell_{m+1} + 1 \geq 0 \quad , \quad \ell_{m+1} \leq \left (
\ell_m - 1 \right ) \wedge \left ( k_{m+1} - \lambda \right )\ , \\ \\ \ell_{m+1} +
\displaystyle{{n\over 2}}< \ell_j + \ell_{m-j} \quad , \quad \ell_{m+1} + \lambda + n/2
< k_j + k_{m+1-j} \quad , \quad 0 \leq j \leq m \ ,\end{array} \right . \eqno(5.15)$$
\noi and for $\varphi_0$
$$\ell_0 + \lambda \leq k \quad , \quad \ell_0 + \lambda + n/2 < 2 k \ . \eqno(5.16)$$

Integrating (5.14) between 1 and $t$ with $\varphi_m(1) = 0$ and using (3.78) (3.76)
and (3.80) (3.84) yields the result for $\varphi_{m+1}$ (and similarly for
$\varphi_0$). \par

We saturate the nonstrict part of (5.13) (5.15) (5.16) by the choice (5.8), where in
addition we optimize (maximize) $\{\ell_m\}$ for given $k$. The strict conditions
then reduce to $k > n/2$, while the condition $k \geq (p+2)\bar{\lambda} - 1$ is simply the
condition $k_{p+1} \wedge (\ell_{p+1} + 1) \geq 0$. \par

Finally, if $(p+2)\gamma > 1$, we integrate (5.14) with $m = p$ between $t$ and
$\infty$ and use (3.77) (3.78) and (3.80) (3.82) with $m = p$. \par

\noi \underbar{Part (2)}. The proof is identical with that of Proposition II.5.1 and
will be omitted. \par

\noi \underbar{Part (3)}. From the fact that the RHS of (2.18)-(2.19) are bilinear,
it follows by induction of $m$ that the difference between two solutions $\{w_m,
\varphi_m\}$ and $\{w'_m, \varphi '_m\}$ associated with $w_+$ and $w'_+$ is
estimated by 
$$|w_{m+1} - w'_{m+1}|_{k_{m+1}} \leq A(a) \ |w_+ - w'_+|_k \ \bar{Q}_m(t) \ ,
\eqno(5.17)$$
$$|\varphi_{m} - \varphi '_{m}|_{\ell_{m}} \leq A(a) \ |w_+ - w'_+|_k \ \bar{N}_m(t) 
\eqno(5.18)$$  
\noi for $0 \leq m \leq p$, and if $(p+2)\gamma > 1$,
$$|\varphi_{p+1} - \varphi '_{p+1}|_{\ell_{p+1}} \leq A(a) \ |w_+ - w'_+|_k \
\bar{P}_p(t) \ , \eqno(5.19)$$
\noi where $a = |w_+|_k \vee |w'_+|_k$. The continuity stated in Part (3) follows from
those estimates. \par \nobreak
\hfill $\sq$ \par

We now turn to the main result of this section, namely to the proof of existence of
asymptotic states $(w_+, \psi_+)$ for solutions of the auxiliary system (2.11)-(2.12).
That result relies heavily on suitable estimates of the remainders
$$q_{m+1}(t) = w(t) - W_m (t) \eqno(5.20)$$ 
$$\psi_{m+1}(t) = \varphi (t) - \phi_m (t) \eqno(5.21)$$
\noi where $W_m$ and $\phi_m$ are defined (see (2.14) (2.15)) by
$$W_m = \sum_{0 \leq j \leq m} w_j \quad , \quad \phi_m = \sum_{0 \leq j \leq m}
\varphi_j \ . \eqno(5.22)$$
\noi In view of Proposition 4.1, we shall consider solutions $(w, \varphi )$ of the
system (2.11)-(2.12) such that $(w, h_0^{-1} \varphi ) \in {\cal Y}_{\rho}^{k, \ell}
([T, \infty ))$, where $h_0$ is a suitable ${\cal C}^1$ positive increasing function of
time. We shall assume in addition that $t^{-\gamma} \leq c \ h'_0$, a property which
occurs naturally in the interesting examples relevant for Section 4. In addition to the
estimating functions $\bar{N}_m$, $\bar{Q}_m$ and $\bar{P}_m$ associated with
$\bar{h}_0$, we shall also need the estimating functions $N_m$, $Q_m$ and $P_m$
associated with $h_0$, defined by (3.68) (3.69) (3.72). From the relation $t^{-\gamma}
\leq c \ h'_0$, it follows that $\bar{N}_m$, $\bar{Q}_m$ and $\bar{P}_m$ are estimated
by $N_m$, $Q_m$, and $P_m$, and more precisely
$$\bar{N}_m \leq c^{m+1} \ N_m \quad , \quad \bar{Q}_m \leq c^{m+1} \ Q_m \quad ,
\quad \bar{P}_m \leq c^{m+2} \ P_m \ . \eqno(5.23)$$

We can now state the main result of this section, which is the extension of
Proposition II.6.1 to the present context. \\  

\noi{\bf Proposition 5.3.} {\it Let $p \geq 0$ be an integer. Let $(k, \ell , k_0)$
satisfy 
$$k > n/2 +  (p-2)\bar{\lambda} + \lambda \vee 0 \quad , \quad k \geq k_0 -
\bar{\lambda} + 2 \quad , \quad \ell \geq k_0 - \lambda  , \eqno(5.24)$$
$$k_0 > n/2 \quad , \quad k_0 \geq (p+2) \bar{\lambda} - 1 \ , \eqno(5.25)$$
\noi and let
$$k_m = k_0 - m \bar{\lambda} \quad , \quad \ell_m = k_0 - \lambda - m \bar{\lambda}
\quad , \quad 0 \leq m \leq p + 1 \ . \eqno(5.26)$$
\noi Let $h_0$ be a ${\cal C}^1$ positive increasing function defined in $[1, \infty
)$, such that $t^{-2} h_0$, $t^{-1} h'_0 \in L^1([1, \infty ))$ and $t^{-\gamma} \leq
c \ h'_0$. Let $T \geq 2$, let $(w, \varphi )$ be a solution of the system (2.11)-(2.12)
such that $(w, h_0^{-1} \varphi ) \in {\cal Y}^{k,\ell}_{\rho}([T, \infty ))$ and
define $a$, $b$ by (5.6). Let $w_+ = \lim\limits_{t \to \infty} w(t) \in K_{\rho}^k$ be
defined by Proposition 5.1, so that in particular 
$$|w(t) - w_+|_{k_1} \to 0 \quad \hbox{when} \ t \to \infty \ .\eqno(5.27)$$
\noi Let $(w_{m+1}, \varphi_m)$, $0 \leq m \leq p$ be defined by Proposition 5.2 and
let $(W_m, \phi_m)$, $0 \leq m \leq p$, be defined by (5.22). Then the following
estimates hold for all $t \in [T, \infty )$~:
$$|w(t) - W_m(t)|_{k_{m+1}} \leq A(a, b) \ Q_m(t) \ , \eqno(5.28)$$  
$$|\varphi (t) - \phi_m(t)|_{\ell_{m+1}} \leq A(a, b) \ N_{m+1}(t) \ . \eqno(5.29)$$
\noi for $0 \leq m \leq p$, and for some estimating function $A(a, b)$. \par

If in addition $(p+2) \gamma > 1$ and if $P_p (1) < \infty$, then the following limit
exists
$$\lim_{t \to \infty} \left ( \varphi (t) - \phi_p (t) \right ) = \psi_+ \eqno(5.30)$$
\noi as a strong limit in $Y_{\rho}^{\ell_{p+1}}$, and the following estimate holds} 
$$|\varphi (t) - \phi_p(t) - \psi_+|_{\ell_{p+1}} \leq A(a, b) P_p(t) \ . \eqno(5.31)$$  

\noi {\bf Remark 5.1.} When applied to solutions $(w, \varphi )$ of the system (2.11)-(2.12)
obtained in Proposition 4.1, the time decay estimates of Proposition 5.3 will take the
following typical form. Take $|\rho '| = t^{-1 - \varepsilon}$ and $h_0 = t^{1 - \gamma +
\varepsilon}$, which is adequate for Propostion 4.1. Then for $(p+1)\gamma < 1$ and
sufficiently small $\varepsilon$ $$N_m(t) \sim t^{1-(m+1)(\gamma - \varepsilon)} \ , \ Q_m
\sim t^{-(m+1)(\gamma - \varepsilon)} \ , \ P_p (t) \sim t^{1-(p+2)(\gamma - \varepsilon)}$$
\noi and the condition $P_p(1) < \infty$ reduces to $(p+2)(\gamma - \varepsilon ) >
1$. For $(p+1)(\gamma - \varepsilon ) >1$, the time decay saturates at $N_p(t) \sim
1$, $Q_p(t) \sim t^{-1}$ and $P_p (t) \sim t^{-\gamma}$, as explained in Section II.3. \\

\noi {\bf Remark 5.2.} We have kept the parameter $k_0$ in the statement of
Proposition 5.3 because it plays a central role in the proof. For a given solution
$(w, \varphi )$, namely for given $(k, \ell )$, we can optimize the results by
maximizing $k_0$ as allowed by (5.24), namely by taking $k_0 = (k + \bar{\lambda} -
2) \vee (\ell + \lambda )$. The condition (5.25) then reduces to lower bounds on
$(k, \ell )$. Those bounds are stronger than (5.4) and therefore allow for the
application of Proposition 5.1. Note also that the regularity obtained for the
remainders $(q_m, \psi_m)$ for $m \geq 1$ is weaker than that of the estimating
functions of the same level $(w_m,\varphi_m)$ since $(k_m, \ell_m)$ in (5.26)
contain only $k_0 \leq k$ whereas $(k_m, \ell_m)$ in (5.8) contain $k$ and $w_+
\in K_{\rho}^k$ in both cases. \\

\noi {\bf Proof of Proposition 5.3.} The proof is essentially the same as that of
Proposition II.6.1 and proceeds by an induction on $m$, the starting point of which is the
estimate for $q_1$. We assume the estimates (5.28) (5.29) to hold for $(q_j, \psi_j)$, $0 \leq j
\leq m$, and we derive them for $(q_{m+1}, \psi_{m+1})$, with $(q_m, \psi_m)$ defined by (5.20)
(5.21) and $(q_0, \psi_0 ) = (w, \varphi )$. \par

We substitute the decompositions $w = W_m + q_{m+1}$ and $\varphi = \phi_m +
\psi_{m+1}$ in the LHS of (2.11)-(2.12) and the decompositions $w = W_{m-1} +
q_m$ and $\varphi = \phi_{m-1} + \psi_m$ in the RHS of the same, thereby obtaining
$$\partial_t \ q_{m+1} = \left ( 2t^2 \right )^{-1} \Big \{ i \Delta w + (2
\nabla \varphi \cdot \nabla + (\Delta \varphi )) q_m$$
$$+ \left ( 2 \nabla \psi_m \cdot \nabla + (\Delta \psi_m) \right ) \ W_{m-1} +
\sum_{0 \leq i,j \leq m - 1\atop{i+j\geq m}} \left ( 2 \nabla \varphi_i \cdot
\nabla + (\Delta \varphi_i) \right ) w_j \Big \} \eqno(5.32)$$      
$$\partial_t \ \psi_{m+1} = \left ( 2t^2 \right )^{-1} \Big \{ \left ( \nabla \varphi
+ \nabla \phi_{m-1} \right ) \cdot \nabla \psi_m + \sum_{0 \leq i,j \leq m - 1\atop{i+j\geq
m}} \nabla \varphi_i \cdot \nabla \varphi_j \Big \}$$  
$$+ t^{-\gamma} \Big \{ g_0(q_m,q_1) + g_0 (q_m, W_{m-1} - w_0) + 2g_0 (q_{m+1}, w_0)
+ \sum_{0 \leq i,j \leq m - 1\atop{i+j\geq m+1}} g_0(w_i, w_j) \Big \} \eqno(5.33)$$ 
\noi for $m \geq 1$ and 
$$\left \{ \begin{array}{l} \partial_t \ q_{1} = \left ( 2t^2 \right )^{-1}
\left \{ i \ \Delta \ w + 2 \nabla \varphi \cdot \nabla w + (\Delta \varphi )
w \right \} \\ \\ \partial_t \ \psi_{1} = \left ( 2t^2 \right )^{-1} \ |\nabla \varphi
|^2 + t^{-\gamma} \ g_0\left (q_1, w + w_+ \right ) \end{array} \right . \eqno(5.34)$$
\noi for $m = 0$. We let the exponents $(k_m, \ell_m)$ be undefined for the moment,
except for the property of being decreasing in $m$ and of being not larger than the
corresponding exponents of Proposition 5.2, which we denote momentarily by $(\bar{k}_m,
\bar{\ell}_m)$ inside this proof, namely 
$$k_m \leq \bar{k}_m = k - m \bar{\lambda} \quad , \quad \ell_m \leq \bar{\ell}_m = k -
\lambda - m\bar{\lambda} \ .$$
\noi We estimate (5.32) by (3.23) (3.24), by (5.9) (5.10) (5.23) and by the induction
assumption, and we estimate similarly (5.33) by (3.26) (3.27) (3.28) (3.29) and the
same other ingredients, thereby obtaining 
$$|\partial_t \ q_{m+1}|_{k_{m+1}} \leq A(a,b) \ t^{-2} \Big \{ 1 + h_0 \ Q_{m-1} +
N_m + \sum_{m \leq i+j \leq 2(m - 1)} N_i \ Q_{j-1} \Big \} \eqno(5.35)$$
$$|\partial_t \ \psi_{m+1}|_{\ell_{m+1}} \leq A(a,b) \Big \{ t^{-2} \Big ( h_0 \
N_m + \sum_{m \leq i+j \leq 2(m - 1)} N_i \ N_{j} \Big ) $$
$$+ t^{-\gamma} \Big ( Q_0 \ Q_{m-1} + Q_m + \sum_{m +1 \leq i+j \leq 2(m - 1)}
Q_{i-1} \ Q_{j-1} \Big ) \Big \} \eqno(5.36)$$
\noi for $m \geq 1$, under the conditions
$$k_{m+1} \leq \left ( k_m - 1 \right ) \wedge \ell_m \quad , \quad \ell_{m+1} \leq
\left ( \ell_m - 1 \right ) \wedge \left ( k_{m+1} - \lambda \right ) \ , \eqno(5.37)$$
$$\left \{ \begin{array}{l} k_{m+1} + n/2 < \left ( \ell + k_m \right ) \wedge \left (
\ell_m + \bar{k}_{m-1} \right ) \\ \\ \ell_{m+1} + n/2 < \left ( \ell + \ell_m \right )
\wedge \left ( \ell_m + \bar{\ell}_{m-1} \right ) \\ \\ \ell_{m+1} + n/2 + \lambda <
\left ( k_m + k_1 \right ) \wedge \left ( k_m + \bar{k}_{m-1} \right ) \wedge \left (
k_{m+1} + k \right ) \ ,  \end{array} \right . \eqno(5.38)$$
 \noi and
$$|\partial_t\ q_1 |_{k_1} \leq C \ t^{-2} \ (a + h_0 \ ab) \eqno(5.39)$$
$$|\partial_t \ \psi_1 |_{\ell_1} \leq C \ t^{-2} \ h_0^2 \ b^2 + C \ t^{-\gamma} \ A(a,b)
 \ Q_0 \eqno(5.40)$$ 
\noi under the conditions
$$k_1 \leq (k - 2) \wedge \ell \quad , \quad \ell_1 \leq (\ell - 1) \wedge (k_1 -
\lambda ) \ , \eqno(5.41)$$
$$k_1 + n/2 < k + \ell \quad , \quad \ell_1 + n/2 < 2 \ell \wedge (k_1 + k - \lambda )
\eqno(5.42)$$
\noi for $m = 0$. \par

As in the proof of Proposition 5.2, we saturate (5.37) and maximize $\ell_m$ with
respect to $k_m$ by the choice (5.26). Then (5.41) reduces to the last two conditions
of (5.24), while (5.38) is easily seen to reduce to $k_0 > n/2$ and to the first
condition of (5.24), and (5.42) follows from (5.24) and (5.25). The conditions $k_m \leq \bar{k}_m$
and $\ell_m \leq \bar{\ell}_m$ follow from the fact that $k \geq k_0$, implied by (5.24). \par

Integrating (5.39) between $t$ and $\infty$ with initial condition
(5.27) yields
$$|q_1(t)|_{k_1} \leq C \left ( a + h_0(1) \right ) t^{-1} + C \ ab \ Q_0(t) \ , \eqno(5.43)$$
\noi namely the estimate (5.28) for $m = 0$ (since $Q_0(t) \geq t^{-1} Q_0(1))$ which is the
starting point of the induction procedure. Integrating (5.35) between $t$ and
$\infty$  with the initial condition $q_{m+1}(\infty ) = 0$ for $m \geq 1$ and using (3.79)
(3.75) yields (5.28). Similarly, integrating (5.40) and (5.36) between $T$ and $t$ and using
(3.78) (3.76) and (3.80) (3.84) yields (5.29). Finally, if $(p + 2) \gamma > 1$ and
$P_p (1) < \infty$, the RHS of (5.36) with $m = p$ is integrable at infinity in time,
which proves the existence of the limit (5.30). Integrating (5.36) between $t$ and $\infty$ and
using (3.78) (3.77) and (3.80) (3.82) yields (5.31). \par \nobreak
\hfill $\sq$ \par

\section{Cauchy problem at infinity and wave operators for the auxiliary system}
\hspace*{\parindent} In this section we derive the main technical result of this paper,
which is in some sense the converse of Proposition 5.3, namely we prove that
sufficiently regular asymptotic states $(w_+, \psi_+)$ generate solutions $(w, \varphi
)$ of the system (2.11)-(2.12) in the sense described in Section 2, thereby allowing
for the definition of the local wave operator at infinity $\Omega_0$~: $(w_+, \psi_+)
\to (w, \varphi )$. As a preliminary, and in order to allow for an easy proof of the
gauge invariance of the construction, we first solve the linear transport equations
(2.23) (2.24) with initial condition (2.25) and derive some asymptotic properties of
their solutions. \par

In all this section, as in Proposition 4.1 and in contrast with Section 5, we are
again solving the system (2.11)-(2.12) and we have to solve the simpler equations
(2.23) (2.24) in the same framework. As a consequence, we again need the full
machinery of Gevrey spaces with time dependent $\rho$ so as to be able to use (3.12),
we need the integral norms in ${\cal X}_{\rho}^{k,\ell}$ and the integration by parts
(3.22) and (3.25). We therefore begin by choosing $|\rho '|$ exactly as in Section 4. Now
however in contrast with Proposition 4.1 where we kept $t_0$ finite, we want to take
$t_0 = \infty$, and therefore we must take $\rho$ to be increasing. Therefore, in all
this section, we take
$$\rho (t) = \rho_{\infty} - \int_t^{\infty} dt_1 |\rho '(t_1)| \ , \eqno(6.1)$$ 
\noi taking $\rho_{\infty}$ sufficiently large for $\rho (t)$ to be nonnegative in the
(asymptotic) region of interest, a sufficient condition for which being that
$\rho_{\infty} \geq \parallel |\rho '|; L^1([1, \infty ))\parallel$. Except for that
condition, $\rho_{\infty}$ will be arbitrary (but fixed), and all subsequent estimates
will be independent of $\rho_{\infty}$, for the same reasons as in the previous
sections.\par

Solving the Cauchy problem either for the system (2.11)-(2.12) or for the transport
equations (2.23) (2.24) with infinite initial time will be done as in II by first
solving that system (or equation) with large but finite initial time $t_0$ and then
letting $t_0$ tend to infinity. In order not to make this paper too cumbersome, we
shall restrict our attention to solving that problem only for $t \leq t_0$ when
$t_0$ is finite. The solutions could easily be extended to $[t_0, \infty )$ with a
modified $\rho$ of the type (4.1), but that extension would be useless in the limit
$t_0 \to \infty$, and we shall refrain from performing it. \par

By analogy with the spaces ${\cal X}_{\rho}^{k,\ell}(I)$ defined by (4.2) (4.3), we
extend the definition (5.1) of the spaces ${\cal Y}_{\rho}^{k,\ell}(I)$ from the
case of constant $\rho$ to the case of variable $\rho$ considered in this section by 
$${\cal Y}_{\rho}^{k,\ell}(I) = \left \{ (w, \varphi ) : \left ( F^{-1} \ f \ \widehat{w} \
, F^{-1} \ f \ \widehat{\varphi} \right ) \in {\cal Y}_0^{k,\ell} (I) \right \} \eqno(6.2)$$
\noi with $f$ defined by (3.1), $\rho$ defined by (6.1) and ${\cal Y}_0^{k,\ell}(I)$
defined by (5.1). Occasionally, we shall have to state that a single function $w$ or
$\varphi$ belongs to the $w$-subspace or to the $\varphi$-subspace of some space
${\cal X}_{\rho}^{k,\ell}$ or ${\cal Y}_{\rho}^{k,\ell}$. In order to avoid
introducing additional notation, we shall then write $(w, 0) \in {\cal
X}_{\rho}^{k,0}(I)$ or $(0, \varphi ) \in {\cal X}_{\rho}^{0,\ell}(I)$, and similarly
with ${\cal X}$ replaced by ${\cal Y}$. \par

We shall have to consider norms of the type $|w_-|_k$ or $|\varphi_-|_{\ell}$ for
$w_-$ or $\varphi_-$ that are differences of functions taken at different times,
possibly leaving in doubt the value of $t$ appearing in $\rho (t)$ in the
definition of the norm. In such cases it will be understood that the value of $t$
appearing in $\rho (t)$ should be the smaller of the times appearing in $w_-$ or
$\varphi_-$, thereby yielding the smaller of the corresponding values of $\rho$.
\par

We begin with the study of the transport equation (2.23). As compared with the
treatment of that problem given in II, however, a new difficulty arises. In order
to compare $V$ with the asymptotic approximation $W_p$ to the anticipated
solution of (2.11)-(2.12), it is no longer sufficient to take for $V$ the initial
condition $V(t_0) = w_+$ when solving (2.23) with finite initial time $t_0$, and
we have to use instead the better initial condition $V(t_0) = W_p (t_0)$. On the
other hand, the results on the Cauchy problem for (2.23) (2.24) do not depend
on detailed properties of $\phi_{p-1}$ and $W_p$. They are therefore stated in
Propositions 6.1 and 6.2 in terms of general functions $\phi$ and $W$, to be
taken as $\phi_{p-1}$ and $W_p$ from Proposition 6.3 on. \par

In all this section, we use systematically the notation $y_{(1)}$, $Y$,
$z_{(1)}$, $Z$ defined by (4.4)-(4.9). \par

We begin with the study of the transport equation (2.23) which we rewrite with
general $\phi$ as
$$\partial_t \ V = \left ( 2t^2\right )^{-1} \left ( 2 \nabla \phi \cdot \nabla
+ (\Delta \phi ) \right ) V \ . \eqno(6.3)$$        
\vskip 3 truemm  

\noi {\bf Proposition 6.1.} {\it Let $(\widetilde{k}, \widetilde{\ell}, \bar{k}, k)$
satisfy
$$\widetilde{\ell} > n/2 - \nu \quad , \quad \widetilde{k} \wedge (\widetilde{\ell} +
\nu /2) \geq \bar{k} \geq k + 1 - \nu \quad , \quad k \geq \nu /2 \ . \eqno(6.4)$$
\noi Let $1 \leq T < \infty$. Let $h_0$ and $h_1$ be ${\cal C}^1$ positive functions
defined in $[T, \infty )$ with $h_0$ nondecreasing, $h_1$ nonincreasing and tending to
zero at infinity, and $h_1 \geq t^{-2} \rho '^{-1}h_0$. Let $w_+ \in
K_{\rho_{_\infty}}^{\widetilde{k}}$. Let $(W , \phi )$ be such that $(W, h_0^{-1} \phi
) \in {\cal Y}_{\rho}^{\widetilde{k}, \widetilde{\ell}}([T, \infty ))$ and that $W(t)$
tends to $w_+$ as $t \to \infty$, with an estimate 
$$|W(t) - w_+|_k \leq c_1 \ h_1 (t) \eqno(6.5)$$
\noi for some constant $c_1$. Let
$$a = \ \mathrel{\mathop {\rm Sup}_{t}} |W(t)|_{\widetilde{k}} \quad , \quad b = \ \mathrel{\mathop
{\rm Sup}_{t}} h_0^{-1}(t) |\phi (t)|_{\widetilde{\ell}} \ .
\eqno(6.6)$$
\noi Then \par

(1) There exist constants $c$ and $C$ such that if
$$b \ h_1(T) \leq c \ , \eqno(6.7)$$
\noi there exists a unique solution $V$ of the equation (6.3) such that $(V, 0) \in
{\cal X}_{\rho}^{\bar{k},0}([T, \infty ))$ and such that the following estimates hold~:
$$Y(V;[T, \infty ), 1, \bar{k}) \leq C \ a \ , \eqno(6.8)$$
$$Y(V - w_+;[T, \infty ), h_1, k) \leq C a\ b \ . \eqno(6.9)$$ 

(2) $V$ is the limit as $t_0 \to \infty$ of solutions $V_{t_0}$ of (6.3) such that
$V_{t_0}(t_0) = W(t_0)$ and $(V_{t_0}, 0) \in {\cal X}_{\rho}^{\bar{k}, 0}([T, t_0])$.
The convergence is in the strong sense in ${\cal X}_{\rho}^{k', 0}([T, T_1])$ for $k' <
\bar{k}$ and in the weak-$*$ sense in ${\cal X}_{\rho}^{\bar{k}, 0}([T, T_1])$ for
every $T_1$, $T < T_1 < \infty$ and the following estimate holds for all $t_0 > T$
$$Y\left ( V- V_{t_0};[T, t_0], 1, k \right ) \leq C \left ( ab + c_1 \right )\ h_1(t_0)
\ . \eqno(6.10)$$

(3) The solution $V$ is unique in $L^{\infty}([T, \infty ), L^2)$ under the
condition that $\parallel V(t) - w_+\parallel_2$ tends to zero when $t \to \infty$.} \\

\noi {\bf Remark 6.1.} From the uniqueness statement of Proposition 6.1, Part (2), it
follows that for given $\phi$ and $w_+$, $V$ is independent of $W$. Actually Parts (1)
and (3) make no reference to $W$ and could be proved by taking $W(t) \equiv w_+$. $W$
appears only in the limiting process of Part (2), which however will be esssential to
derive the more accurate estimates of Proposition 6.3. \\

\noi {\bf Proof.} We prove Parts (1) and (2) together.

We first solve (6.3) with initial data $V_{t_0}(t_0) = W(t_0)$ at finite $t_0$. This is
a linear transport equation with ${\cal C}^{\infty}$ vector field $\nabla \phi$ and
${\cal C}^{\infty}$ initial data and the existence and uniqueness of a solution, for
instance with value in $H^N$, is a standard result. We concentrate on the Gevrey
estimates and on the subsequent limit $t_0 \to \infty$. \par

In the same way as in the proof of Lemma 3.5, from (3.22) (3.24) we obtain for $t \leq
t_0$
$$\partial_t |V_{t_0}|_{\bar{k}}^2 \geq 2 \rho ' |V_{t_0}|_{\bar{k} + \nu /2}^2 - C \
t^{-2} |V_{t_0}|_{\bar{k}+ \nu /2}^2 \ |\phi |_{\widetilde{\ell}} \eqno(6.11)$$
\noi under the conditions
$$\widetilde{\ell} > n/2 - \nu \quad , \quad \bar{k} \geq \nu /2 \quad , \quad
\widetilde{\ell} + 1 \geq \bar{k} - \nu /2$$
\noi and
$$\bar{k} + \widetilde{\ell} + \nu /2 > \bar{k} - \nu /2 + n/2 \quad , \quad
\widetilde{\ell} \geq \bar{k} - \nu /2$$
\noi which follow from (6.4). Integrating (6.11) over time and using (6.6), we obtain
as in the proof of Lemma 4.1
$$y^2 \vee y_1^2 \leq y_0^2 + C\ b \ y_1^2 \ h_1(t) \eqno(6.12)$$
\noi where $y_{(1)} = y_{(1)}(V_{t_0};[t,t_0], 1, \bar{k})$ and $y_0 =
|W(t_0)|_{\bar{k}}$, which under the condition (6.7) yields
$$Y\left ( V_{t_0};[T,t_0], 1, \bar{k} \right ) \leq C \ a \ . \eqno(6.13)$$

We next estimate the difference $v_1(t) \equiv V_{t_0}(t) - V_{t_0}(t_0) \equiv V_{t_0} (t) -
W(t_0)$, which satisfies the equation 
$$\partial_t \ v_1 = \left ( 2t^2\right )^{-1} \left ( 2 \nabla \phi \cdot \nabla +
(\Delta \phi ) \right ) V_{t_0} \eqno(6.14)$$
\noi with initial condition $v_1(t_0) = 0$. Let $\widetilde{v}_1 = h_1^{-1}v_1$. From
(3.23) (3.24) with $m = k - \nu /2$, $k \to \bar{k} + \nu /2$ and $\ell =
\widetilde{\ell}$, we obtain
$$\partial_t |\widetilde{v}_1|_k^2 \geq 2 \rho ' \ |\widetilde{v}_1|_{k+ \nu /2}^2 - C
\ t^{-2} \ h_1^{-1} \ |\widetilde{v}_1|_{k+ \nu /2} \ |V_{t_0}|_{\bar{k}+ \nu /2} \
|\phi |_{\widetilde{\ell}} \eqno(6.15)$$
\noi under the conditions
$$\bar{k} + \widetilde{\ell} > k + n/2 - \nu \quad , \quad \widetilde{\ell} \geq k -
\nu /2 \geq 0 \quad , \quad \bar{k} + \nu /2 \geq k + 1 - \nu /2 \ , $$
\noi which follow from (6.4). \par

Defining $y_{(1)} = y_{(1)} (v_1;[T, t_0], h_1, k)$, integrating over time and  using
(6.6) (6.13) and the Schwarz inequality we obtain in the same way as before
$$y^2 \vee y_1^2 \leq C \ a\ b \ y_1 \left \{ \mathrel{\mathop {\rm Sup}_{t}}t^{-2} \
\rho '^{-1} \ h_0 \ h_1^{-1} \right \} = C\ a\ b\ y_1$$
\noi and therefore
$$Y \left ( V_{t_0} - W(t_0); [T, t_0], h_1, k \right ) \leq C\ a\ b\ . \eqno(6.16)$$

We now take the limit $t_0 \to \infty$, and for that purpose we estimate the
difference $V_{t_1} - V_{t_0}$ of two solutions corresponding to $t_0$ and $t_1$,
with $T <t_0 \leq t_1$. Since the equation (6.3) is linear in $V$, the difference of
two solutions is estimated in the same way as a single solution. In the same way as
in the proof of (6.13), we obtain 
$$Y\left ( V_{t_1} - V_{t_0};[T, t_0], 1, k\right ) \leq C \left | V_{t_1} (t_0) -
V_{t_0}(t_0) \right |_k$$
$$\leq C \left \{ \left | V_{t_1}(t_0) - W(t_1) \right |_k + \left | W(t_1) -
W(t_0)\right |_k \right \} \ .\eqno(6.17)$$
\noi We estimate the first norm in the last member of (6.17) by (6.16) with $t_1$
replacing $t_0$, and where we use the pointwise estimate taken at $t= t_0$, and we
estimate the second norm by (6.5) used both for $t = t_0$ and $t = t_1$, thereby
obtaining
$$Y\left ( V_{t_1} - V_{t_0};[T, t_0], 1, k \right ) \leq C \left ( ab + c_1 \right )
h_1(t_0) \ .\eqno(6.18)$$
\noi From (6.18) it follows that when $t_0 \to \infty$, $V_{t_0}$ converges to some $V$
such that $(V,0) \in {\cal X}_{\rho}^{k,0}([T, \infty ))$ strongly in ${\cal
X}_{\rho}^{k,0}([T, T_1])$ for all $T_1$, $T < T_1 < \infty$. By standard compactness
arguments and by (6.13) $(V, 0) \in {\cal X}_{\rho}^{\bar{k},0}([T, \infty ))$, $V$
satisfies (6.8) and the convergence holds in the sense of Part (2) of the proposition.
Furthermore, taking the limit $t_1 \to \infty$ in (6.18) yields (6.10). \par

It remains only to prove (6.9). For that purpose, we take again $T < t_0 < t_1$ and we
estimate
$$Y\left ( V - w_+;[T, t_0],h_1, k \right ) \leq Y \left ( V_{t_1} - W(t_1);[T, t_0],
h_1, k \right ) $$
$$+ h_1 (t_0)^{-1} \left \{ Y \left ( V - V_{t_1};[T, t_0], 1, k \right ) + Y \left (
W(t_1) - w_+; [T, t_0], 1, k \right ) \right \}$$
$$\leq C \ a\ b + h_1(t_0)^{-1} \left \{ C\left ( ab + c_1 \right ) h_1(t_1) + C \ c_1 \
h_1 (t_1) \left ( \rho (t_1) - \rho (t_0) \right )^{-1/2} \right \}$$
 \noi by (6.16) and (6.10)
with $t_1$ replacing $t_0$, and by (6.5) and the fact that $\rho (t)$ is strictly
increasing, so that $$\int_T^{t_0} dt \ \rho '|W(t_1) - w_+|_{\rho (t), k + \nu /2}^2
\leq C |W(t_1) - w_+|_{\rho (t_1), k}^2 \ |\rho (t_1) - \rho (t_0)|^{-1} \ \parallel
\rho ' \parallel_1$$ \noi where $|\cdot |_{\rho , k}$ denotes the norm in $K_{\rho}^k$.
Taking the limits $t_1 \to \infty$ and $t_0 \to \infty$ in that order yields (6.9), with
the same constant as in (6.16). \par Part (3) follows from an elementary energy
estimate for the $L^2$ norm of the difference of two solutions, namely $$\parallel
V_1(t) - V_2(t) \parallel_2 \ \leq \ \parallel V_1(t') - V_2(t') \parallel_2 \ \exp \left
( C \ b |h(t) - h(t')| \right )$$  \noi where $h(t) = \int_t^{\infty} dt_1 \ t_1^{-2} \
h_0(t_1)$. \par \nobreak \hfill $\sq$ \par

We now turn to the transport equation (2.24), which we rewrite with general $\phi$ as
$$\partial_t \ \chi = t^{-2} \ \nabla \phi \cdot \nabla \chi \ . \eqno(6.19)$$
\vskip 3 truemm

\noi {\bf Proposition 6.2.} {\it Let $(\widetilde{\ell}, \bar{\ell} , \ell )$ satisfy}
$$\widetilde{\ell} > n/2 - \nu \quad , \quad \widetilde{\ell} + \nu /2 \geq \bar{\ell}
+ 1 \quad , \quad \bar{\ell} \geq \ell + 1 - \nu \quad , \quad \ell + 1 \geq \nu /2 \
. \eqno(6.20)$$
\noi {\it Let $1 \leq T < \infty$ and let $h_0$, $h_1$ be as in Proposition 6.1. Let
$\psi_+ \in Y_{\rho_{_\infty}}^{\bar{\ell}}$ and let $\beta = |\psi_+|_{\bar{\ell}}$.
Let $\phi$ be such that $(0, h_0^{-1}\phi ) \in Y_{\rho}^{0,\widetilde{\ell}}([T, \infty ))$
with 
$$b = \ \mathrel{\mathop {\rm Sup}_{t}}h_0(t)^{-1} \ |\phi (t) |_{\widetilde{\ell}} \ . 
\eqno(6.21)$$
\noi Then \par

(1) There exist constants $c$ and $C$ such that if (6.7) holds, there exists a unique
solution $\chi$ of the equation (6.19) such that $(0, \chi ) \in {\cal
X}_{\rho}^{0,\bar{\ell}}([T, \infty ))$ and such that the following estimates hold
$$Z\left ( \chi ; [T, \infty ), 1, \bar{\ell} \right ) \leq C\  \beta \ ,
\eqno(6.22)$$   
$$Z\left ( \chi - \psi_+;[T, \infty ), h_1, \ell \right ) \leq C \ b\ \beta \ .
\eqno(6.23)$$

(2) $\chi$ is the limit as $t_0 \to \infty$ of solutions $\chi_{t_0}$ of (6.19) such
that $\chi_{t_0}(t_0) = \psi_+$ and $(0, \chi_{t_0}) \in {\cal
X}_{\rho}^{0,\bar{\ell}}([T, t_0))$. The convergence is in the strong sense in ${\cal
X}_{\rho}^{0,\ell '}([T, T_1])$ for $\ell ' \leq \bar{\ell}$ and in the weak-$*$ sense
in ${\cal X}_{\rho}^{0,\bar{\ell}}([T, T_1])$ for every $T_1$, $T < T_1 < \infty$, and
the following estimate holds
$$Z \left ( \chi - \chi_{t_0};[T,t_0], 1, \ell \right ) \leq C \ b\ \beta\ h_1(t_0) \
. \eqno(6.24)$$

(3) The solution $\chi$ is unique in $L^{\infty}(I, L^{\infty})$ under the condition
that $\parallel \chi (t) - \psi_+ \parallel_{\infty}$ tends to zero as $t \to \infty$.
\par

(4) Let in addition $(\widetilde{k}, \bar{k}, k)$ satisfy (6.4), let $w_+ \in
K_{\rho_{_\infty}}^{\widetilde{k}}$ and let $V$ be defined by Proposition 6.1 for some
$W$ (for instance $W(t) \equiv w_+$). Then for fixed $\phi$, $V \exp (- i \chi )$ is
gauge invariant in the following sense. If $(V, \chi )$ and $(V', \chi ')$ are the
solutions obtained from $(w_+, \psi_+)$ and $(w'_+, \psi '_+)$ with $w_+ \exp (- i
\psi_+) = w'_+ \exp (-i \psi '_+)$, then $V(t) \exp (-i \chi (t)) = V'(t) \exp (-i
\chi '(t))$ for all $t \in I$.} \\

\noi {\bf Proof.} \underbar{Parts (1) and (2)}. The proof is very similar to that of
Proposition 6.1 and we concentrate again on the Gevrey estimates and on the limit $t_0
\to \infty$. Let $\chi_{t_0}$ be the solution of (6.19) with initial data
$\chi_{t_0}(t_0) = \psi_+$. In addition to (6.19), it is convenient to use also the
equation
$$\partial_t \ \tau = t^{-2} \left ( S\cdot \nabla \tau + \tau \cdot \nabla S \right )
\eqno(6.25)$$
\noi satisfies by $\tau = \nabla \chi$, with $S = \nabla \phi$. \par

We estimate $\partial_t|\chi_{t_0}|_{\bar{\ell}}^2$ in the same way as in the proof of
Lemma 3.5. We estimate the contribution of $S\cdot \nabla \tau_{t_0}$ from (6.25) by
(3.25) with $(s, s', \ell, \ell ')$ replaced by $(S, \tau_{t_0}, \widetilde{\ell},
\bar{\ell})$ and the contribution of $\tau_{t_0}\cdot \nabla S$ by (3.26) with $(s,
s', \ell, \ell ', m)$ replaced by $(\tau_{t_0}, S, \bar{\ell} + \nu /2,
\widetilde{\ell}, \bar{\ell} - \nu /2)$, thereby obtaining
$$\partial_t |\chi_{t_0}|_{\bar{\ell}}^2 \geq 2 \rho '|\chi_{t_0}|_{\bar{\ell} + \nu
/2}^2 - C \ t^{-2} \ |\chi_{t_0}|_{\bar{\ell} + \nu /2}^2 \
|\phi|_{\widetilde{\ell}} \eqno(6.26)$$ 
\noi under the conditions
$$\widetilde{\ell} > n/2 - \nu \quad , \quad \bar{\ell} + 1 \geq \nu /2 \quad , \quad
\widetilde{\ell} \geq \bar{\ell} + 1 - \nu /2$$
\noi which follow from (6.20). Introducing $z_{(1)} = z_{(1)} (\chi_{t_0};[T, t_0], 1,
\bar{\ell})$ and integrating (6.26) over time, we obtain in the same way as before
$$z^2 \vee z_1^2 \leq z_0^2 + C \ b \ h_1(T) z_1^2$$
\noi with $z_0 = \beta$, which under the conditions (6.7) implies
$$Z\left ( \chi_{t_0};[T, t_0], 1, \bar{\ell} \right ) \leq C\ \beta \ .
\eqno(6.27)$$
\noi We next estimate the difference $\chi_1(t) \equiv \chi_{t_0}(t) -
\chi_{t_0}(t_0) \equiv \chi_{t_0}(t) - \psi_+$ which satisfies the equations
$$\partial_t \ \chi_1 = t^{-2} \ \nabla\phi \cdot \nabla \chi_{t_0}$$
$$\partial_t \ \tau_1 = t^{-2} \left ( S\cdot \nabla \tau_{t_0} + \tau_{t_0}\cdot
\nabla S \right )\ .$$
\noi Let $\widetilde{\chi}_1 = h_1^{-1}\chi_1$. Using again (3.26) with $(s, s') = (S,
\tau_{t_0})$ or $(\tau_{t_0},S)$, with $m = \ell - \nu /2$ and with $(\ell , \ell ') =
(\widetilde{\ell} , \bar{\ell} + \nu /2)$ or $(\bar{\ell} + \nu /2,
\widetilde{\ell})$, we obtain 
$$\partial_t \ |\widetilde{\chi}_1|_{\ell}^2 \geq 2 \rho ' |\widetilde{\chi}_1|_{\ell
+ \nu /2}^2 - C \ t^{-2} \ h_1^{-1} \ |\widetilde{\chi}_1|_{\ell + \nu /2}\ 
|\chi_{t_0}|_{\bar{\ell} + \nu /2} \ |\phi |_{\widetilde{\ell}} \eqno(6.28)$$
\noi under the conditions
$$\widetilde{\ell} + \bar{\ell} > \ell + n/2 - \nu \quad , \quad \widetilde{\ell} \geq
\ell + 1 - \nu /2 \quad , \quad \bar{\ell} \geq \ell + 1 - \nu \quad , \quad \ell + 1
\geq \nu /2$$
\noi which also follow from (6.20). Defining now $z_{(1)} = z_{(1)} (\chi_1 ;[T, t_0],
h_1, \ell)$, integrating over time and using (6.21) (6.27), we obtain in the same way
as before
$$z^2 \vee z_1^2 \leq C\ b\ \beta \ z_1\ \mathrel{\mathop {\rm Sup}_{t}} \left \{ t^{-2} \
\rho '^{-1} \ h_0 \ h_1^{-1} \right \}  = C\ b\ \beta \ z_1$$
\noi and therefore
$$Z \left ( \chi_{t_0} - \psi_+;[T, t_0], h_1, \ell \right ) \leq C \ b \ \beta \ .
\eqno(6.29)$$
Starting from the basic estimates (6.27) and (6.29), the end of the proof is the same
as that of Proposition 6.1 based on (6.13) (6.16), with the simplification that the
initial condition at $t_0$ is given by a fixed $\psi_+$ instead of a time dependent
$W$. The difference between two solutions $\chi_{t_1}$ and $\chi_{t_0}$ with $T <
t_0 < t_1$ is estimated with the help of the extension of (6.27) to that difference
and of the pointwise part of (6.29) with $t_0$ replaced by $t_1$ taken at time $t =
t_0$ as
$$Z \left ( \chi_{t_1} - \chi_{t_0};[T, t_0], 1, \ell \right ) \leq C
\left | \chi_{t_1}(t_0) - \psi_+ \right |_{\ell} \leq C \ b \ \beta \ h_1(t_0)
\eqno(6.30)$$
\noi from which the existence of $\chi$ with the properties and convergences stated in
Part (2) follow. Taking the limit $t_1 \to \infty$ in (6.30) yields (6.24), while
taking the limit $t_0 \to \infty$ in (6.27) (6.29) yields (6.22) (6.23) in the same way
as in the proof of Proposition 6.1. \par

\noi \underbar{Part (3)} follows from elementary estimates together with the estimate
on $\phi$ expressed by (6.21). \par
   
\noi \underbar{Part (4)}. It follows from (6.3) and (6.19) that $V \exp (-i \chi )$
also satisfies (6.3), with gauge invariant initial condition $V(\infty ) \exp (-i
\chi (\infty )) = w_+ \exp (-i \psi_+)$. The result then follows from the uniqueness
statement of Proposition 6.1, part (3). \par \nobreak
\hfill $\sq$ \par 

\noi {\bf Remark 6.2.} Because of the linearity of the equations (6.3) (6.19), the
solutions $V$ and $\chi$ constructed in Propositions 6.1 and 6.2 have obvious
continuity properties with respect to $w_+$ and $\psi_+$ respectively. The continuity
with respect to $\phi$ is more delicate and will not be considered here. \par

We shall use the results of Propositions 6.1 and 6.2 in the special case where $\phi =
\phi_{p-1}$ and $W = W_p$ as defined by (5.22). In that case, $V$ satisfies an
asymptotic estimate which is much more accurate than (6.9) and which shows that $V$ is
a good approximation to $W_p$. In order to state that result, we need the results of
Proposition 5.2. We need in particular the special function $\bar{h}_0$ defined by
(5.7) and the associated functions $\bar{N}_m$ and $\bar{Q}_m$ associated with it
according to (3.68) (3.69). The result can be stated as follows. \\

\noi {\bf Proposition 6.3.} {\it Let $p \geq 1$ be an integer. Let $(k_+, \bar{k}, k)$
satisfy
$$k_+ > n/2\ , \ \widetilde{\ell} > n/2 - \nu \ , \ \widetilde{k} \geq k + 1 - \nu /2
\ , \ \widetilde{k} \geq \bar{k} \geq k + 1 - \nu \ , \ k \geq \nu /2 \eqno(6.31)$$
\noi where $\widetilde{\ell} (\equiv \ell_{p-1}) = k_+ - \lambda - (p -1)
\bar{\lambda}$ and $\widetilde{k} (\equiv k_p) = k_+ - p \bar{\lambda}$. Let $w_+ \in
K_{\rho_{_\infty}}^{k_+}$ and let $a_+ = |w_+|_{k_+}$. Let $\bar{h}_1$ and $h_2$ be ${\cal
C}^1$ positive nonincreasing functions defined in $[1, \infty )$ and tending to zero
at infinity, with $\bar{h}_1 \geq t^{-2} \rho'^{-1} \bar{h}_0$ and $h_2 \geq t^{-2}
\rho'^{-1} \bar{N}_p$. Let $\phi = \phi_{p-1}$ and $W = W_p$ be defined by (5.22), and
let 
$$b = \ \mathrel{\mathop {\rm Sup}_{t}} \bar{h}_0 (t)^{-1} \ |\phi
|_{\widetilde{\ell}} \ \equiv \ \mathrel{\mathop {\rm Sup}_{t}} \bar{h}_0(t)^{-1} \
|\phi_{p-1}(t)|_{\ell_{p-1}} \ . \eqno(6.32)$$
\noi Let $T, 1 < T < \infty$ satisfy
$$b \ \bar{h}_1 (T) \leq c \eqno(6.33)$$
\noi for a suitable constant $c$ (see (6.7)) and let $V$ be the solution of the
equation (6.3) constructed in Proposition 6.1, and such that $(V, 0) \in {\cal
X}_{\rho}^{\bar{k},0}([T, \infty ))$. Then $V$ satisfies the following estimate}
$$Y\left ( V - W_p;[T, \infty ), h_2, k \right ) \leq A(a_+) \ . \eqno(6.34)$$ 
\vskip 3 truemm

\noi {\bf Proof.} One checks easily that the assumptions of Proposition 6.3 imply the
relevant assumptions of Propositions 5.2 and 6.1. In particular (6.31) implies (6.4)
and the exponents $\widetilde{k}$ and $\widetilde{\ell}$ are those given by
Proposition 5.2. Let now $V_{t_0}$ be defined in Proposition 6.1, part (2), let $r =
V_{t_0} - W_p$ so that in particular $r(t_0) = 0$, and $r$ satisfies the equation 
$$\partial_t r = \left ( 2t^2\right )^{-1} \left \{ \left ( 2 \nabla \phi_{p-1}
\cdot \nabla + (\Delta \phi_{p-1}) \right ) r + \sum_{i\leq p - 1, j\leq p\atop{i+j
\geq p}} \left ( 2 \nabla \varphi_i \cdot \nabla + (\Delta \varphi_i)\right ) w_j
\right \} \eqno(6.35)$$ 
\noi obtained by taking the difference between (6.3) and the appropriate sum of
(2.18). Let $\widetilde{r} = h_2^{-1}r$. We now estimate $\partial_t
|\widetilde{r}|_k^2$. The contribution of the terms containing $r$ in the RHS are
estimated exactly as in the proof of Proposition 6.1. The remaining terms are
estimated by (3.23) (3.24) with $m = k - \nu /2$, $k \to \widetilde{k}$ and $\ell \to
\widetilde{\ell}$ under the conditions $\widetilde{\ell} > n/2 - \nu$, $0 \leq k -
\nu /2 \leq \widetilde{k} - 1$, which follow from (6.31). We obtain
$$\partial_t \ |\widetilde{r}|_k^2 \geq 2 \rho '\  |\widetilde{r}|_{k + \nu /2}^2 -
C \ t^{-2} \ |\widetilde{r}|_{k + \nu /2}^2 \ |\phi |_{\widetilde{\ell}} - C\ t^{-2}
\ h_2^{-1} \ |\widetilde{r}|_{k + \nu /2} \sum |\varphi_i|_{\widetilde{\ell}}\
|w_j|_{\widetilde{k}}$$
$$\geq 2 \rho ' |\widetilde{r}|_{k + \nu /2}^2 - C \ b \ t^{-2} \ \bar{h}_0(t) \
|\widetilde{r}|_{k + \nu /2}^2 - A(a_+) \ t^{-2} \ h_2^{-1} \ \bar{N}_p \
|\widetilde{r}|_{k + \nu /2} \eqno(6.36)$$
\noi by (6.32) (5.9) (5.10) (3.70) and (3.79). Defining $y_{(1)} = y_{(1)} (r,
[T, t_0], h_2, k)$ and integrating (6.36) over time with the initial condition
$r(t_0) = 0$, we obtain in the same way as before
$$y^2 \vee y_1^2 \leq C\ b\ \bar{h}_1(t) \ y_1^2 + A(a_+) \ y_1$$
\noi by using the properties of $\bar{h}_1$ and $h_2$, and therefore, under the
condition (6.33)
$$Y\left ( r; [T, t_0], h_2, k \right ) \leq A(a_+) \ . \eqno(6.37)$$
\noi We now take $T < t_0 < t_1$, and take the limit $t_1 \to \infty$, $t_0 \to
\infty$ in that order in the estimate
$$Y\left ( V_{t_1} - W; [T, t_0], h_2, k \right ) \leq A(a_+) $$     
\noi which follows from (6.37). The estimate (6.34) then follows from the convergence
(6.10). \par \nobreak
\hfill $\sq$ \par

\noi {\bf Remark 6.2.} In order to appreciate the improvement of the asymptotic
accuracy of Proposition 6.3, especially (6.34) over Proposition 6.1, especially (6.9),
it is useful to consider again the special case $\rho ' = t^{-1 - \varepsilon}$.
Proposition 6.1 with $\phi = \phi_{p-1}$ is now applied with $h_0 = \bar{h}_0 \sim
t^{1-\gamma}$ and therefore $h_1 = \bar{h}_1 \sim t^{-\gamma + \varepsilon}$, so that
the pointwise part of (6.9) states that
$$|V(t) - w_+|_k \ \lsim \ C \ ab\  \ t^{-\gamma +
\varepsilon} \ .$$
\noi On the other hand $\bar{N}_p \sim t^{1-(p+1)\gamma}$ for $(p+1)\gamma < 1$ thereby allowing
for $h_2 \sim t^{-2} \rho'^{-1} \bar{N}_p \sim t^{-(p+1)\gamma + \varepsilon}$ in Proposition 6.3,
so that the pointwise part of (6.34) states that
$$|V(t) - W_p(t)|_k \ \lsim \ A(a_+) \ t^{-(p+1)\gamma + \varepsilon} \ .$$

That improvement will play an essential role in the estimates of Proposition 6.4 below. \\

We now turn to the construction of solutions $(w, \varphi )$ of the system (2.11)-(2.12) with
given asymptotic states $(w_+, \psi_+)$. For that purpose, we first take a large positive
$t_0$ and we construct a solution $(w_{t_0}, \varphi_{t_0})$ of (2.11)-(2.12) with initial
data $(V(t_0), \phi_p (t_0) + \chi (t_0))$ at $t_0$. The solution $(w, \varphi )$ will be
obtained therefrom by taking the limit $t_0 \to \infty$, as explained in Section 2. \\

\noi {\bf Proposition 6.4.} {\it Let $(k, \ell)$ satisfy (3.48) and $k \geq 1 - \nu
/2$, let $p$ be an integer such that $(p+2)\gamma > 1$ and let $k_+$ and $\ell_+$
satisfy 
$$k_+ > n/2 \ , \ k_+ \geq (k + 2 - \nu ) \vee (\ell + \lambda + 1) + p \bar{\lambda}\ ,
\ \ell_+ \geq \ell + 1 \ , \eqno(6.38)$$
\noi where $\lambda = \mu - n + 2$ and $\bar{\lambda} = \lambda \vee 1$. Let
$\bar{h}_0$ be defined by (5.7) and let $\bar{N}_m$, $\bar{Q}_m$ be the associated
estimating functions defined by (3.68) (3.69). Let $\bar{h}_1$, $h_1$, $h_2$ and $h_3$
be positive nonincreasing ${\cal C}^1$ functions defined in $[1, \infty )$, tending
to zero at infinity, and satisfying 
$$\bar{h}_1 \geq t^{-2} \ \rho '^{-1} \ \bar{h}_0 \ , \ h_2 \geq t^{-2} \ \rho '^{-1}\ 
\bar{N}_p \ , \ h_2 \geq \bar{Q}_p \ {\it if} \ p \geq 1 \ , \eqno(6.39)$$
$$h_3 \geq t^{-\gamma} \ \rho '^{-1} \ h_2 \ , \ h_1 \geq t^{-2} \ \rho '^{-1} \ h_3 \
h_2^{-1} \ , \ h_3 \geq C\ \bar{h}_1 \ . \eqno(6.40)$$
\noi Let $w_+ \in K_{\rho_{_\infty}}^{k_+}$ and let $(W_m, \phi_m)$, $0 \leq m \leq p$,
be defined by (5.22) so that $(W_m, \bar{h}_0^{-1} \phi_m) \in {\cal
Y}_{\rho_{_\infty}}^{k_m,\ell_m}([1, \infty ))$ by Proposition 5.2. Let $V$ be defined
by Proposition 6.1 with $(W, \phi ) = (W_p, \phi_{p-1})$ and $(\widetilde{k} ,
\widetilde{\ell}) = (k_p, \ell_{p-1})$. Let $\psi_+ \in Y_{\rho_{_\infty}}^{\ell_+}$
and let $\chi$ be defined by Proposition 6.2 with the same $(W, \phi )$. Let
$$a_+ = |w_+|_{\rho_{_\infty},k_+} \quad , \quad b_+ = |\psi_+|_{\rho_{_\infty},
\ell_+} \ .$$
\noi Then there exists $T$, $1 \leq T < \infty$, depending only on $(\gamma ,p, a_+,
b_+)$ such that for all $t_0 \geq T$, the system (2.11)-(2.12) with initial data
$w(t_0) = V(t_0)$, $\varphi (t_0) = \phi_p (t_0) + \chi (t_0)$ has a unique solution
$(w_{t_0}, \varphi_{t_0}) \in {\cal X}_{\rho}^{k, \ell} ([T,t_0])$. One can define
$T$ by a condition of the type
$$A(a_+, b_+) \left ( h_1(T) \vee \bar{h}_1(T) \vee h_2(T) \right ) = 1 \eqno(6.41)$$
\noi and the solution satifies the estimates}
$$Y\left ( w_{t_0} - V ; [T, t_0], h_2, k \right ) \vee Y\left ( w_{t_0} - W_p;[T,
t_0],h_2,k \right ) \leq A(a_+, b_+) \eqno(6.42)$$
  $$Z\left ( \varphi_{t_0} - \phi_p - \chi  ; [T, t_0], h_3, \ell \right ) \vee Z\left (
\varphi_{t_0} - \phi_p - \psi_+;[T, t_0],h_3, \ell \right ) \leq A(a_+, b_+)
\eqno(6.43)$$ 
$$Y\left ( w_{t_0} ; [T, t_0], 1, k \right ) \leq A(a_+, b_+) \ , \eqno(6.44)$$
$$Z \left ( \varphi_{t_0}; [T, t_0], \bar{h}_0 , \ell \right ) \leq A(a_+, b_+) \ .
\eqno(6.45)$$
\vskip 3 truemm

\noi {\bf Remark 6.3.} As mentioned previously, in the same way as in Propositions 6.1
and 6.2, we could easily (but we shall not) extend the solution $(w_{t_0},
\varphi_{t_0})$ to the interval $[T, \infty )$. \\

\noi {\bf Remark 6.4.} In order to understand the time decay estimates implied by
(6.42) (6.43), it is useful to consider again the special case $\rho ' = t^{-1 -
\varepsilon}$. Saturating as far as possible the inequalities in (6.39) (6.40), we
obtain 
$$\bar{h}_1 \sim t^{-\gamma + \varepsilon} \ , \ h_2 \sim t^{-(p+1)\gamma +
\varepsilon} \ , \ h_3 \sim t^{1-(p+2)\gamma + 2 \varepsilon} \ , \ h_1 = t^{-\gamma
+ 2 \varepsilon} $$
\noi for $(p+1)\gamma < 1$, and the assumptions are satisfied for $\varepsilon$
sufficiently small, namely $2 \varepsilon < (p+2) \gamma - 1$. In particular the
condition that $h_3$ be decreasing in $t$ essentially imposes the condition
$(p+2)\gamma > 1$.  \par

Note also that in the condition (6.41) $h_2$ yields only a marginal restriction as
soon as $p \geq 1$, while for $p = 0$ it is natural to take $h_2 = \bar{h}_1$ (see
(6.39) with $\bar{N}_0 = \bar{h}_0)$. Finally, the condition $h_3 \geq C \bar{h}_1$
will in general be automatically satisfied for any reasonable choice of the
estimating functions $\bar{h}_1$ and $h_3$. \\

\noi {\bf Proof.} The proof follows the same pattern as that of Proposition 4.1,
involving a parabolic regularization, possibly a regularization of the initial data,
the local resolution of the regularized system by a fixed point method, the
derivation of a priori estimates uniform in the regularization, and a limiting
procedure. The only difference lies in the a priori estimates of the solutions in
${\cal X}_{\rho}^{k,\ell} ([T, t_0])$. Those estimates are much more elaborate than
previously and will ensure in particular that $T$ can be taken independent of $t_0$,
contrary to what happened in Lemma 4.1 (see especially (4.13)). We concentrate on
the proof of those estimates, omitting the parabolic regularization terms for
brevity. Their contribution will be briefly discussed at the end of that proof. \par

Let $(w_{t_0}, \varphi_{t_0}) \in {\cal X}_{\rho}^{k,\ell } ([T, t_0])$ be a
solution of the system (2.11)-(2.12) with initial condition $(w_{t_0}(t_0),
\varphi_{t_0}(t_0)) = (V(t_0), \phi_p(t_0) + \chi (t_0))$ where $V$ and $\chi$ are
defined by Propositions 6.1 and 6.2. Instead of estimating $(w_{t_0},
\varphi_{t_0})$ directly, we estimate the differences $q = w_{t_0} - V$ and $\psi
= \varphi_{t_0} - \phi_p - \chi$. For convenience we also introduce the gradients
$\sigma = \nabla \psi$, $s_{t_0} = \nabla \varphi_{t_0}$, $\tau = \nabla \chi$, as
well as $s_m = \nabla \varphi_m$, $S_m = \nabla \phi_m$ for $0 \leq m \leq p$.
Comparing the equations (2.11)-(2.12) and (6.3) (6.19) with $\phi = \phi_{p-1}$,
we obtain 
$$\partial_tq = \left ( 2t^2\right )^{-1} \Big \{ i \Delta w_{t_0} +
2s_{t_0}\cdot \nabla q + 2 (\sigma + s_p + \tau ) \cdot \nabla V + (\nabla \cdot
\sigma ) w_{t_0} + (\nabla \cdot ( S_p + \tau )) q + (\nabla \cdot (s_p + \tau
))V \Big \} \eqno(6.46)$$      
$$\partial_t\psi = \left ( 2t^2\right )^{-1} \Big \{ |\sigma|^2 + 2 \sigma \cdot ( S_p
+ \tau) + |\tau|^2 + 2 \tau \cdot s_p + \sum_{i,j\leq p\atop{i+j\geq p}} s_i \cdot s_j
\Big \}$$
$$+ t^{-\gamma} \Big \{ g_0(q,q) + 2g_0(q, V) + q_0(V- W_p, V + W_p) + \sum_{i,j \leq
p\atop{i+j>p}} g_0 (w_i, w_j) \Big \} \ . \eqno(6.47)$$
\noi It is convenient to write also the equation for $\sigma = \nabla \psi$, namely
$$\partial_t\sigma = t^{-2} \Big \{ s_{t_0} \cdot \nabla \sigma + \sigma \cdot \nabla
(S_p + \tau ) + (\tau + s_p) \cdot \nabla \tau + \tau \cdot \nabla s_p + \sum_{i,j\leq
p\atop{i+j \geq p}} s_i\cdot \nabla s_j \Big \}$$
$$+ t^{-\gamma} \ \nabla \Big \{ g_0(q,q) + 2g_0(q, V) + q_0(V- W_p, V + W_p) + \sum_{i,j
\leq p\atop{i+j>p}} g_0 (w_i, w_j) \Big \} \ . \eqno(6.48)$$
\noi We define $\widetilde{q} = h_2^{-1}q$ and $\widetilde{\psi} = h_3^{-1}\psi$ and
we estimate $\partial_t|\widetilde{q}|_k^2$ and $\partial_t|\widetilde{\psi}|_{\ell}^2$
by exactly the same method as in Lemmas 3.5, 3.6, 3.7, based on Lemma 3.4, using in
particular (3.22) (3.23) (3.24) for $\widetilde{q}$ and (3.25) (3.26) (3.28) for
$\widetilde{\psi}$, and omitting the terms with $h'_2$ and $h'_3$. We obtain for $t
\leq t_0$
$$\partial_t|\widetilde{q}|_k^2 \geq 2 \rho '\ |\widetilde{q}|_{k+ \nu /2}^2 - C \
t^{-2} \ |\widetilde{q}|_{k+ \nu /2} \Big \{ h_2^{-1}\  |V|_{k+2- \nu /2}$$
$$+ |\widetilde{q}|_{k+ \nu /2} \ |\psi + \phi_p + \chi |_{\ell} + h_2^{-1} \
|V|_{k+1 - \nu /2} \ |\psi + \varphi_p + \chi |_{\ell}$$ 
$$+ \left ( |\widetilde{q}|_k + h_2^{-1} \ |V|_k \right ) |\psi |_{\ell + \nu /2} +
|\widetilde{q}|_k \ |\phi_p + \chi|_{\ell + \nu /2} + h_2^{-1} \ |V|_k \ |\varphi_p
+ \chi |_{\ell + \nu /2} \Big \} \eqno(6.49)$$   
$$\partial_t|\widetilde{\psi}|_\ell^2 \geq 2 \rho '\ |\widetilde{\psi}|_{\ell+ \nu /2}^2 -
C \ t^{-2} \ |\widetilde{\psi}|_{\ell + \nu /2} \Big \{ |\widetilde{\psi}|_{\ell + \nu /2} \ 
|\psi + \phi_p + \chi |_{\ell}$$
$$+ |\widetilde{\psi}|_{\ell} \ |\phi_p +  \chi |_{\ell + 1 - \nu /2} + h_3^{-1} \
|\varphi_p + \chi |_{\ell} \ |\chi|_{\ell + 1 - \nu /2} + h_3^{-1} \ |\chi|_{\ell} \
|\varphi_p|_{\ell + 1 - \nu /2}$$
$$+ h_3^{-1} \sum_{i,j\leq p\atop{i+j \geq p}} |\varphi_i|_{\ell} \ |\varphi_j|_{\ell +
1 - \nu /2} \Big \} - C \ t^{-\gamma} \ h_3^{-1} \ |\widetilde{\psi}|_{\ell+ \nu /2}
\Big \{ h_2 |\widetilde{q}|_{k+ \nu /2}$$
$$\left ( |q|_k + |V|_k \right ) + h_2 |\widetilde{q}|_k \ |V|_{k + \nu /2} + |V -
W_p|_k \ |V + W_p|_{k + \nu/2}$$
$$+ |V - W_p|_{k + \nu /2} \ |V + W_p|_k + \sum_{i,j \leq p\atop{i+j>p}} |w_i|_k \
|w_j|_{k + \nu /2} \Big \} \ . \eqno(6.50)$$
\noi In order to continue the estimates, we need some information on $V$, $\chi$ and on the
$(\varphi_m, w_m)$. Applying Proposition 6.1 with $\phi = \phi_{p-1}$ and Proposition
6.3, both with $\bar{k} = k + 2 - \nu$, we rewrite (6.8) and (6.34) as
$$Y\left ( V;[T, \infty ), 1, k+2 - \nu \right ) \leq a \eqno(6.51)$$
$$Y\left ( V- W_p;[T, \infty ), h_2, k\right ) \leq a \eqno(6.52)$$
\noi for some $a$ depending on $a_+$, under the conditions
$$k_+ > n/2 \ , \ \ell_{p-1} > n/2 - \nu \ , \ k_p \geq k + 2 - \nu \eqno(6.53)$$
\noi where (cf (5.8))
$$k_m = k_+ - m \ \bar{\lambda} \quad , \quad \ell_m = k_+ - \lambda - m\ \bar{\lambda}
\ .$$
\noi Similarly, applying Proposition 6.2 with $\phi = \phi_{p-1}$, $\bar{\ell} = \ell +
1$, $h_0 = \bar{h}_0$ and therefore $h_1 = \bar{h}_1$, we rewrite (6.22) (6.23) as
$$Z\left ( \chi ;[T, \infty ), 1, \ell +1\right ) \leq b \eqno(6.54)$$
$$Z\left ( \chi - \psi_+;[T, \infty ), \bar{h}_1, \ell \right ) \leq b \eqno(6.55)$$
\noi for some $b$ depending on $a_+$ and $b_+$, under the conditions
$$\ell_{p-1} > n/2 - \nu \ , \ \ell_{p-1} \geq \ell + 2 - \nu /2 \ , \ \ell_+ \geq
\ell + 1 \ . \eqno(6.56)$$
\noi Finally, from Proposition 5.2, we obtain
$$|w_j(t)|_{k + \nu /2} \leq a \ \bar{Q}_{j-1}(t)    \qquad {\rm for} \ 1 \leq j \leq p \
, \eqno(6.57)$$
$$|W_p(t)|_{k + \nu /2} \leq a \ , \eqno(6.58)$$
\noi for some $a$ depending on $a_+$, under the conditions
$$k_+ > n/2 \quad , \quad k_p \geq k + \nu /2 \ , \eqno(6.59)$$
\noi and
$$|\varphi_j(t)|_{\ell + 1 - \nu /2} \leq b \ \bar{N}_j(t) \qquad {\rm for} \ 0 \leq j
\leq p \ , \eqno(6.60)$$
$$|\phi_p(t)|_{\ell + 1 - \nu /2} \leq b\ \bar{h}_0 (t) \eqno(6.61)$$
\noi for some $b$ depending on $a_+$, under the conditions
$$k_+ > n/2 \quad , \quad \ell_p \geq \ell + 1 - \nu /2 \ . \eqno(6.62)$$
\noi In the estimates (6.51) (6.52) (6.54) (6.55) (6.57) (6.58) (6.60) (6.61), we use
two common letters $a$ and $b$ to refer to estimates of amplitudes and phases
respectively. The constant $a$ depends only on $a_+$, while $b$ depends on $a_+$ and on
$b_+$. The conditions (6.53) (6.56) (6.59) (6.62) are implied by (6.38), which is the
statement of 
$$k_+ > n/2 \ , \ k_p \geq k + 2 - \nu \ , \ \ell_p \geq \ell + 1 \ , \ \ell_+ \geq
\ell + 1 \ .$$

With the estimates (6.51)-(6.61) available, we continue to estimate $(q, \psi )$ by
integrating (6.49) (6.50) over time in the interval $[t, t_0]$. We define $y_{(1)} =
y_{(1)} (q;[t,t_0], h_2,k)$ and $z_{(1)} = z_{(1)}(\psi ; [t, t_0], h_3, \ell )$. We
proceed exactly as in the proofs of Lemmas 4.1, 4.2 and 4.3, using the Schwarz
inequality for the time integrals whenever necessary. We use furthermore the fact that
(6.49) (resp. (6.50)) contains $|\widetilde{q}|_{k + \nu /2}$ (resp.
$|\widetilde{\psi}|_{\ell + \nu /2})$ as a factor in its RHS, thereby yielding a
factor $y_1$ (resp. $z_1$) after integration, and the elementary fact that 
$$y^2 \vee y_1^2 \leq A\ y_1 \Rightarrow y \vee y_1 \leq A$$
\noi and its analogue for $(z, z_1)$. We then obtain the following estimates, where
Sup means that the Supremum of the function of time that follows is taken in the
interval $[t, t_0]$, and where an overall constant $C$ is omitted for brevity. 
$$y \vee y_1 \leq a \ {\rm Sup} \left ( t^{-2} \ \rho'^{-1}\ h_2^{-1} \right ) + ab \
{\rm Sup} \left ( t^{-2} \ \rho'^{-1} \ \bar{N}_p \ h_2^{-1} \right )$$
$$+ b(y + y_1) \ {\rm Sup} \left ( t^{-2} \ \rho'^{-1} \ \bar{h}_0 \right ) + a (z +
z_1) \ {\rm Sup} \left ( t^{-2} \ \rho '^{-1} \ h_2^{-1} \ h_3 \right )$$
$$+ ( y z_1 + y_1 z ) \ {\rm Sup} \left ( t^{-2} \ \rho '^{-1} \ h_3 \right ) \ ,
\eqno(6.63)$$
$$z \vee z_1 \leq b^2 \ {\rm Sup} \left ( t^{-2} \ \rho'^{-1}\ \bar{h}_0 \ h_3^{-1}
\bar{N}_p \right ) + b(z + z_1)  \ {\rm Sup} \left ( t^{-2} \ \rho'^{-1} \ \bar{h}_0 \
\right)$$
$$+ z \ z_1 \ {\rm Sup} \left ( t^{-2} \ \rho'^{-1} \ h_3 \right ) + a (y +
y_1) \ {\rm Sup} \left ( t^{-\gamma} \ \rho '^{-1} \ h_3^{-1} \ h_2 \right )$$
$$+  y \ y_1  \ {\rm Sup} \left ( t^{-\gamma } \ \rho '^{-1} \ h_3^{-1} \ h_2^2 \right )
+ p a^2 \ {\rm Sup} \left ( t^{-\gamma} \ \rho'^{-1} \ h_3^{-1} (h_2 + \bar{Q}_p)
\right ) \eqno(6.64)$$
\noi where the factor $p$ in the last term simply means that that term is absent for
$p = 0$. The various Sup in time are estimated in an obvious way with the help of the
conditions (6.39) (6.40) which are taylored for that purpose. The only non-obvious
term is the coefficient of $b^2$ in (6.64), which is estimated by (6.39) (6.40) as  
$$\left ( t^{-2} \ \rho'^{-1} \ \bar{N}_p \ h_2^{-1} \right ) \left ( t^{-\gamma} \
\rho'^{-1} \ h_2 \ h_3^{-1} \right ) \left ( t^{\gamma} \ \rho ' \ \bar{h}_0 \right )
\leq \ \parallel \rho '\parallel_1$$
\noi since 
$$\bar{h}_0 = \int_1^t dt_1 \ t_1^{-\gamma} \ \rho '(t_1)^{-1} \ \rho '(t_1) \leq
t^{-\gamma} \ \rho '(t)^{-1} \int_1^t dt_1 \ \rho '(t_1) \leq \parallel \rho '
\parallel_1 \ t^{-\gamma} \ \rho '^{-1} \eqno(6.65)$$
\noi by the monotony of $t^{-\gamma} \rho '^{-1}$. Absorbing the factor $\parallel
\rho ' \parallel_1$ in the (again omitted) overall constant and defining as
previously $Y = y \vee y_1$ and $Z = z \vee z_1$, we end up with
$$Y \leq a + ab + b \ Y \ \bar{h}_1 + a\ Z\ h_1 + Y \ Z \ h_1 \ h_2 \eqno(6.66)$$
$$Z \leq b^2 + pa^2 + b \ Z\ \bar{h}_1 + Z^2 \ h_1 \ h_2 + a\ Y + Y^2 \ h_2
\eqno(6.67)$$
\noi where the functions $h_1$, $\bar{h}_1$ and $h_2$ are taken at time $t$ where
they take their Supremum in $[t, t_0]$, since they are assumed to be decreasing.\par

In order to conclude the estimates, we impose the conditions
$$4 \ b \ \bar{h}_1(T) \leq 1 \ , \ 4(1 + b)\ h_2(T) \leq 1 \ , \ Z\ h_1 (t) \leq ( 1 +
b) \eqno(6.68)$$
\noi which together imply $4 Z \ h_1(t) \ h_2(t) \leq 1$ for all $t \geq T$. We then
obtain
$$Y \leq 2a (1 + b) + 2a \ Z\ h_1 \leq 4a(1 + b) \eqno(6.69)$$
$$Z \leq 2(b^2 + pa^2) + 2a \ Y + 2Y(Y \ h_2) \leq 2(b^2 + pa^2) + 4a Y \leq 2(b^2 +
pa^2) + 16a^2(1 + b) \ .\eqno(6.70)$$
\noi The condition $Zh_1 \leq 1 + b$ then reduces to
$$\left ( 2(b^2 + pa^2) (1+b)^{-1} + 16a^2 \right ) \ h_1 \leq 1$$
\noi which is implied by
$$\left ( 2b + (16 + 2p) a^2 \right ) \ h_1(T) \leq 1 \ . \eqno(6.71)$$
\noi The condition (6.71) together with the first two conditions of (6.68) then take
the form (6.41), while the estimates (6.69) (6.70) yield the estimate of the first
terms in the LHS of (6.42) (6.43). The estimates of the second terms follow from those of the
first ones and from (6.52) and (6.55), together with the condition $h_3 \geq C
\bar{h}_1$ from (6.40). \par

Finally the estimates (6.44) (6.45) follow immediately from (6.42) (6.43), from (6.51)
or (6.58) and from (6.61) (6.54). \par 
   
We now discuss briefly the contribution of the parabolic regularization terms in the
previous proof. We regularize the system (2.11)-(2.12) in the same way as in the proof
of Proposition 4.1 (see (4.51), where however the sign of the regularizing terms should
be changed since we are now solving the equations for decreasing $t$ starting from
$t_0$). Instead of (6.49) (6.50), we then obtain
$$\partial_t \ |\widetilde{q}|_k^2 \geq 2 \theta \left ( |\nabla \widetilde{q}|_k^2 +
h_2^{-1} \ {\rm Re} \ <\nabla \widetilde{q}, \nabla V>_k \right ) + \hbox{previous
terms} \ ,$$ 
$$\partial_t \ |\widetilde{\psi}|_{\ell}^2 \geq 2 \theta \left ( |\nabla
\widetilde{\psi}|_{\ell}^2 + h_3^{-1} \  <\nabla \widetilde{\psi}, \nabla (\phi_p +
\chi )>_{\ell} \right ) + \hbox{previous terms} $$ 
\noi where $<\cdot , \cdot >_k$ and $<\cdot , \cdot >_{\ell}$ denote the scalar
products in $K_{\rho}^k$ and $Y_{\rho}^{\ell}$. The scalar products are controlled
since we have assumed that $(V, \phi_p + \chi ) \in {\cal X}_{\rho , loc}^{k+1, \ell +
1}$ by imposing $\bar{k} = k + 2 - \nu$, $\bar{\ell} = \ell + 1$ and $\ell_p \geq \ell
+ 1$. They produce additional terms in the final estimates which are uniformly bounded
in $\theta$ in a neighborhood of zero. \par

With the a priori estimates (6.44) (6.45) replacing Lemma 4.1, the proof of
Proposition 6.4 proceeds in the same way as that of Proposition 4.1, as mentioned
above. In particular, the required regularity estimates  and difference estimates are
provided by Parts 3 of Lemmas 4.2 and 4.3, which are taylored for that purpose. The
assumption on $(w, \varphi )$ made in those parts follow from (6.44) (6.45) and from
the relation 
$$\bar{h}_0 \leq \ \parallel \rho '\parallel_1 \ h_0$$
\noi which follows from (6.65), while the regularity assumptions on the initial data
required in Lemma 4.2 are ensured by the previous regularity of $V$ and $\phi_p +
\chi$. \par \nobreak
\hfill $\sq$ \par

We can now take the limit $t_0 \to \infty$ of the solution $(w_{t_0}, \varphi_{t_0})$
constructed in Proposition 6.4, for fixed $(w_+, \psi_+)$. \\

\noi {\bf Proposition 6.5.} {\it Let the assumptions of Proposition 6.4 be satisfied.
Then\par

(1) There exists $T$, $1 \leq T < \infty$, depending only on $(\gamma , p, a_+, b_+)$
and there exists a unique solution $(w, \varphi )$ of the system (2.11)-(2.12) in the
interval $[T, \infty )$ such that $(w, \bar{h}_0^{-1} \varphi ) \in {\cal
X}_{\rho}^{k,\ell} ([T, \infty ))$ and such that the following estimates hold
$$Y\left ( w - V;[T, \infty ), h_2, k\right ) \vee Y\left ( w - W_p;[T, \infty ),
h_2, k\right ) \leq A(a_+, b_+) \ , \eqno(6.72)$$
$$Z\left ( \varphi - \phi_p - \chi ; [T, \infty ), h_3, \ell \right ) \vee Z \left (
\varphi - \phi_p - \psi_+;[T, \infty ), h_3, \ell \right ) \leq A(a_+, b_+) \ ,
\eqno(6.73)$$
$$Y \left ( w; [T, \infty ),1,k \right ) \leq A(a_+, b_+) \ , \eqno(6.74)$$
$$Z \left ( \varphi ; [T, \infty ), \bar{h}_0, \ell \right ) \leq A(a_+, b_+) \ .
\eqno(6.75)$$
\noi One can define $T$ by a condition of the type (6.41). \par

(2) Let $(w_{t_0}, \varphi_{t_0}) \in {\cal X}_{\rho}^{k,\ell}([T, t_0])$ be the
solution of the system (2.11)-(2.12) constructed in Proposition 6.4. Then $(w_{t_0},
\varphi_{t_0})$ converges to $(w, \varphi )$ in norm in ${\cal X}_{\rho}^{k',\ell
'}([T, T_1])$ for $k' < k$, $\ell ' < \ell$ and in the weak-$*$ sense in ${\cal
X}_{\rho}^{k,\ell}([T, T_1])$ for any $T_1$, $T < T_1 < \infty$, and in the weak-$*$
sense in $K_{\rho}^k \oplus Y_{\rho}^{\ell}$ pointwise in $t$ for all $t \geq T$. \par

(3) The map $(w_+, \psi_+) \to (w, \varphi )$ defined in Part (1) is continuous from
the norm topology of $(w_+, \psi_+)$ in $K_{\rho_{_\infty}}^{k_+} \oplus
Y_{\rho_{_\infty}}^{\ell_+}$ to the norm topology of $(w, \varphi )$ in ${\cal
X}_{\rho}^{k',\ell '}([T, T_1])$ for $k' < k$, $\ell ' < \ell$ and to the weak-$*$
topology in ${\cal X}_{\rho}^{k,\ell}([T, T_1])$ for any $T_1$, $T < T_1 < \infty$,
and to the weak-$*$ topology in $K_{\rho}^k \oplus Y_{\rho}^{\ell}$ pointwise in $t$
for all $t \geq T$.}\\

\noi {\bf Remark 6.5.} For simplicity, we have not stated the strongest continuity
properties that would follow by tracking more accurately the exponents in the proof
of Part (3). Actually the required topology on $(w_+, \psi_+)$ could be weakened to
the norm topology of $K_{\rho_{_\infty}}^{k'} \oplus Y_{\rho_{_\infty}}^{\ell '}$ on
the bounded sets of $K_{\rho_{_\infty}}^{k_+} \oplus Y_{\rho_{_\infty}}^{\ell_+}$
for suitable $(k', \ell ')$ smaller than $(k_+, \ell_+)$. \\

\noi {\bf Proof.} Parts (1) and (2) will follow from the convergence of $(w_{t_0},
\varphi_{t_0})$ when $t_0 \to \infty$ in the topologies stated in Part (2). Let $T
\leq t_0 \leq t_1$, let $(w_{t_0}, \varphi_{t_0})$ and $(w_{t_1}, \varphi_{t_1})$ be
the corresponding solutions of the system (2.11)-(2.12) obtained in Proposition 6.4,
and let $(w_-, \varphi_-) = (w_{t_0} - w_{t_1}, \varphi_{t_0} - \varphi_{t_1})$.
From (6.42) (6.43) and their analogues for $(w_{t_1}, \varphi_{t_1})$, it follows
that      
$$\left \{ \begin{array}{l} Y\left ( w_-;[T, t_0], h_2, k \right ) \leq A(a_+, b_+) \\
\\ Z \left ( \varphi_-;[T,t_0], h_3, \ell \right ) \leq A(a_+, b_+) \end{array} \right .
\eqno(6.76)$$
\noi so that in particular
$$\left \{ \begin{array}{l} \left | w_-(t_0) \right |_k \leq A(a_+, b_+) h_2(t_0) \\ \\
\left | \varphi_- (t_0)\right |_{\ell} \leq A(a_+, b_+) h_3(t_0) \ . \end{array} \right .
\eqno(6.77)$$

We now apply Lemma 4.3, part (3) to $(w_-, \varphi_-)$. For that purpose we take $h_0 =
t^{-\gamma} \rho'^{-1}$, so that by (6.40)
$$h_1 \geq t^{-2} \ \rho '^{-1} \ h_0 \qquad , \qquad h_2 \ h_0 \leq h_3 \eqno(6.78)$$
\noi and in particular $(h_0, h_1)$ satisfy the assumptions of Lemma 4.3. The
assumptions (4.35) (res\-tric\-ted to the relevant interval $[T, t_0])$ and (4.38) follow
from (6.44) (6.45) and from (6.65), while the condition (4.26) can be included in
(6.41). We now apply (4.41) with $k' = k - 1 + \nu$ $\ell ' = \ell - 1 + \nu$, together
with (6.77) (6.78), thereby obtaining 
$$Y \left ( w_-;[T,t_0], h_0^{-1}, k - 1 + \nu \right ) \leq A(a_+,b_+) \Big \{
h_0(t_0)|w_-(t_0)|_{k - 1 + \nu} + h_1(T) |\varphi_-(t_0)|_{\ell - 1 + \nu} \Big \}$$
$$ \leq
A(a_+, b_+) h_3(t_0) \eqno(6.79)$$
$$Z \left ( \varphi_-;[T,t_0], 1, \ell  - 1 + \nu \right ) \leq A(a_+,b_+) \Big \{
|\varphi_-(t_0)|_{\ell - 1 + \nu} + h_0(t_0) |w_-(t_0)|_{k - 1 + \nu} \Big \}$$
$$ \leq
A(a_+, b_+) h_3(t_0) \ . \eqno(6.80)$$
\noi From (6.79) (6.80) and from the fact that $h_3$ tends to zero at infinity, it
follows that there exists $(w, \varphi ) \in {\cal X}_{\rho , loc}^{k-1 + \nu , \ell -
1 + \nu} ([T, \infty ))$ such that $(w_{t_0}, \varphi_{t_0})$ converges to $(w,
\varphi )$ in norm in ${\cal X}_{\rho}^{k - 1 + \nu , \ell - 1 + \nu}([T, T_1 ])$ for all
$T_1$, $T < T_1 < \infty$. From that convergence, from (6.42)-(6.45) and from standard
compactness, continuity and interpolation arguments, it follows that $(w,
\bar{h}_0^{-1}\varphi ) \in {\cal X}_{\rho}^{k , \ell }([T, \infty ))$, that $(w,
\varphi )$ satisfies the estimates (6.72)-(6.75) and that $(w_{t_0}, \varphi_{t_0})$
converges to $(w, \varphi )$ in the other topologies stated in Part (2). Furthermore
$(w, \varphi )$ satisfies the system (2.11)-(2.12), and uniqueness of $(w, \varphi )$
under the conditions (6.72) (6.73) follows from Proposition 4.3 and from the fact that
$h_3(t)$ tends to zero at infinity.\par

\noi \underbar{Part (3)}. Let $(w_+, \psi_+)$ and $(w'_+, \psi '_+)$ belong to a
fixed bounded set of $K_{\rho_{_\infty}}^{k_+} \oplus
Y_{\rho_{_\infty}}^{\ell_+}$, so that (5.17) (5.18) and (6.72)-(6.75) hold with
fixed $A$. Let $(W_p, \phi_p)$ and $(W'_p, \phi '_p)$ be the associated
functions defined by (5.22) and let $(w, \varphi )$ and $(w', \varphi ')$ be the
associated solutions of the system (2.11)-(2.12) defined in Part (1). We assume
that $(w'_+, \psi '_+)$ is close to $(w_+, \psi_+)$ in the sense that 
$$\left | w_+ - w'_+ \right |_{k_+} \leq \varepsilon \qquad , \quad \left |
\psi_+ - \psi '_+ \right |_{\ell} \leq \varepsilon_0 \ . \eqno(6.81)$$
\noi Let $w_- = w - w'$, $\varphi_- = \varphi - \varphi '$ and let $t_0 > T$ be
defined by $h_2(t_0) = \varepsilon$, which can be done for any (sufficiently
small) $\varepsilon > 0$ and which implies that $t_0 \to \infty$ when
$\varepsilon \to 0$. It follows from (6.72) and (5.17) that
$$|w_-(t_0)|_k \leq A \ h_2(t_0) + |W_p(t_0) - W'_p(t_0)|_k \leq A \left ( h_2(t_0) +
\varepsilon \right ) \leq A \ h_2(t_0) \eqno(6.82)$$ 
\noi and from (6.73) and (5.18) that
$$\begin{array}{ll} |\varphi_-(t_0)|_{\ell} &\leq Ah_3(t_0) + |\phi_p(t_0) -
\phi '_p(t_0)|_{\ell} + |\psi_+ - \psi '_+|_{\ell} \\ \\ &\leq A\left (
h_3(t_0) + \bar{h}_0(t_0) \varepsilon \right ) + |\psi_+ - \psi '_+|_{\ell}
\leq Ah_3(t_0) + \varepsilon_0 \ . \end{array} \eqno(6.83)$$

We now apply Lemma 4.3, part (3) with $k' = k - 1 + \nu$, $\ell ' = \ell - 1 + \nu$,
thereby obtaining
$$\begin{array}{ll} Y\left ( w_-;[T, t_0],h_0^{-1},k-1 + \nu \right ) &\leq A \left (
h_0(t_0) \ h_2(t_0) + h_1(T) \right . \left ( h_3(t_0) + \varepsilon_0 )\right ) \\ \\
&\leq A \left ( h_3 (t_0) + \varepsilon_0 \right ) 
\end{array} \eqno(6.84)$$ 
$$\begin{array}{ll} Z \left ( \varphi_-;[T,t_0], 1, \ell - 1 + \nu \right )
&\leq A \left \{ h_3(t_0) + \varepsilon_0 + h_0(t_0) h_2(t_0) \right \} \\ \\ 
&\leq A \left ( h_3(t_0) + \varepsilon_0 \right ) \ . \end{array}
\eqno(6.85)$$ \noi When $\varepsilon$ tends to zero, $h_3(t_0)$ tends to zero, and
therefore (6.84) (6.85) imply norm continuity in ${\cal X}_{\rho}^{k-1 + \nu , \ell - 1
+ \nu}([T, T_1])$ for all $T_1$, $T < T_1 < \infty$. The other continuities follow
therefrom, from the estimates (6.74) (6.75) and from standard continuity, interpolation
and compactness arguments. \par \nobreak \hfill $\sq$

\section{Asymptotics and wave operators for $u$}
\hspace*{\parindent} In this section we complete the construction of the wave operators for the
equation (1.1) and we derive asymptotic properties of solutions in their range. The construction
relies in an essential way on those of Section 6, especially Proposition 6.5, and will
require a discussion of the gauge invariance of those constructions. This section follows
closely Section II.7. \par

We first define the wave operator for the auxiliary system (2.11)-(2.12).\\

\noi {\bf Definition 7.1.} We define the wave operator $\Omega_0$ as the map
$$\Omega_0 : (w_+, \psi_{+} ) \to (w, \varphi ) \eqno(7.1)$$
\noi from $K_{\rho_{_\infty}}^{k_+} \oplus Y_{\rho_{_\infty}}^{\ell_+}$ to the space of $(w,
\varphi )$ such that $(w, \bar{h}_0^{-1} \varphi ) \in {\cal X}_{\rho}^{k,\ell}([T, \infty))$ for
some $T$, $1 \leq T < \infty$, where $\rho$, $k_+$, $\ell_+$, $k$, $\ell$ satisfy the assumptions
of Proposition 6.5 and $(w, \varphi )$ is the solution of the system (2.11)-(2.12) obtained
in Part (1) of that proposition. \\

In order to study the gauge invariance of $\Omega_0$, we need some information on the Cauchy
problem at finite times for the equation (1.1). We define the operator $J \equiv J(t) = x + it
\nabla$, which satisfies the commutation relation
$$i \ M \ D \ \nabla = J \ M \ D\ , \eqno(7.2)$$
\noi where $M$ and $D$ are defined by (2.4) (2.5). For any interval $I \subset [1, \infty )$,
any nonnegative integer $k$ and any nonnegative ${\cal C}^1$ function $\rho$ defined in $I$, we
define the space
$$\begin{array}{ll} X_{\rho}^k(I) &= \left \{ u: D^*M^*u \in {\cal C} (I, K_{\rho}^k)
\right \} \\ \\
& = \left \{ u : <J(t)>^k \ f(J(t)) u \in {\cal C}(I, L^2) \right \} \ .\end{array} \eqno(7.3)$$  
\noi where $<\lambda > = (1 + \lambda^2)^{1/2}$ for any real number or self-adjoint operator
$\lambda$ and where the second equality follows from (7.2) and from Remark 3.1. (The space
$X_0^k(I)$ was denoted ${\cal X}^k(I)$ in II). We recall the following result (see Proposition
I.7.1). \\

\noi {\bf Proposition 7.1.} {\it Let $k$ be a positive integer and let $0 < \mu < 2k$.
Then the Cauchy problem for the equation (1.1) with initial data $u(t_0) = u_0$ such that
$<J(t_0)>^k$ $u_0 \in L^2$ at some initial time $t_0 \geq 1$ is locally well posed in
${\cal X}_0^k(\cdot )$, namely} \par
{\it (1) There exists $T > 0$ such that (1.1) has a unique solution with initial data $u(t_0) = u_0$
in ${\cal X}^k_0 ([1 \vee (t_0 - T), t_0 + T])$.} \par
{\it (2) For any interval $I$, $t_0 \in I \subset [1 , \infty )$, (1.1) with initial data $u(t_0) =
u_0$ has at most one solution in ${\cal X}_0^k(I)$.} \par
{\it (3) The solution of Part (1) depends continuously on $u_0$ in the norms considered there.} \\

We come back from the system (2.11)-(2.12) to the equation (1.1) by reconstructing $u$ from
$(w, \varphi )$ by (2.7) and accordingly we define the map
$$\Lambda : (w, \varphi ) \to u = M \ D \exp (-i \varphi ) w \quad . \eqno(7.4)$$

\noi It follows immediately from Lemma 3.3 that the map $\Lambda$ satisfies the following
property.\\

\noi {\bf Lemma 7.1.} {\it Let $\ell + 2 > n/2$ and $0 \leq k \leq \ell + 2$. Then for any
interval $I \subset [1, \infty )$ and any nonnegative ${\cal C}^1$ function $\rho$ defined in
$I$, the map $\Lambda$ is bounded and continuous from ${\cal Y}_{\rho , loc}^{k, \ell}(I)$
(defined by (5.1) (6.2)) to $X_{\rho}^k(I)$, with norm estimates in compact intervals independent
of $\rho$.} \\

We now give the following definition\\

\noi {\bf Definition 7.2.} Two solutions $(w, \varphi )$ and $(w', \varphi ')$ of the system
(2.11)-(2.12) in ${\cal Y}_{\rho, loc}^{k, \ell}(I)$ for some $k$, $\ell$, $\rho$ and some
interval $I \subset [1, \infty )$ are said to be gauge equivalent if $\Lambda (w, \varphi ) =
\Lambda (w', \varphi ')$, or equivalently if
$$\exp (- i \ \varphi (t)) \ w(t) = \exp (- i \ \varphi '(t)) \ w'(t) \eqno(7.5)$$
\noi for all $t \in I$. \par

The following sufficient condition for gauge equivalence follows immediately from Lemma 7.1
and from the uniqueness statement of Proposition 7.1, part (2). \\

\noi {\bf Lemma 7.2.} {\it Let $\ell + 2 > n/2$ and $0 \leq k \leq \ell + 2$. Let $I \subset [1,
\infty )$ be an interval, let $\rho$ be a strictly positive ${\cal C}^1$ function defined in
$I$, and let $(w, \varphi )$ and $(w', \varphi ')$ be two solutions of the system
(2.11)-(2.12) in ${\cal Y}_{\rho , loc}^{k, \ell}(I)$. In order that $(w, \varphi )$ and
$(w', \varphi ')$ be gauge equivalent, it is sufficient that (7.5) holds for one $t \in I$. }
\\

We now turn to the study of the gauge equivalence of asymptotic states (if any) for solutions of
the system (2.11)-(2.12) such as those obtained in Proposition 4.1. For that purpose, we need
$\rho$ to be decreasing, namely to be defined by (4.1) for some $t_0 \geq 1$. \\

\noi {\bf Proposition 7.2.} {\it Let $k \geq 0$ and $\ell > n/2 - \nu$. Let $\rho$ be defined by
(4.1) for some $t_0 \geq 1$ and let $h_0$ and $h_1$ be as in Proposition 4.1. Let $(w, \varphi )$
and $(w', \varphi ')$ be two gauge equivalent solutions of the system (2.11)-(2.12) such that
$(w, h_0^{-1} \varphi )$, $(w', h_0^{-1} \varphi ') \in {\cal X}_{\rho}^{k,\ell} ([T, \infty
))$ for some $T$, $t_0 \leq T < \infty$. Let 
$$b = Z\left ( \varphi ; [T, \infty ),h_0, \ell \right ) \vee Z \left ( \varphi ' ;[T, \infty ),
h_0 , \ell \right ) \ . \eqno(7.6)$$
\noi Then \par

(1) There exists $\omega \in Y_{\rho(\infty)}^{\ell - 1 + \nu}$ such that $\varphi '(t) - \varphi
(t)$ converges to $\omega$ strongly in $Y_{\rho ( \infty)}^{\ell '}$ for $\ell ' < \ell - 1 + \nu$
and weakly in $Y_{\rho (\infty)}^{\ell - 1 + \nu}$. The following estimate holds~:
$$\left |\varphi '(t) - \varphi (t) - \omega \right |_{\ell - 2 + \nu} \leq C \ b\left |\varphi
'(T_1) - \varphi (T_1) \right |_{\ell - 1 + \nu} \ h_1(t) \eqno(7.7)$$
\noi for $t \geq T_1$, for some $T_1$ sufficiently large, namely satisfying
$$b\ h_1(T_1) \leq c \eqno(7.8)$$
\noi for some constant $c$. \par

(2) Assume in addition that $1 - \nu /2 \leq k \leq \ell + 1$ and let $w_+$ and $w'_+$ be the
limits of $w(t)$ and $w'(t)$ as $t \to \infty$ obtained in Proposition 4.2. Then $w'_+ =
w_+ \exp (- i \omega )$. \par

(3) Let $p \geq 0$ be an integer. Assume in addition that $w_+$, $w'_+ \in K_{\rho (\infty
)}^{k_+}$ where $k_+$ satisfies
$$k_+ > n/2 \qquad , \qquad k_+ \geq (p+1)\bar{\lambda} - 1 \eqno(7.9)$$
\noi and let $\phi_p$, $\phi '_p$ be associated with $(w_+, w'_+)$ according to (5.22) and
Proposition 5.2. Assume that the following limits exist
$$\lim_{t \to \infty} \left ( \varphi (t) - \phi_p (t) \right ) = \psi_+ \quad , \quad \lim_{t \to
\infty} \left ( \varphi '(t) - \phi'_p(t)\right ) = \psi '_+ \eqno(7.10)$$
\noi as strong limits in $L^{\infty}$. Then $\psi '_+ = \psi_+ + \omega$.} \\

\noi {\bf Proof.} \underbar{Part (1)}. Let $\varphi_{\pm}(t) = \varphi '(t) \pm \varphi (t)$. By
the same estimates as in the proof of Lemma 3.7, we obtain
$$\partial_t \ |\varphi_-|_{\ell '}^2 - 2 \rho ' \ |\varphi_-|_{\ell ' + \nu /2}^2 \leq C \ t^{-2}
\ |\varphi _-|_{\ell ' + \nu /2} \left \{ | \varphi_-|_{\ell ' + \nu /2} \ |\varphi_+|_{\ell} +
|\varphi_-|_{\ell '} \ |\varphi_+|_{\ell + \nu /2} \right \}  \eqno(7.11)$$ \noi for $\nu /2
\leq \ell ' + 1 \leq \ell + \nu$. Defining $z_{(1)} = z_{(1)}(\varphi_-;[T_1,t], 1, \ell ')$
with $T \leq T_1 < t$, integrating over time and using (7.6), we obtain
$$\begin{array}{ll}
z^2 \vee z_1^2 &\leq  z_0^2 + C \left ( \displaystyle{\mathrel{\mathop {\rm Sup}_{[T_1,t]}}}\ 
t^{-2} \ |\rho '|^{-1} \ h_0 \right ) b\ z_1(z+ z_1) \\ \\ &\leq  z_0^2 + C\ b\ h_1 (T_1) \ z_1(z +
z_1) \end{array}$$
\noi where $z_0 = |\varphi_- (T_1)|_{\ell '}$, which under the condition (7.8) with suitable $c$,
yields
$$Z\left ( \varphi_-;[T_1, \infty ) , 1, \ell ' \right ) \leq C \ |\varphi_-(T_1)|_{\ell '} \ .
\eqno(7.12)$$
\noi From
$$\partial_t \ \varphi_- = \left ( 2t^2 \right )^{-1} \left ( \nabla \varphi_- \cdot \nabla
\varphi_+ \right )$$
\noi we next estimate directly
$$\begin{array}{ll} |\partial_t \ \varphi_-|_{\ell '-1} &\leq C\ t^{-2}\
|\varphi_-|_{\ell '} \ |\varphi_+|_{\ell} \\ \\ &\leq C \ t^{-2} \ h_0 \ b
\ |\varphi_-(T_1)|_{\ell '}\end{array}\eqno(7.13)$$ 
\noi by (7.6) and (7.12). The last member of (7.13) is integrable in time since
$$\int_t^{\infty} dt_1 \ t_1^{-2} \ h_0(t_1) \leq \ \parallel \rho '\parallel_1 \ h_1(t) \ ,
\eqno(7.14)$$
\noi which for $\ell ' = \ell - 1 + \nu$ proves the existence of the limit $\omega \in Y_{\rho
(\infty )}^{\ell - 2 + \nu}$ for $\varphi_-$ together with the estimate (7.7). The fact that
actually $\omega \in Y_{\rho (\infty )}^{\ell - 1 + \nu}$ and the other convergences stated in
Part (1) follow therefrom and from (7.12) by standard compactness and interpolation arguments. \par

 \underbar{Parts (2) and (3)}. The proof is identical with that of the corresponding
statements in Proposition II.7.2 and will be omitted.\par \nobreak
\hfill $\sq$ \par

\noi {\bf Remark 7.1.} The assumptions made on $(k, \ell )$ in Parts (1) and (2) of Proposition
7.2 are implied by the assumptions (3.48) and $k \geq 1 - \nu /2$ of Proposition 4.1, so that
Parts (1) and (2) apply directly to the solutions of the system (2.11)-(2.12) constructed in
that proposition. In Part (3), the assumption (7.9) is that required in Proposition 5.2 to
ensure the existence and appropriate estimates of $\phi_p$, $\phi '_p$. That assumption,
together with the existence of the limits (7.10), is ensured under the assumptions of
Proposition 5.3, namely if $(k , \ell )$ satisfy in addition (5.24) for some $k_0$
satisfying (5.25) and if $(p + 2)\gamma > 1$ and $P_p(1) < \infty$.\\

Proposition 7.2 prompts us to make the following definition. \\

\noi {\bf Definition 7.3.} Two pairs of asymptotic states $(w_+, \psi_+)$ and $(w'_+, \psi' _+)$
are gauge equivalent if $w_+ \exp (-i \psi_+) = w'_+ \exp (- i \psi '_+)$. \\

With this definition, Proposition 7.2 implies that two gauge equivalent solutions of the system
(2.11)-(2.12) in ${\cal R}(\Omega_0)$ are images of two gauge equivalent pairs of asymptotic
states. One should however not overlook the following fact. The wave operator $\Omega_0$ is
defined through Proposition 6.5 which uses an \underbar{increasing} $\rho$ defined by (6.1)
whereas Proposition 7.2 uses a \underbar{decreasing} $\rho$ (essentially in the proof of (7.12)).
In order to apply Proposition 7.2 to the solutions constructed in Proposition 6.5 with increasing
$\rho$, one has therefore to take some large $t_0$, to define
$$\widetilde{\rho}(t) = \rho (t_0) - \int_{t_0}^t \rho '(t_1) \ dt_1$$
\noi and to apply Proposition 7.2 with that new $\widetilde{\rho}$, thereby ending up with
information on $(\omega, w_+, w'_+,$ $\psi_+, \psi '_+)$ in spaces $K$ or $Y$ associated with
$\widetilde{\rho}(\infty ) = \rho (\infty) - 2 \int_{t_0}^{\infty} \rho '(t) dt < \rho (\infty
)$. That fact of course does not impair the algebraic relations expressing gauge invariance. \par

We now turn to the converse result, namely to the fact that gauge equivalent asymptotic states
generate gauge equivalent solutions through $\Omega_0$. \\

\noi {\bf Proposition 7.3.} {\it Let $(k, \ell )$ and $(k_+ , \ell_+)$ satisfy the assumptions
of Proposition 6.5, namely (3.48), $k \geq 1 - \nu /2$ and (6.38). Let $\rho$ and the
estimating functions of time also satisfy the assumptions of Proposition 6.5. Let $(w_+,
\psi_+)$, $(w'_+, \psi '_+) \in K_{\rho_{_\infty}}^{k_+} \oplus Y_{\rho_{_\infty}}^{\ell_+}$
be gauge equivalent and let $(w, \varphi)$, $(w', \varphi ')$ be their images under
$\Omega_0$. Then $(w, \varphi )$ and $(w', \varphi ')$ are gauge equivalent.}\\

\noi {\bf Proof.} The proof is identical with that of Proposition II.7.3. We reproduce it for
completeness. \par

Let $t_0$ be sufficiently large and let $(w_{t_0}, \varphi_{t_0})$ and
$(w'_{t_0}, \varphi '_{t_0})$ be the solutions of the system (2.11)-(2.12) constructed by
Proposition 6.4. From the initial conditions
$$w_{t_0}(t_0) = V(t_0) \qquad , \qquad w'_{t_0} (t_0) = V'(t_0) \ ,$$
$$\varphi_{t_0}(t_0) = \phi_p(t_0) + \chi (t_0) \quad , \quad \varphi '_{t_0}(t_0) = \phi '_p(t_0)
+ \chi ' (t_0)  , $$
\noi from the fact that $\phi_p = \phi '_p$ by Proposition 5.2, part (2) and that $V \exp (-
i \chi ) = V' \exp (- i \chi ')$ by Proposition 6.2, part (4), it follows that
$$w_{t_0}(t_0) \exp  (- i \varphi_{t_0}(t_0) ) = w'_{t_0} (t_0) \exp ( - i
\varphi '_{t_0} (t_0) )$$
\noi and therefore by Lemma 7.2, $(w_{t_0}, \varphi_{t_0})$ and $(w'_{t_0}, \varphi '_{t_0})$ are
gauge equivalent, namely
$$w_{t_0} (t) \exp  ( - i \varphi_{t_0}(t)  ) = w'_{t_0}(t) \exp ( - i \varphi
'_{t_0}(t)  ) \eqno(7.15)$$

\noi for all $t$ for which both solutions are defined. \par

We now take the limit $t_0 \to \infty$ for fixed $t$ in (7.15). By Proposition 6.5, part
(2), for fixed $t$, $(w_{t_0}, \varphi_{t_0})$ and $(w'_{t_0}, \varphi '_{t_0})$ converge
respectively to $(w, \varphi )$ and $(w', \varphi ')$ in $K_{\rho}^{k'} \oplus
Y_{\rho}^{\ell '}$. By Lemma 3.3, one can take the limit $t_0 \to \infty$ in (7.15), thereby
obtaining (7.5), so that $(w, \varphi )$ and $(w', \varphi ')$ are gauge equivalent. \par
\nobreak \hfill $\sq$ \par

We can now define the (local) wave operators (at infinity) for $u$. \\

\noi {\bf Definition 7.4.} The wave operator $\Omega$ is defined as the map 
$$\Omega : u_+ \to u = (\Lambda \circ \Omega_0) (Fu_+,0) \eqno(7.16)$$
\noi where $\Omega_0$ and $\Lambda$ are defined by Definition 7.1 and by (7.4). \\

We collect in the following proposition the information on $\Omega$ that follows from the
previous study, in particular from Propositions 6.5 and 7.3. \\

\noi {\bf Proposition 7.4.} {\it Let $p \geq 0$ be an integer with $(p+ 2)\gamma > 1$. Let
$\rho$ (defined by (6.1)) and the estimating functions of time satisfy the assumptions of
Proposition 6.5 (see especially (6.39) (6.40) and Remark 6.4). Let
$$\lambda = \mu - n + 2 \leq 2 \nu \ , \eqno(7.17)$$
$$k \geq 1 - \nu /2 \qquad , \quad k > 1 - \nu + \mu /2
\ . \eqno(7.18)$$
\noi Let $k_+$ satisfy (6.38) for some $\ell$ satisfying $\ell > n/2 - \nu$, $\ell \geq k -
\nu$. Then\par

(1) The wave operator $\Omega$ maps $F K_{\rho_{_\infty}}^{k_+}$ to $X_{\rho}^k([T, \infty ))$
for some $T$, $1 \leq T < \infty$. (Actually $T$ depends on $u_+$).\par

(2) $\Omega$ is injective.} \\

\noi {\bf Proof.} \underbar{Part (1)} follows from the definitions, from Lemma 7.1 and from
Proposition 6.5 with $\psi_+ = 0$. The only point to be checked is the fact that (7.17)
(7.18) imply the existence of $\ell$ such that $(k, \ell)$ satisfies (3.48). Now the
$\ell$ dependent part of (3.48) reduces to
$$k - \nu \leq \ell \leq k + \nu - \lambda \eqno(7.19)$$
$$n/2 - \nu < \ell < 2k-\lambda + \nu - n/2 \eqno(7.20)$$
\noi and the compatibility of (7.19) (7.20) for $\ell$ is easily seen to reduce to (7.17)
and to the second inequality in (7.18). \par

 \underbar{Part (2)} follows from Proposition 7.2 and from the fact that a gauge
equivalence class of asymptotic states contains at most one element with $\psi_+ = 0$.\par
\nobreak
\hfill $\sq$  \\

\noi {\bf Remark 7.2.} One may wonder whether the restriction to asymptotic states with
$\psi_+ = 0$ in (7.16) restricts the range of $\Omega$ as compared with that of $\Lambda
\circ \Omega_0$. From Proposition 7.3, it follows that
$$(\Lambda \circ \Omega_0) (w_+, \psi_+) = (\Lambda \circ \Omega_0) (w_+ \exp (- i
\psi_+), 0)$$
\noi in so far as $w_+ \exp (- i \psi_+)$ has the regularity needed to apply Proposition
6.5. This is the case if $\ell_+ \geq k_+ - 2$ by Lemma 3.3. That condition however is
significantly stronger than the condition on $\ell_+$ contained in (6.38), especially for
large $p$, i.e. for small $\gamma$. Therefore, there is actually a restriction of the
range for regularity reasons, in spite of the (algebraic) gauge invariance of the
construction expressed by Proposition 7.3. \\

We next collect the information and in particular the asymptotic estimates obtained for
the solutions of the equation (1.1) in ${\cal R}(\Omega )$. \\

\noi {\bf Proposition 7.5.} {\it Let $0 < \mu \leq n  - 2 + 2 \nu \leq n$. Let $0 < \gamma \leq 1$
and let $p \geq 0$ be an integer with $(p+2) \gamma > 1$. Let $\rho$ (defined by 6.1) and the
estimating functions of time satisfy the assumptions of Proposition 6.5 (especially (6.39)
(6.40)). Let $(k, \ell , k_+)$ satisfy (7.18) (3.48) (6.38). Let $u_+ \in F
K_{\rho_{_\infty}}^{k_+}$ and $a_+ = \parallel Fu_+;K_{\rho_{_\infty}}^{k_+}\parallel$. Then \par

(1) There exists $T$, $1 \leq T < \infty$, and there exists a unique solution $u \in
X_{\rho}^k([T, \infty ))$ of the equation (1.1) which can be represented as 
$$u = M\ D \exp (- i \varphi ) w$$
\noi where $(w, \varphi )$ is a solution of the system (2.11)-(2.12) such that $(w,
\bar{h}_0^{-1} \varphi ) \in {\cal X}_{\rho}^{k, \ell}([T, \infty ))$ and such that
$$|w(t) - Fu_+|_{k-1 + \nu} \ h_0(t) \to 0 \eqno(7.21)$$
$$|\varphi (t) - \phi_p(t)|_{\ell - 1 + \nu} \to 0 \eqno(7.22)$$
\noi when $t \to \infty$. The time $T$ can be defined by a condition of the type (6.41) with $b_+
= 0$.\par

(2) The solution is obtained as $u = \Omega u_+$, following Definition 7.4. \par

(3) The map $\Omega$ is continuous from $FK_{\rho_{_\infty}}^{k_+}$ to the norm topology of
$X_{\rho}^{k'}(I)$ for $k' < k$ and to the weak-$*$ topology of $X_{\rho}^k(I)$ for any compact
interval $I \subset [T, \infty )$, and to the weak topology of $MDK_{\rho}^k$ pointwise in
$t$.\par

(4) The solution $u$ satisfies the following estimate for $t\geq T$~: 
$$\parallel <J(t)>^k \ f(J(t)) \left ( \exp [i \ \phi_p(t,x/t)] u(t) - M(t) D(t) Fu_+ \right )
\parallel_2 \ \leq A(a_+)\ h_3(t) \eqno(7.23)$$
\noi for some estimating function $A(a_+)$. \par

(5) Let $2 \leq r < \infty$. The solution $u$ satisfies the following estimate~:
$$\parallel u(t) - \exp [ - i \ \phi_p(t, x/t) ] M(t) \ D(t) \ Fu_+\parallel_r \ \leq A(a_+)
\ \rho (t)^{-\beta} \ t^{-\delta (r)} \ h_3 (t)  \eqno(7.24)$$ \noi where $\delta (r) = n/2 - n/r$,
$\beta = \nu^{-1} (\delta (r) - k) \vee 0$ if $r < \infty$ or $k > n/2$, $\beta = \nu^{-1} (n/2 - k
+ \varepsilon )$ if $r = \infty$ and $k \leq n/2$.  } \\

\noi {\bf Proof.} \underbar{Parts (1) (2) (3)} follow from Propositions 6.5 and 4.3, from
Definition 7.4 and from Proposition 7.4. \par

\underbar{Part (4)}. From the commutation relation (7.2), from Remark 3.1 and from Lemma
3.3, it follows that the LHS of (7.23) is estimated by
$$\parallel \cdot \parallel_2 \ \leq C \left | \exp (i (\phi_p - \varphi )) w- Fu_+\right |_k \leq
C \left \{ |w - Fu_+|_k + \left | \exp (i (\phi_p - \varphi )) - 1 \right |_{\ell} \ |w|_k \right
\}$$
$$\leq C \left \{ \left | w - W_p \right |_k + \left | W_p - w_+ \right |_k + \exp \left (C
\left | \phi_p - \varphi \right |_{\ell} \right ) \left | \phi_p - \varphi \right |_{\ell}
\ |w|_k \right \}$$
\noi and the result follows from the estimates (6.72) (5.9) (6.73) (6.74).\par

\underbar{Part (5)} follows from Part (4) and from the inequality
$$\begin{array}{ll}\parallel v\parallel_r \ &= t^{-\delta (r)} \parallel D^*\ M^*\ v \parallel_r \
\leq C\ t^{-\delta (r)}\ \parallel <\nabla>^m \ D^*\ M^*\ v \parallel_2 \\ \\ &= C\ t^{-\delta (r)}
\ \parallel <J(t)>^m v \parallel_2 \end{array}$$
\noi which follows from (7.2) and Sobolev inequalities with $m = \delta (r)$ if $r < \infty$, $m =
n/2 + \varepsilon$ if $r = \infty$, together with the fact that 
$$|\xi|^{m-k} \ f(\xi)^{-1} \leq C\ \rho^{-\beta} \ .$$
\hfill $\sq$ \\

\noi {\bf Remark 7.3.} In (7.23) and (7.24), one could replace $MDFu_+$ by $U(t)  u_+$ since the
difference is small in the relevant norms. One could also replace $Fu_+$ by $W_p$, but this would
not produce any improvement in the estimates, since the main contribution comes from the phase.
The estimate (7.24) is a rather weak one, since we have omitted the function $f$ which generates
the Gevrey regularity. It is only one example of a large number of similar estimates exhibiting
the typical $t^{-\delta (r)}$ decay associated with $L^r$ norms. \\

The final step of the standard construction of the wave operators for the equation (1.1) would
consist in extending the solutions $u$ to arbitrary finite times, and defining the maps
$\Omega_1 : u_+ \to u(1)$ where $u = \Omega u_+$. In order not to waste the Gevrey regularity
of the local solutions at infinity, this would require a treatment of the global Cauchy
problem at finite times for arbitrarily large data in the same Gevrey framework. This is a
rather different problem and we shall refrain from considering it here. \\

\noi {\bf Acknowledgements.} Part of this work was done while one of the authors (G.V.) was
visiting the Institut des Hautes Etudes Scientifiques (IHES), Bures-sur-Yvette, France and the
Laboratoire de Physique Th\'eorique (LPT), Universit\'e de Paris-Sud, France. He is very
grateful to Professor Jean-Pierre Bourguignon, Director of the IHES, and to Professor Dominique
Schiff, Director of the LPT, for the kind hospitality extended to him.

\section*{Appendix A}
\hspace*{\parindent} In this appendix we derive a number of properties of the function
$\widetilde{f}$ defined by (3.2), which make it a possible substitute for the function $f_0$
defined by (3.1) in the definition of the spaces ${\cal X}_{\rho}^{k,\ell}$. In all this
appendix, we assume $0 < \nu \leq 1$ and we take $\rho = 1$. The parameter $\rho$ can be
reintroduced easily by scaling. Accordingly we define 
$$\widetilde{f}(\xi ) = \sum_{j\geq 0} (j!)^{-1/\nu} \ |\xi |^j \ , \eqno({\rm A.1})$$
$$F(\xi ) = \sum_{j\geq 0} (j + 1)^{-1} \ (j!)^{-1/ \nu} \ |\xi|^{j+1} \ . \eqno({\rm A.2})$$
\noi We also use (A.1) and (A.2) to define $\widetilde{f}$ and $F$ when applied to $\xi \in
{I\hskip - 1 truemm R}^+$. With that convention, we have $\widetilde{f}(\xi ) =
\widetilde{f}(|\xi |)$ and $F(\xi ) = F(|\xi |)$ for all $\xi \in {I \hskip - 1 truemm R}^n$.
Furthermore $\widetilde{f} = dF/d|\xi |$. \par

We first derive some preliminary estimates which allow in particular for a comparison of
$\widetilde{f}$ and $f_0$. \\

\noi {\bf Lemma A.1.} {\it The following estimates hold for all $\xi \in {I \hskip - 1 truemm
R}^n$ :
$$\sum_{j\geq 1} j (j!)^{-1/\nu} \ |\xi |^j \leq |\xi |^{\nu} \sum_{j\geq 0} (j!)^{-1/ \nu}\ |\xi
|^j \leq \sum_{j \geq 0} (j+1) (j!)^{-1/ \nu} \ |\xi|^j \ , \eqno({\rm A.3})$$
$$\sum_{j\geq 0} (j + 1)^{-1}\  (j!)^{-1/\nu} \ |\xi |^{j+1} \leq |\xi |^{1 - \nu} \sum_{j\geq 0}
(j!)^{-1/ \nu}\ |\xi |^j  \ , \eqno({\rm A.4})$$
$$\widetilde{f}^{-1} (d\widetilde{f}/d|\xi |) \leq |\xi|^{\nu - 1} \leq F^{-1}\  \widetilde{f} =
F^{-1} (dF/d|\xi |) \ , \eqno({\rm A.5})$$
$$|\xi |^{\nu} \ \widetilde{f} \leq d (|\xi | \widetilde{f})/d|\xi | \ , \eqno({\rm A.6})$$
$$F(a) (|\xi| \vee a)^{\nu - 1} \ \exp \left ( \nu^{-1} (|\xi |^{\nu} - a^{\nu}) \right ) \leq
\widetilde{f}(\xi ) \leq \exp (\nu^{-1} \ |\xi |^{\nu}) \eqno({\rm A.7})$$
\noi for all $a > 0$,}
$$\widetilde{f}(\xi ) = (2 \pi )^{(\nu - 1)/2 \nu} \ \nu^{1/2} \ |\xi|^{(\nu - 1)/2} \ \exp
\left ( \nu^{-1}\  |\xi|^{\nu} \right ) (1 + o(1)) \quad {\it when} \ |\xi | \to \infty \ .
\eqno({\rm A.8})$$
\vskip 3 truemm

\noi {\bf Proof.} (A.3a) and (A.4) follow from the H\"older inequality on $\Z^+$ with the
measure $(j!)^{-1/ \nu} |\xi |^j $ and the exponents $1/ \nu$ and $1/(1 - \nu )$, applied
respectively to the pairs of functions $(j, 1)$ and $((j+1)^{-1}, 1)$. \par

(A.3b) follows similarly from the H\"older inequality with the measure
$(j+1)(j!)^{-1/\nu}|\xi|^j$ and the exponents $\nu^{-1}(1 + \nu )$ and $1 + \nu$, applied to
the functions $(j+1)^{-1}$ and 1. \par

(A.5) is a rewriting of (A.3a) and (A.4), while (A.6) is a rewriting of (A.3b).\par

(A.7) follows from (A.5). In fact (A.5) states that the functions $\widetilde{f}(\xi ) \exp (-
|\xi|^{\nu}/\nu )$ and $F(\xi ) \exp (-|\xi |^{\nu}/ \nu)$ are respectively decreasing (and
therefore less than one) and increasing in $|\xi |$. The first fact yields (A.7b) while both of
them together with $\widetilde{f} \geq |\xi|^{\nu - 1}F$ yield (A.7a). Note also that (A.7b)
follows directly from the definition (A.1) and from the fact that 
$$\widetilde{f}(\xi ) = \parallel (j!)^{-1}\ |\xi|^{\nu j} ; \ell^{1/\nu}\parallel^{1/\nu}
\ \leq \ \parallel (j!)^{-1}\ |\xi|^{\nu j} ; \ell^1 \parallel^{1/ \nu} \ = \exp (\nu^{-1} \
|\xi |^{\nu})$$
\noi by the standard embedding of $\ell^p$ spaces. \par

(A.8) is proved in \cite{15r} (see (8.07) p. 308) in the special case $\nu \geq 1/4$, but the
proof extends easily to the whole range $0 < \nu \leq 1$. \par \nobreak
\hfill $\sq$ \par

The estimates (A.7) and the asymptotic property (A.8) compare $\widetilde{f}$ with $f_0$ by
stating that essentially $\widetilde{f} \sim f_0^{1/ \nu}$. On the other hand (A.6) is the
analogue of the fact that $df_0/d|\xi | = \nu |\xi |^{\nu - 1}f_0$ and allows for the
construction of function space norms satisfying an inequality which can replace (3.12) in the
subsequent estimates. \par

We next show that $\widetilde{f}$ satisfies estimates similar to those of Lemma 3.1.\\

\noi {\bf Lemma A.2.} {\it Let $(\xi , \eta ) \in {I\hskip-1truemm R}^n$. Then $\widetilde{f}$
satisfies the estimates  
$$\widetilde{f}(\xi ) \leq \widetilde{f}(\xi - \eta ) \ \widetilde{f}(\eta
) \qquad {\it for \ all} \  (\xi , \eta ) \ , \eqno({\rm A.9})$$ 
$$\widetilde{f}(\xi ) \leq \widetilde{f}
(\xi - \eta ) \exp (|\eta |^{\nu}) \qquad {\it for} \ |\xi | \wedge |\eta | \leq |\xi - \eta | \ ,
\eqno({\rm A.10})$$
 $$(\widetilde{f}(\eta ) - \widetilde{f}(\xi )) \ |\eta |^{1 - \nu} \leq |\xi - \eta |
\widetilde{f}(\eta ) \qquad {\it for} \ |\xi | \leq |\eta | \ , \eqno({\rm A.11})$$
$$|\widetilde{f}(\xi ) - \widetilde{f}(\eta ) | \ |\eta |^{1-\nu} \leq |\xi - \eta |^{1 - \nu} \
\widetilde{f}(\xi - \eta ) \ \widetilde{f}(\eta ) \qquad {\it for \ all} \ (\xi , \eta ) \ ,
\eqno({\rm A.12})$$
$$|\widetilde{f}(\xi ) - \widetilde{f}(\eta )| \ |\eta |^{1 - \nu} \leq |\xi - \eta|^{1 - \nu}
(\exp (|\xi - \eta |^{\nu}) - 1) \widetilde{f}(\eta ) \qquad {\it for } \ |\xi| \wedge |\xi - \eta
| \leq |\eta | \ , \eqno({\rm A.13})$$
$$|\widetilde{f}(\xi ) - \widetilde{f}(\eta )| \ |\eta |^{1 - \nu} \leq C \ |\xi - \eta|^{1 - \nu}
(\widetilde{f}(\xi - \eta ) -1) \exp (|\eta |^{\nu}) \qquad {\it for } \ |\xi| \wedge |\eta
| \leq |\xi - \eta | \ . \eqno({\rm A.14})$$
\noi In (A.14) one can take $C = 1$, except in the region $|\xi | \leq |\xi - \eta | \leq
|\eta |$ where $C = 2^{1 - \nu}$.\par

The function $\widetilde{f}(\xi ) \vee \widetilde{f} (1)$ satisfies the same estimates as
$\widetilde{f} (\xi )$. } \\

\noi {\bf Proof.} (A.9) is trivial except if $|\xi | \geq |\eta | \vee |\xi - \eta |$. In
all cases, we estimate 
$$\begin{array}{ll} \widetilde{f}(\xi ) &\leq \displaystyle{\sum_{j,k \geq 0}} (j!)^{-1} \ (k!)^{-1}
\ ((j + k)!)^{1 - 1/\nu} \ |\xi - \eta |^j \ |\eta |^k\\
\\
&\leq \displaystyle{\sum_{j,k\geq 0}} (j!)^{-1 / \nu} \ (k!)^{-1/\nu} \ |\xi - \eta |^j \ |\eta |^k
= \widetilde{f} (\xi - \eta ) \ \widetilde{f}(\eta ) \end{array} \eqno({\rm A.15})$$ \noi since $(j
+ k)! \geq j! \ k!$. \par

(A.10) is trivial in the allowed region except if $|\eta | \leq |\xi - \eta | \leq |\xi |$.
In that case, we rewrite the first inequality in (A.15) as        
$$\begin{array}{ll}\widetilde{f}(\xi ) &\leq \displaystyle{\sum_{k\geq 0}} \left \{
\displaystyle{\sum_{j\geq 0}} \left ( (j!)^{-1/ \nu} \ |\xi - \eta |^j \right )^{\nu} \left
( ((j+k)!)^{-1/ \nu} \ |\xi - \eta|^{j+k} \right )^{1 - \nu} \right \}  (k!)^{-1} \
|\eta |^k\ |\xi - \eta |^{(\nu - 1) k}\\ \\ &\leq \widetilde{f}(\xi - \eta ) \exp \left (
|\eta | \ |\xi - \eta |^{\nu - 1} \right ) \leq \widetilde{f}(\xi - \eta ) \exp (|\eta
|^{\nu}) \end{array} \eqno({\rm A.16})$$ \noi by the H\"older inequality applied to the sum
over $j$ for fixed $k$ and the fact that $|\eta | \leq |\xi - \eta |$. \par

(A.11). We estimate
$$\widetilde{f}(\eta) - \widetilde{f}(\xi ) = \sum_{j\geq 1} (j!)^{-1 / \nu} \ (|\eta |^j - |\xi
|^j) \leq |\xi - \eta | \sum_{j\geq 1} (j!)^{-1/\nu}\ j |\eta|^{j-1} \leq |\xi - \eta | \
|\eta|^{\nu - 1} \ \widetilde{f}(\eta )$$
\noi by (A.3a). \par

(A.12) follows from (A.11) and from $|\xi - \eta |^{\nu} \leq \widetilde{f} (\xi - \eta )$ if
$|\xi | \leq |\eta |$. If $|\xi | \geq |\eta |$, we estimate
$$\widetilde{f}(\xi ) - \widetilde{f}(\eta ) \leq \sum_{j\geq 1,k\geq 0} (j!)^{-1} \ (k!)^{-1} \
((j+k)!)^{1- 1/ \nu} \ |\xi - \eta |^j \ |\eta |^k$$
$$= \sum_{j,k\geq 0} ((j+1)!)^{-1} \ (k!)^{-1} \ ((j+k+1)!)^{1 - 1/ \nu} \ |\xi - \eta |^{j+1} \
|\eta |^k$$
$$\leq \sum_{j,k\geq 0} (j+1)^{-1} \ (j!)^{-1/ \nu} \  |\xi - \eta |^{j+1} \
(k+1)((k+1)!)^{-1 / \nu} \ |\eta |^k$$ 
\noi since $(j+k+1)! \geq j!(k+1)!$,
$$\cdots \leq |\xi - \eta |^{1 - \nu} \ \widetilde{f}(\xi - \eta ) \ |\eta |^{\nu - 1} \
\widetilde{f}(\eta )$$
\noi by (A.3a) and (A.4). \par

(A.13) follows from (A.11) and $|\xi - \eta |^{\nu} \leq (\exp (|\xi - \eta |^{\nu}) - 1)$ if
$|\xi | \leq |\eta |$. If $|\xi - \eta | \leq |\eta | \leq |\xi |$, we estimate in the same way as
in (A.16)
$$(\widetilde{f}(\xi ) - \widetilde{f}(\eta )) \ |\eta |^{1 - \nu} \leq |\eta |^{1 - \nu}
\sum_{j\geq 1} \left \{ \sum_{k \geq 0} ((k!)^{-1/\nu} \ |\eta |^k )^{\nu} 
\left ( ((j+k)!)^{-1/\nu} \ |\eta |^{j+k} \right )^{1 - \nu} \right \}$$
$$\times  (j!)^{-1} \ |\xi -
\eta |^j \ |\eta |^{(\nu - 1)j}$$
$$\leq \widetilde{f} (\eta ) \Big \{ |\eta |^{1 - \nu} (\exp (|\xi - \eta | \ |\eta |^{\nu - 1} )
- 1 ) \Big \}  \eqno({\rm A.17})$$ \noi by the H\"older inequality applied to the sum over
$k$ for fixed $j$. Now for fixed $|\xi - \eta |$, the last bracket in (A.17) is a decreasing
function of $|\eta |$ and for $|\eta | \geq |\xi - \eta |$ is therefore less than its value
for $|\eta | = |\xi - \eta |$, which yields (A.13) in that case. \par

(A.14). If $|\xi| \vee |\eta | \leq |\xi - \eta |$, (A.14) with $C = 1$ follows from $|\eta |
\leq |\xi - \eta |$ and from 
$$|\widetilde{f}(\xi ) - \widetilde{f}(\eta )| \leq \widetilde{f} (\xi - \eta ) - 1 \ .$$
\noi If $|\eta | \leq |\xi - \eta | \leq |\xi |$, (A.14) with $C = 1$ follows from $|\eta | \leq
|\xi - \eta |$ and from a minor variant of (A.16) with the sum over $j$ restricted to $j \geq 1$,
so that  
$$\widetilde{f}(\xi ) - \widetilde{f}(\eta ) \leq (\widetilde{f}(\xi - \eta ) - 1) \exp (|\eta
|^{\nu}) \ . \eqno({\rm A.18})$$
\noi If $|\xi | \leq |\xi - \eta | \leq |\eta |$, (A.14) with $C = 2^{1- \nu}$ follows from $|\eta
| \leq 2 |\xi - \eta |$ and from (A.18) with $\xi$ and $\eta$ interchanged. \par

The last statement of Lemma A.2 is obvious as regards (A.9) (A.10) and follows from the fact that
$$|\widetilde{f}(\xi ) \vee a - \widetilde{f}(\eta ) \vee a | \leq |\widetilde{f} (\xi ) -
\widetilde{f} (\eta )|$$
\noi for all $\xi$, $\eta$ and $a > 0$ as regards (A.11)-(A.14), in the same way as in Lemma 3.1.
\par \nobreak
\hfill $\sq$\par

It follows from Lemma A.2 that $\widetilde{f} (\xi )$ and $\widetilde{f} (\xi ) \vee
\widetilde{f}(1)$ satisfy the basic estimates (A.9) and (A.12) which are used throughout this paper,
thereby making those functions into suitable substitutes for $f_0$ and $f$ in the definition of the
spaces. Note also that by (A.7) (A.8), $\exp (|\eta |^{\nu}) \sim \widetilde{f}(\eta )^{\nu}$ (up to
a small power of $\eta$), which makes (A.10) (A.13) (A.14) into close analogues of (3.4) (3.6)
(3.7).\par

We finally use the function $\widetilde{f}$ to relate the definition of spaces such as $K_{\rho}^k$
or $Y_{\rho}^{\ell}$ to more standard definition of Gevrey spaces. Since this part is meant
to be only illustrative, we restrict our attention to space dimension $n = 1$. The Gevrey
class $G_{1/\nu}$ can be defined as the vector space of ${\cal C}^{\infty}$ functions $u$
such that there exists a constant $C$ such that 
$$c_j \ C^{-j}\ (j!)^{-1/\nu} \ \partial^j u \in \ell_j^p (L_x^q) \eqno({\rm A.19})$$
\noi where $1 \leq p$, $q \leq \infty$ and $\{c_j\}$ is a sequence of positive numbers with at most
polynomial increase or decrease at infinity. Of course if one allows for all possible $C > 0$, the
parameters $p$, $q$, and $\{c_j\}$ are irrelevant, and one can take $p = q = \infty$ and $c_j
\equiv 1$, which yields the standard definition. Here however we fix $C$, in fact $C = \rho^{-1/
\nu}$, so as to obtain a Banach space, and those parameters become important. Since moreover we
want a Hilbert space in order to apply the energy method in a convenient way, we take $p = q=2$,
thereby obtaining the Gevrey Hilbert space $X$ with norm
$$\parallel u; X \parallel = \parallel b_j \ \rho^{j/\nu} \ \partial^j u;\ell_j^2(L_x^2)\parallel
\eqno({\rm A.20})$$
\noi where $b_j = c_j (j!)^{-1/ \nu}$. Let now
$$f(\xi ) = f(|\xi |) = \sum_{j\geq 0} a_j |\xi|^j = f_+ (\xi ) + f_- (\xi )$$
\noi where $\{a_j\}$ is a sequence of positive numbers and where $f_+$ and $f_-$ denote the sums
over even and odd $j$ respectively. We claim that for a suitable relation between $\{a_j\}$ and
$\{b_j\}$, the norm in $X$ is equivalent to the norm
$$\parallel u\parallel_* = \parallel f(\rho^{1/\nu} \xi ) \widehat{u} \parallel_2 \ .
\eqno({\rm A.21})$$
\noi In fact
$$f_+(\xi )^2 + f_-(\xi )^2 \leq f(\xi )^2 \leq 2\left (f_+(\xi )^2 + f_-(\xi )^2 \right )$$
\noi so that if we define $\{b_j\}$ by
$$f_+(\xi )^2 + f_-(\xi )^2 = \sum_{j\geq 0} b_j^2 \ |\xi |^{2j} \eqno(A.22)$$
\noi then
$$\parallel u; X\parallel \ \leq \ \parallel u \parallel_* \ \leq \sqrt{2} \parallel u; X \parallel
\eqno(A.23)$$
\noi which proves the equivalence. The relation (A.22) can be rewritten as
$$b_k^2 = \sum_{0 \leq j \leq 2k} a_j \ a_{2k-j} \ .$$
\noi In the special case of $\widetilde{f}$, where $a_j = (j!)^{-1/\nu}$, it is obvious that
$$a_j \leq b_j \leq \sqrt{2j+1} \ a_j$$
\noi and one can show that actually when $j \to \infty$
$$b_j = a_j (\pi \nu j)^{1/4} \ (1 + o(1))  \eqno(A.24)$$
\noi which together with (A.8) makes it possible to define equivalent norms of the form (A.20) for
the spaces $K_{\rho}^k$ and $Y_{\rho}^{\ell}$.

\section*{Appendix B}
\hspace*{\parindent} In this appendix we prove a Lemma which exemplifies the fact that for $\rho >
0$ and $\nu < 1$, the lower condition $\ell + 2 > n/2$ on $\ell$ needed in Lemma 3.3 in order to
make $Y_{\rho}^{\ell}$ into an algebra can be relaxed by using the estimate (3.4) instead of (3.3)
(but not uniformly in $\nu$ and $\rho$). Following Remark 3.1, we consider the Hilbert space $K$
with norm
$$\parallel u; K \parallel \ = \ \parallel \bar{f} \widehat{u}\parallel_2$$
\noi where $\bar{f}$ is either $ff_1$ or $f_0f_1$ with $f_0$, $f$ defined by (3.1) and $f_1$ defined
by (3.10) for some $k_<$, $k_> \in {I \hskip - 1 truemm R}^+$. \\

\noi {\bf Lemma B.1}. {\it Let $K$ be as above, with $0 < \nu < 1$, $\rho > 0$, $k_> \geq 0$ and $0
\leq k_< < n/2$. Then $K$ is an algebra, namely there exists a constant $C$ such that for all $u_1$,
$u_2 \in K$
$$\parallel u_1 \ u_2;K \parallel \ \leq C\parallel u_1;K\parallel \ \parallel u_2;K\parallel \ .
\eqno({\rm B.1})$$
\noi One can take
$$C^2 = \int d \eta \ \bar{f}(\eta )^{-2} \left ( 1 + 2^{2k} \ f_0(\eta )^{2\nu} \right )
\eqno({\rm B.2})$$
\noi where $k = k_< \vee k_>$ and the integral converges under the assumptions made on $\nu$,
$\rho$ and $k_<$.} \\

\noi {\bf Proof.} By the Schwarz inequality, we estimate
$$\parallel \bar{f} \widehat{u_1u_2} \parallel_2^2 \ = \int d\xi \ \bar{f}(\xi )^2 \left | \int d
\eta \  \widehat{u}_1(\eta ) \ \widehat{u}_2(\xi - \eta ) \right |^2
\leq \int d \xi \ \bar{f}(\xi)^2 \left \{ \int d \eta \ \bar{f}(\eta )^{-2} \ \bar{f} (\xi - \eta
)^{-2} \right \}$$
$$ \times \left \{ \int d \eta \ \bar{f} (\eta )^2 \ \bar{f}(\xi - \eta )^2 \ |\widehat{u}_1(\eta )
|^2 \ |\widehat{u}_2 (\xi - \eta )|^2 \right \}
                  \leq C^2 \parallel \bar{f} \widehat{u}_1 \parallel_2^2 \ \parallel \bar{f}
\widehat{u}_2 \parallel_2^2$$ 
\noi where
\begin{eqnarray*}
C^2 &=& \mathrel{\mathop {\rm Sup}_{\xi}} \bar{f}(\xi )^2 \int d\eta \ \bar{f}(\eta )^{-2} \
\bar{f}(\xi - \eta )^{-2}\\
&=& 2 \mathrel{\mathop {\rm Sup}_{\xi}} \bar{f}(\xi )^2 \int_{|\eta | \leq |\xi - \eta |} \bar{f}
(\eta )^{-2} \ \bar{f}(\xi - \eta )^{-2} \ .
\end{eqnarray*}
\noi Now for $|\xi | \leq |\xi - \eta|$, $\bar{f}(\xi) \leq \bar{f}(\xi - \eta )$ while for
$|\eta | \leq |\xi - \eta | \leq |\xi |$ one has $|\xi| \leq 2 |\xi - \eta |$ so that
$$f_1(\xi ) \leq 2^{k_< \vee k_>} \ f_1(\xi - \eta ) = 2^k \ f_1 (\xi - \eta ) \ ,$$
\noi and
$$f_{(0)}(\xi ) \leq f_{(0)}(\xi - \eta ) \ f_0(\eta )^{\nu}$$
\noi by (3.4). Therefore
$$C^2 \leq 2 \mathrel{\mathop {\rm Sup}_{\xi}} \left \{ \int_{|\eta | \vee |\xi|\leq |\xi - \eta | }
d\eta \ \bar{f}(\eta)^{-2} + \int_{|\eta | \leq |\xi - \eta |\leq |\xi |} d\eta \ \bar{f}
(\eta)^{-2} \ 2^{2k} \ f_0(\eta )^{2 \nu} \right \} \ .$$
\noi The Supremum over $\xi$ is easily seen to be the limit $|\xi| \to \infty$, namely
\begin{eqnarray*}
C^2 &\leq & 2 \left \{ \int_{\xi \cdot \eta \leq 0} d \eta \ \bar{f}(\eta)^{-2} + \int_{\xi \cdot
\eta \geq 0} d \eta \ \bar{f}(\eta)^{-2} \ 2^{2k} \ f_0 (\eta )^{2 \nu} \right \} \\ && \\
&= & \int d \eta \ \bar{f}(\eta)^{-2} \left ( 1 + 2^{2k} \ f_0(\eta )^{2 \nu} \right )
\end{eqnarray*} \noi and the integral converges for small $\eta$ by the condition $k_< < n/2$ and for
large $\eta$ by the conditions $\rho > 0$ and $\nu < 1$. \par \nobreak
\hfill $\sq$\par

The same result with essentially the same proof holds if one replaces $f_0$ by $\widetilde{f}$ in
the definition of $K$. One then has to replace (3.4) by (A.10) in the proof. The same result also
holds for arbitrary $k_> \in {I \hskip - 1 truemm R}$ (still with $0 \leq k_< < n/2)$, but the proof
is more cumbersome for $k_> < 0$.

\end{document}